\documentclass[11pt,reqno]{amsart}      
\usepackage[letterpaper, margin=1.2in]{geometry}
\usepackage{mathtools}
\usepackage[dvipsnames]{xcolor}
\usepackage{hyperref}
\usepackage{dsfont}
\usepackage{enumerate}
\usepackage{dsfont}
\usepackage[pdftex]{graphicx}
\usepackage{caption}
\usepackage{subcaption}
\usepackage{tikz}
\usetikzlibrary{positioning,decorations.pathreplacing}
\allowdisplaybreaks

\usepackage{amsmath, amsthm, amssymb}
\usepackage{latexsym,graphicx,multicol,mathrsfs,xypic}

\newtheorem{theorem}{Theorem}[section]

\newtheorem{proposition}[theorem]{Proposition}
\newtheorem{lemma}[theorem]{Lemma}
\newtheorem{corollary}[theorem]{Corollary}

\newtheorem{remark}[theorem]{Remark}
\newtheorem{assumption}[theorem]{Assumption}

\newtheorem*{proposition*}{Proposition}

\DeclareMathOperator{\var}{var}
\DeclareMathOperator{\cov}{cov}

\usepackage{autonum}

\thanks{This research has been partially supported by NSF grant DMS-1407558.}

\begin{document}

\title[Strict convexity of the free energy]{Strict convexity of the free energy of the canonical ensemble under decay of correlations}

\author{Younghak Kwon}
\address{Department of Mathematics, University of California, Los Angeles}
\email{yhkwon@math.ucla.edu}

\author{Georg Menz}
\address{Department of Mathematics, University of California, Los Angeles}
\email{gmenz@math.ucla.edu}

\subjclass[2010]{Primary: 82B05, Secondary: 60F05, 82B26.}
\keywords{Canonical ensemble, equivalence of ensembles, local Cram\'er theorem, strong interaction, phase transition}

\date{\today}

\begin{abstract}
We consider a one-dimensional lattice system of unbounded, real-valued spins. We allow arbitrary strong, attractive, nearest-neighbor interaction. We show that the free energy of the canonical ensemble converges uniformly in~$C^2$ to the free energy of the grand canonical ensemble. The error estimates are quantitative. A direct consequence is that the free energy of the canonical ensemble is uniformly strictly convex for large systems. Another consequence is a quantitative local Cram\'er theorem which yields the strict convexity of the coarse-grained Hamiltonian. With small adaptations, the argument could be generalized to systems with finite-range interaction on a graph, as long as the degree of the graph is uniformly bounded and the associated grand canonical ensemble has uniform decay of correlations. 
\end{abstract}

\maketitle

\section{Introduction} \label{intro} 

The broader scope of this article is the study of phase transitions. A phase transition occurs if a microscopic change in a parameter leads to a fundamental change in one or more properties of the underlying physical system. The most well-known phase transition is when water becomes ice. Many physical and non-physical systems and mathematical models have phase transitions. For example, liquid-to-gas phase transitions are known as vaporization. Solid-to-liquid phase transitions are known as melting. Solid-to-gas phase transitions are known as sublimation. More examples are the phase transition in the 2-d Ising model (see for example~\cite{Sel16}), the Erd\"os-Renyi phase transition in random graphs (see for example~\cite{ErRe60},~\cite{ErRe61} or~\cite{KrSu13}) or phase transitions in social networks (see for example~\cite{Fro07}). \\

We are interested in studying a one-dimensional lattice systems of unbounded real-valued spins. The system consists of a finite number of sites~$i \in \Lambda \subset \mathbb{Z}$ on the lattice~$\mathbb{Z}$. For convenience, we assume that the set~$\Lambda$ is given by~$\left\{1, \ldots, K \right\}$. At each site~$i \in \Lambda$ there is a spin~$x_i$. In the Ising model the spins can take on the value~$0$ or~$1$. In this article, we consider real-valued spins~$x_i \in \mathbb{R}$. A configuration of the lattice system is given by a vector~$x \in \mathbb{R}^{K}$. The energy of a configuration~$x$ is given by the Hamiltonian~$H: \mathbb{R}^K \to \mathbb{R}$ of the system. For the detailed definition of the Hamiltonian~$H$ we refer to Section~\ref{s_setting_and_main_results}. We consider arbitrary strong, attractive, nearest-neighbor interaction.\\

We consider two ensembles of the lattice system. The first ensemble is the grand-canonical ensemble which is given by the Gibbs measure
\begin{align}\label{d_gc_ensemble}
  \mu^\sigma(dx) = \frac{1}{Z} \exp\left( \sigma \sum_{i=1}^K x_i- H(x) \right) dx.
\end{align}
Here,~$Z$ is a generic normalization constant making the measure~$\mu^\sigma$ a probability measure. The constant~$\sigma \in \mathbb{R}$ is interpreted as an external field. The second ensemble is the canonical ensemble. It emerges from the grand-canonical ensemble by conditioning on the mean spin
\begin{align}
  m = \frac{1}{K} \sum_{i =1}^K x_i.
\end{align}
The canonical ensemble is given by the probability measure
\begin{align}
  \mu_m(dx)&= \mu^{\sigma} \left( dx \ | \    \frac{1}{K} \sum_{i=1}^K x_i= m  \right) \\
  &= \frac{1}{Z} \  \mathds{1}_{\left\{ \frac{1}{K}\sum_{i=1}^K x_i= m \right\}  } \ \exp(\sigma \sum_{i=1}^K x_i - H(x)) \mathcal{L}^{K-1}(dx),
\end{align} 
where~$\mathcal{L}^{K-1}$ denotes the~$(K-1)$-dimensional Hausdorff measure.\\

The grand-canonical ensemble has a phase transition on the two-dimensional lattice (see for example~\cite{Pei36}). However, on the one-dimensional lattice the grand-canonical ensemble does not have a phase transition if the interaction decays fast enough (see for example~\cite{Isi25,Dob68,Dob74,Rue68,MeNi14}). More precisely, in this work a system has no phase transition if the infinite-volume Gibbs measure of the system is unique. It is a natural question if the canonical ensemble~$\mu_m$ also does not have a phase transition on the one-dimensional lattice. This is a non-trivial question since there are known examples where the grand canonical ensemble has no phase transition but the canonical ensemble has (see for example~\cite{ScSh96,BiChKo02,BiChKo03}).\\

If the spins are $\left\{0,1 \right\}$-valued there is no phase transition for the canonical ensemble on a one dimensional lattice with nearest-neighbor interaction. The authors could not find a proof of that statement in the literature but it follows from a result by Cancrini, Martinelli and Roberto~\cite{CMR02}. There, a logarithmic Sobolev inequality is deduced for the canonical ensemble on lattices of arbitrary dimension, provided the grand canonical ensemble satisfies a mixing condition. The mixing condition used in~\cite{CMR02} is that the grand canonical ensemble has an exponential decay of correlation that is uniform in the external field~$\sigma$. This hypothesis is satisfied if the underlying lattice is one-dimensional. In our article we will use a similar mixing condition. \\

Up to the authors knowledge, this question is still open if the spins are real-valued and unbounded. We conjecture that this is true i.e.~the infinite-volume Gibbs measure of the canonical ensemble should be unique. A first step toward verifying this conjecture is to study the equivalence of the grand-canonical and canonical ensemble.
In equivalent ensembles, properties usually transfer from one ensemble to the other. The equivalence of ensembles in one-dimensional lattice system was deduced by Dobrushin~\cite{DoTi77} for discrete (or bounded) spin values or by Georgii~\cite{Ge95} for quadratic Hamiltonians. However, our case where the spin values are unbounded and the Hamiltonian is not quadratic is still open.\\

There are many different notions of equivalence of ensembles. We only consider the most simple type, namely the equivalence of thermodynamic quantities (see for example~\cite{Ad06}). This means that as the system size goes to infinity the free energy of the grand canonical ensemble converges to the free energy of the canonical ensemble (for more details see Section~\ref{s_setting_and_main_results} below). \\

In the main result of this article, i.e.~in Theorem~\ref{p_equivalence_of_ensembles} below, we show that the grand canonical and canonical ensemble are equivalent. In fact, we show that free energies converge uniformly in~$C^2$ as the system size goes to infinity. The rate of convergence in Theorem~\ref{p_equivalence_of_ensembles} is explicit. We therefore extend and refine the results of Dobrushin~\cite{DoTi77} and Georgii~\cite{Ge95}.\\

Our argument is quite general and should apply to more general situations. The argument does not use that the lattice is one-dimensional. Instead, it only uses that the grand canonical ensemble on a one-dimensional lattice has an uniform exponential decay of correlations (see for example~\cite{MeNi14} and~\cite{Zeg96}). Under the assumption of an uniform exponential decay of correlation, one should be able to use similar calculations to deduce the local Cram\'er theorem for spin systems on arbitrary graphs, as long as the degree is uniformly bounded and the interaction has finite range. However, we only consider the one-dimensional lattice with nearest-neighbor interaction because less notational burden is better for explaining ideas and presenting the calculations.\\

A consequence of Theorem~\ref{p_equivalence_of_ensembles} is that the free energy of the canonical ensemble is uniformly strictly convex and quadratic for large enough systems (see Corollary~\ref{p_strict_convexity_A_ce}). Strict convexity of the free energy rules out phase coexistence which corresponds to flat parts in the free energy. The most prominent example of phase coexistence is that under ordinary pressure water and ice can coexist at 0 degree Celsius. We want to point out that our result already applies to large but finite systems. In the infinite-volume limit, ordinary equivalence of ensembles (and not equivalence in $C^2$) would suffice to conclude that the free energy of the canonical ensemble is strictly convex.\\

Closely related to the free energy~$A_{ce}$ of the canonical ensemble is the notion of the coarse-grained Hamiltonian $\bar H$ (cf.~is \eqref{e_relation_A_c_bar_H} and~\cite{GrOtViWe09}). As in~\cite{GrOtViWe09}, we derive from Theorem~\ref{p_equivalence_of_ensembles} a local Cram\'er theorem (see Theorem~\ref{p_local_cramer}). The local Cram\'er theorem shows that the coarse-grained Hamiltonian converges in~$C^2$ to the Legendre transform of the free energy of the grand canonical ensemble. It is a direct consequence of the~$C^2$-local Cram\'er theorem that the coarse-grained Hamiltonian~$\bar H$ is also uniformly strictly convex for large enough system size~$|\Lambda|$ (cf.~Corollary~\ref{p_uniform_convexity_of_coarse_grained_Hamiltonian}).\\

The coarse-grained Hamiltonian~$\bar H$ plays an important role when studying the Kawasaki dynamics. The Kawasaki dynamics is natural drift diffusion process on our lattice system that conserves the mean spin of the system. The canonical ensemble is the stationary and ergodic distribution of the Kawasaki dynamics. The strict convexity of~$\bar H$ is a central ingredient for deducing a uniform logarithmic Sobolev inequality (LSI) for the canonical ensemble via the two-scale approach~\cite{GrOtViWe09}.
The LSI characterizes the speed of convergence of the Kawasaki dynamics to the canonical ensemble. With the equivalence of dynamic and static phase transitions (see~\cite{MeNi14} or~\cite{Yos03} a uniform LSI would also yield the absence of a phase transition and verify our conjecture  (i.e.~that the infinite-volume Gibbs measure is unique). Additionally, a uniform LSI is one of the main ingredients when deducing hydrodynamic limit of the Kawasaki dynamic via the two scale approach (see next paragraph). The uniform LSI for the canonical ensemble with no interaction is a well-known result (see for example~~\cite{Cha03,LaPaYa02,GrOtViWe09}). For weak interaction the uniform LSI was deduced in~\cite{Me11}. The question if the canonical ensemble satisfies a uniform LSI for strong nearest-neighbor interaction is still open. For $\left\{0,1 \right\}$-valued spins the answer is yes (see~\cite{CMR02}). The authors believe that this should also be the case for unbounded real-valued spins. \\

The strict convexity of the coarse-grained Hamiltonian also plays a crucial role when deducing the hydrodynamic limit of the Kawasaki dynamic. The hydrodynamic limit is a law of large numbers for processes. It states that under the correct scaling the Kawasaki dynamics (which is a stochastic process) converges to the solution of a non-linear heat equation (which is deterministic). It is conjectured by H.T.~Yau that the hydrodynamic limit also holds for strong finite-range interactions on a one-dimensional lattice. So far, this conjecture is still wide open. The strict convexity of the coarse-grained Hamiltonian, which is deduced in this article, is an important cornerstone to tackle this problem with the help of the two-scale approach (see~\cite{GrOtViWe09}).\\

Let us now comment on how the $C^2$-equivalence of ensembles is deduced. The motivation for our approach comes from the proof of the local Cram\'er theorem in~\cite{GrOtViWe09} and~\cite{Me11}. By using Cram\'er's trick of an exponential shift is suffices to show~$C^2$-bounds on the density of a sum of random variables ~$X_i$ (see also Proposition~\ref{p_main computation} below). Those desired bounds were derived ~\cite{GrOtViWe09} and~\cite{Me11} via a local central limit theorem (clt) for independent random variables. Our situation is a lot more subtle: Instead of deducing a local clt for independent random variables we would have to deduce a local clt for dependent random variables. At this point one could hope to use existing methods to deduce the local clt. Let us mention for example the approach of Dobrushin~\cite{Dob74}, the approach of Bender~\cite{Ben73} or the approach of Wang and Woodroofe~\cite{WaWo90}. Unfortunately this does not help. All methods --at least the ones that are known to the authors-- use the following principle (see also~\cite{DeSaMe16}):
\begin{align}
  \mbox{integral clt} \quad + \quad  \mbox{regularity} \ \Rightarrow \ \mbox{local clt}.
\end{align}
The first ingredient, namely the integral clt for the dependent random variables~$X_i$ is relatively easy to deduce. There are a lot of methods available. Let us mention for example Stein's method (see for example~\cite{ChGoSh11}), methods that are based on mixing, or methods that are based on Donsker's theorem (see for example~\cite{Dur2010}). Deducing the second ingredient is tricky, not to mention that Dobrushin~\cite{Dob74} carried out that step only for discrete or bounded random variables.\\

All in all, this approach has two fundamental problems. The first one is that we need not only to control the density itself but also the first and second derivative. As a consequence, one would need very detailed information about the regularity of the density. We also believe that showing this regularity is as hard as directly deducing the local central limit theorem. Let us turn to the second problem. In order to deduce Theorem~\ref{p_equivalence_of_ensembles} the local central limit theorem must be quantitative. Using the principle from above yields suboptimal rates of convergence. For deducing Theorem~\ref{p_equivalence_of_ensembles} one has to iteratively apply the principle three times; and in each iteration the convergence rate gets worse. One would have to hope that in the end the convergence rate is still good enough for deducing Theorem~\ref{p_equivalence_of_ensembles}.\\

Instead of using the principle from above, we generalize a well-known method for proving the local clt for independent random variables to dependent ones. We generalize the method that is based on characteristic functions and Fourier inversion (see~\cite{Fel71} and~\cite{GrOtViWe09}). Calculations get quite evolved and lengthy. We do not deduce a local clt for dependent random variables in this work. Instead, we only deduce bounds that are needed to deduce Theorem~\ref{p_equivalence_of_ensembles} (cf.~Proposition~\ref{p_main computation} below). However, one could use our calculations as a guideline for deducing a quantitative, local clt for dependent random variables. When doing so, one would have to substitute some of our arguments that use the specific structure of our lattice model. We use the following special structure:
\begin{itemize}
\item Exponential decay of correlations (see~Lemma~\ref{l_exponential decay of correlations}).\\[-2ex]

\item The interaction has finite range~$R$. More precisely, we use that two spins ~$x_i$ and~$x_j$ become independent if~$|i-j|>R$ and one conditions on the spin values~$(x_k)_{i< k <j}$ between them (see Section~\ref{s_setting_and_main_results}).\\[-2ex]

\item The Hamiltonian is quadratic. More precisely, we use the following consequence. For all~$i \in \Lambda$ the conditional variances~$\var(X_i|X_j, \ |j-i|\leq l )$ is bounded from above and below uniformly in the values~$X_j, \ |j-i|\leq l$ (see Section~\ref{s_setting_and_main_results}).\\[-2ex]

\item Higher moments of~$X_i$ conditioned on~$X_j$,~$ j \neq i$ are uniformly controlled by lower moments. This fact is used to show that the characteristic functions of~$X_i$ conditioned on~$X_j$,~$ j \neq i$ have a uniform decay (see Lemma~\ref{l_variance estimates} and Lemma~\ref{c_amgm}).
\end{itemize}

As mentioned before, the~$C^2-$local Cram\'er theorem (see Theorem~\ref{p_local_cramer}) is deduced by generalizing the argument of~\cite{GrOtViWe09} for independent random variables to dependent random variables. This adds a lot more complexity to the task. We overcome the technical challenges of considering dependent random variables by using two strategies. The first strategy is to induce artificial independence by conditioning on even or odd random variables. The second strategy is to handle dependencies as a perturbation. We morally treat large blocks as single sites of a coarse-grained system. Because there is a big distance between the blocks, the blocks are only weakly dependent. Then, the error term can be controlled by using the decay of correlations. For more details we refer to the comments after Proposition~\ref{p_main computation} and at the beginning of Section~\ref{s_main_computation}.\\

Let us shortly discuss possible generalizations of our main result. We expect that one can generalize our method with only slight modifications to the following situation:
\begin{itemize}
  \item instead of nearest-neighbor interaction to finite range interaction.\\[-2ex]
  
  \item instead of exponential decay of correlations to sufficiently fast algebraic decay.\\[-2ex]
  
  \item instead of a 1d lattice to any lattice or graph with bounded degree, as long as the grand canonical ensemble~$\mu^\sigma$ has sufficient decay of correlations, uniformly in the system size and the external field~$\sigma$.\\[-2ex]
  
  \item instead of attractive interaction to repulsive and mixed interactions, as long as the estimate
  \begin{align}
   \var_{\mu^\sigma} \left[ \sum_{i=1}^{K} X_i \right] \geq C K 
\end{align}
is satisfied. For attractive interaction this estimate is deduced in Lemma~\ref{l_variance estimates}.
\end{itemize}

More challenging, it would be very interesting to study the local Cram\'er theorem for the following changes:
\begin{itemize}   
\item Instead of finite-range interaction one could consider infinite-range interaction. In the case of the grand canonical ensemble on a one-dimensional lattice, there is no phase transition if the interaction~$M_{ij}$ decays algebraically faster than~$(1+ |i-j|)^{2+\varepsilon }$. The decay condition is sharp (see for example~\cite{MeNi14} and references therein). It would be very interesting to know if the~$C^2$-local Cram\'er theorem (see Theorem~\ref{p_local_cramer}) also holds for this decay or a stronger algebraic decay is needed.\\[-2ex]

\item In our model we need a quadratic single-site potential. Inspired from~\cite{FaMe14}, it is natural to ask if the local Cram\'er theorem also holds for super-quadratic or polynomially increasing single-site potentials.\\[-2ex]

\item Inspired by~\cite{Dob74} or~\cite{Ge95} it would be interesting to study more general interaction than pairwise-quadratic interaction.
\end{itemize}

We conclude the introduction by giving a short overview over the remaining article. In Section~\ref{s_setting_and_main_results} we introduce the precise setting and formulate the main results. In Section~\ref{s_proof_local_cramer} we deduce the main results of this article up to Proposition~\ref{p_main computation}. The main computations are done in Section~\ref{s_main_computation} where we give the proof of Proposition~\ref{p_main computation}.

\section*{Conventions and Notation}

\begin{itemize}
\item The symbol~$T_{(k)}$ denotes the term that is given by the line~$(k)$.
\item With uniform we mean that a statement holds uniform in the system size~$\Lambda$,  the mean spin~$m$ and the external field~$s$.
\item We denote with~$0<C<\infty$ a generic uniform constant. This means that the actual value of~$C$ might change from line to line or even within a line. 
\item $C(n)$ denotes a constant that only depends on~$n$.
\item $a \lesssim b$ denotes that there is a uniform constant~$C$ such that~$a \leq C b$.
\item $a \sim b$ means that~$a \lesssim b$ and~$b \lesssim a$.
\item $\mathcal{L}^{k}$ denotes the~$k$-dimensional Hausdorff measure.
\item $Z$ is a generic normalization constant. It denotes the partition function of a measure.  
\item For a function~$f : \mathbb{R}^K \rightarrow \mathbb{R}$,~$\text{supp} \ f  = \{i_1, i_2, \cdots, i_k \}$ denotes the minimal subset of~$\{1, 2, \cdots K\}$ such that~$f(x) = f(x_{i_1}, \cdots, x_{i_k} )$.
\end{itemize}

\section{Setting and main results}
\label{s_setting_and_main_results}
We start with explaining the details of our model. We consider the sublattice~$\left\{1, \ldots, K \right\} \subset \mathbb{Z}$. The Hamiltonian~$H:\mathbb{R}^K \to \mathbb{R}$ of the system is defined as
\begin{align}\label{e_d_hamiltonian}
H(x) = \sum_{i=1}^K \left(\psi (x_i) + s_i x_i - J  x_i x_{i+1} \right),
\end{align}
where~$x_{K+1}$ is defined to be~$0$. We make the following assumptions:
\begin{itemize}
\item The single-site potential~$\psi: \mathbb{R} \to \mathbb{R}$ can be written as
\begin{align}\label{e_single_site_potential}
 \psi (z) = \frac{1}{2} z^2 + \psi_{b}(z),
 \end{align}
where the function~$\psi_b: \mathbb{R} \to \mathbb{R}$ satisfies
 \begin{align}\label{e_nonconvexity_bounds_on_perturbation}
 |\psi_b|_{\infty} + |\psi'_b|_{\infty}  + |\psi''_b|_{\infty} < \infty. 
 \end{align}
 
\item The numbers~$s=(s_i) \in \mathbb{R}^K$ are arbitrary. They model the interaction of the system with an external field or the boundary.\\[-2ex]

\item The number~$J \in \left(-\frac{1}{4}, \frac{1}{4}\right)$ is arbitrary. It models the strength of the interaction. The interaction is attractive if~$J>0$. The interaction is repulsive if~$J<0$ . 
\end{itemize}

Now, let us turn to the first main result of this article, namely the equivalence of ensembles (see Theorem~\ref{p_equivalence_of_ensembles} from below). The grand canonical ensemble (gce)~$\mu^\sigma$ is a probability measure on~$\mathbb{R}^K$ given by the Lebesgue density
\begin{align} \label{e_mu sigma def}
\mu^{\sigma}(dx) : = \frac{1}{Z} \exp\left( \sum_{i=1}^K \sigma x_i -H(x) \right) dx.
\end{align}
The free energy of the gce~$\mu^\sigma$ is given by (cf.~\eqref{d_gc_ensemble} and~\cite{GrOtViWe09})
\begin{align}
 A_{gce}(\sigma):= \widehat{ \mathcal{H}}_K := \frac{1}{K} \ln \int_{\mathbb{R}^K} \exp\left( \sigma \sum_{i=1}^K x_i - H(x)\right) dx.
\end{align}
We observe that~$A_{gce}$ is uniformly strictly convex. More precisely, it holds:
\begin{lemma}\label{p_convexity_free_energy_gce}
Let~$\left(X_1, X_2, \cdots, X_K \right)$ be a real-valued random variable distributed according to the gce~$\mu^\sigma$. Assume that
\begin{align}\label{e_var gtr K}
\var\left(\sum_{i=1}^K X_i \right) \gtrsim K.
\end{align} Then the free energy~$A_{gce}$ of the gce~$\mu^\sigma$ is uniformly strictly convex in the sense that there exists a constant~$C>0$ such that for all~$\sigma \in \mathbb{R}$
\begin{align}
  \frac{1}{C} \leq \frac{d^2}{d \sigma^2} A_{gce} (\sigma) \leq C.
\end{align}
\end{lemma}

The proof of Lemma~\ref{p_convexity_free_energy_gce} is given in Section~\ref{s_proof_local_cramer}. The core ingredient of the argument is a uniform Poincar\'e inequality. The additional assumption~\eqref{e_var gtr K} is not very restrictive. For example, it is automatically satisfied if the interaction is attractive (see Lemma~\ref{l_variance estimates} below). \\

Let us turn to the canonical ensemble (ce)~$\mu_m$. It emerges from the gce by conditioning (i.e.~fixing) on the mean spin
\begin{align}
  m = \frac{1}{K} \sum_{i =1}^K x_i.
\end{align}
The ce is given by the probability measure
\begin{align} \label{d_c_ensemble}
  \mu_m(dx)&= \mu^{\sigma} \left( dx \ | \    \frac{1}{K} \sum_{i=1}^K x_i= m  \right) \\
  &= \frac{1}{Z} \  \mathds{1}_{\left\{ \frac{1}{K}\sum_{i=1}^K x_i= m \right\}  } \ \exp(\sigma \sum_{i=1}^K x_i - H(x)) \mathcal{L}^{K-1}(dx),
\end{align} 
where~$\mathcal{L}^{K-1}$ denotes the~$(K-1)$-dimensional Hausdorff measure. The free energy of the ce~$\mu_{m}$ is given by 
\begin{align} \label{e_free_energy_ce}
 A_{ce}(\sigma)= \frac{1}{K} \ln \int_{\left\{\frac{1}{K} \sum_{i=1}^K x_i =m  \right\}} \exp\left( \sigma \sum_{i=1}^K x_i - H(x)\right) \mathcal{L}^{K-1}(dx).
\end{align}

Equivalence of ensembles only holds if the external field~$\sigma$ of the gce~$\mu^\sigma$ and the mean spin~$m$ of the ce~$\mu_m$ are related in the following way.

\begin{assumption}
We then choose~$\sigma=\sigma(m) \in \mathbb{R}$ and~$m= m (\sigma) \in \mathbb{R}$ such that the following relation is satisfied:
\begin{align}\label{e_relation_m_sigma}
  \frac{d}{d\sigma} A_{gce} (\sigma) =m.
\end{align}
By the strict convexity of~$A_{gce}$ (see Lemma~\ref{p_convexity_free_energy_gce}) there exists for any~$m \in \mathbb{R}$ a unique~$\sigma= \sigma(m) \in \mathbb{R}$ that satisfies the relation~\eqref{e_relation_m_sigma} or vice versa.
\end{assumption}

Now, let us formulate our first main result, namely the equivalence of the free energies in~$C^2$.
\begin{theorem}[Equivalence of ensembles]\label{p_equivalence_of_ensembles}
Let~$\left(X_1, X_2, \cdots, X_K \right)$ be a real-valued random variables distributed according to
\begin{align}
\mu^{\sigma}(dx) : = \frac{1}{Z} \exp\left( \sum_{i=1}^K \sigma x_i -H(x) \right) dx.
\end{align}
Assume that
\begin{align}
\var\left(\sum_{i=1}^K X_i \right) \gtrsim K.
\end{align}
Then it holds that 
\begin{align}\label{e_local_cramer}
 \lim_{K \to \infty} \left|   A_{gce}  - A_{ce} \right|_{C^2} = 0,
\end{align}
where the convergence is uniform in the mean spin~$m$ and the external field~$s$. More precisely, given a constant~$\varepsilon>0$, there is an integer~$K_0\in \mathbb{N}$ such that for all~$K \geq K_0$
\begin{align}
\sup_{\sigma \in \mathbb{R}}\left|   A_{gce}(\sigma)  - A_{ce}(\sigma) \right| &\lesssim \frac{1}{K} ,\label{e_eoe_C_0}\\
\sup_{\sigma \in \mathbb{R}}\left| \frac{d}{d \sigma} A_{gce}(\sigma)  - \frac{d}{d \sigma} A_{ce}(\sigma) \right| &\lesssim \frac{1}{K^{1-\varepsilon}},\label{e_eoe_C_1} \\
\sup_{\sigma \in \mathbb{R}}\left| \frac{d^2}{d \sigma^2} A_{gce}(\sigma)  - \frac{d^2}{d \sigma^2}  A_{ce} (\sigma) \right| &\lesssim \frac{1}{K^{\frac{1}{2} - \varepsilon}} \label{e_eoe_C_2}.
\end{align}
\end{theorem}
We would like to emphasize that Theorem~\ref{p_equivalence_of_ensembles} contains explicit rates of convergence. We give the proof of Theorem~\ref{p_equivalence_of_ensembles} in Section~\ref{s_proof_local_cramer}. \\

A direct consequence of Lemma~\ref{p_convexity_free_energy_gce} and Theorem~\ref{p_equivalence_of_ensembles} is that the free energy~$A_{ce}$ is uniformly strictly convex for large enough systems.

\begin{corollary}\label{p_strict_convexity_A_ce}
There is a uniform constant~$0<C<\infty$ and an integer~$K_0\in \mathbb{N}$ such that for all~$K \geq K_0$ and all~$\sigma \in \mathbb{R}$
\begin{align}\label{e_strict_convexity_A_ce}
  \frac{1}{C} \leq \frac{d^2}{d \sigma^2} A_{ce} (\sigma) \leq C.
\end{align}
\end{corollary}

Let us turn to the second main result of this article, the local Cram\'er theorem. For that purpose let us introduce~$\mathcal{H}_K$ which denotes the Legendre transform of the free energy~$A_{gce}$ (also denoted by~$\widehat{ \mathcal{H}}_K$) i.e.
\begin{align}\label{e_d_cramer_transform}
\mathcal{H}_K (m) = \sup_{\sigma \in \mathbb{R}} \left( \sigma  m -\widehat{ \mathcal{H}}_K (\sigma) \right).
\end{align}
It follows from elementary observations that~$\mathcal{H}_K$ is uniformly strictly convex. 

\begin{lemma}\label{p_strict_convexity_Legendre_transform}
For any~$m \in \mathbb{R}$
\begin{align}
  \mathcal{H}_K (m) =  \sigma(m)   m -\widehat{ \mathcal{H}}_K (\sigma (m)).
\end{align}
Additionally, under the same assumptions as in Theorem~\ref{p_equivalence_of_ensembles}, it holds that~$\mathcal{H}_K$ is uniformly strictly convex in the sense that there is a uniform constant~$0<C< \infty$ such that for all~$\sigma \in \mathbb{R}$
\begin{align}
  \frac{1}{C} \leq \frac{d^2}{dm^2}\mathcal{H}_K (m) \leq C.
\end{align}
\end{lemma}
We give the proof of Lemma~\ref{p_strict_convexity_Legendre_transform} in Section~\ref{s_proof_local_cramer}.\\

The coarse-grained Hamiltonian~$\bar H: \mathbb{R} \to \mathbb{R}$ is defined as
\begin{align}\label{d_bar_H}
\bar{H}(m)= - \frac{1}{K} \ln \int_{\left\{ \frac{1}{K} \sum_{i=1}^K  x_i = m  \right\}} \exp (- H(x)) \mathcal{L}^{K-1} (dx).
\end{align}
Hence, we can rewrite the free energy of the ce as
\begin{align}\label{e_relation_A_c_bar_H}
  A_{ce}(\sigma)= \sigma m - \bar H(m).
\end{align}
It follows that the difference of the free energies~$A_{gce}$ and~$A_{ce}$ can be expressed as
\begin{align}
A_{gce}(\sigma) - A_{ce}(\sigma)& =  \hat{\mathcal{H}}_K (\sigma) - \sigma m + \bar H(m) \\
& = \bar{H}(m)-\mathcal{H}_K (m).   \label{e_rewriting_difference_free_energies}
\end{align}

From Theorem~\ref{p_equivalence_of_ensembles} we deduce the following local Cram\'er theorem.

\begin{theorem}[~$C^2$-local Cram\'er theorem]\label{p_local_cramer}
Let~$\left(X_1, X_2, \cdots, X_K \right)$ be a real-valued random variables distributed according to
\begin{align}
\mu^{\sigma}(dx) : = \frac{1}{Z} \exp\left( \sum_{i=1}^K \sigma x_i -H(x) \right) dx.
\end{align}
Assume that
\begin{align}
\var\left(\sum_{i=1}^K X_i \right) \gtrsim K.
\end{align}
Then it holds that 
\begin{align}
 \lim_{K \to \infty} \left|  \bar H (m) - \mathcal{H}_K (m) \right|_{C^2} = 0,
\end{align}
where the convergence is uniform in the mean spin~$m$ and the external field~$s$. More precisely, given a constant~$\varepsilon>0$, there is an integer~$K_0\in \mathbb{N}$ such that for all~$K \geq K_0$
\begin{align}
\sup_{m \in \mathbb{R}}\left|  \bar H (m) - \mathcal{H}_K (m) \right| &\lesssim \frac{1}{K} , \label{e_theorem1 1}\\
\sup_{m \in \mathbb{R}}\left| \frac{d}{d m} \bar H (m) - \frac{d}{d m}  \mathcal{H}_K (m) \right| &\lesssim \frac{1}{K^{1-\varepsilon}} ,\label{e_theorem1 2} \\
\sup_{m \in \mathbb{R}}\left| \frac{d^2}{d m^2} \bar H (m) - \frac{d^2}{d m^2} \mathcal{H}_M (m) \right| &\lesssim \frac{1}{K^{\frac{1}{2} - \varepsilon}}. \label{e_theorem1 3}
\end{align}
\end{theorem}
Theorem~\ref{p_local_cramer} is an extension of the 
local Cram\'er theorems that were deduced in Proposition~31 in~\cite{GrOtViWe09}, Theorem~4 in~\cite{Me11} and~\cite{MeOt13}. The proof of Theorem~\ref{p_local_cramer} is stated in Section~\ref{s_proof_local_cramer}. The main ingredient is Theorem~\ref{p_equivalence_of_ensembles}.\\

An important consequence of Lemma~\ref{p_strict_convexity_Legendre_transform} and of Theorem~\ref{p_local_cramer} is that for large enough systems the coarse-grained Hamiltonian~$\bar{H}$ is uniformly strictly convex.
\begin{corollary}\label{p_uniform_convexity_of_coarse_grained_Hamiltonian} Under the assumptions of Theorem~\ref{p_local_cramer} there is an positive integer~$K_0$ such that for all~$K\geq K_0$ the coarse-grained Hamiltonian~$\bar H: \mathbb{R} \to \mathbb{R}$ is uniformly strictly convex. More precisely, there is a uniform constant~$0<C< \infty$ such that for all~$m \in \mathbb{R}$
\begin{align}
\frac{1}{C} \leq \frac{d^2}{d m^2} \bar H(m) \leq C.
\end{align}
\end{corollary}

\section{Proof of the main results}\label{s_proof_local_cramer}

In this section we prove the main results of this article. 
\begin{assumption}
From now on we assume that~$X=(X_1, X_2, \cdots , X_K)$ is a real-valued random vector distributed according to
\begin{align}
\mu^{\sigma}(dx) : = \frac{1}{Z}\exp\left(\sum_{i=1}^{K} \sigma x_i -H(x) \right)dx.
\end{align}
\end{assumption}
We begin with simple auxiliary lemma.
\begin{lemma} \label{l_variance estimates}
There exists a uniform constant~$C$ such that
\begin{align}
\var\left(  \sum_{i=1}^{K} X_i \right) \leq CK.
\end{align}
Moreover, if~$J$ in~\eqref{e_d_hamiltonian} is nonnegative, then the condition~\eqref{e_var gtr K} in Lemma~\ref{p_convexity_free_energy_gce} is satisfied.
\end{lemma}
The lemma from above shows that the variance of the the mean spin of the gce~$\mu^\sigma$ is well behaved. \medskip

\textsc{Proof of Lemma~\ref{l_variance estimates}}. \ Let us begin with the proof of the upper bounds. It is known that gce~$\mu^{\sigma}$ satisfies a uniform logarithmic Sobolev inequality (LSI) (see for example~\cite[Theorem 1.6]{MeNi14}). It is also well known that logarithmic Sobolev inequality implies Poincar\'e inequality (PI). Therefore it holds that
\begin{align}\label{e_upper_bound_variance_gce}
\var \left(\sum_{i=1}^K X_i \right) \leq \frac{1}{\rho} \int \left| \nabla \left(\sum_{i=1}^K X_i \right) \right|^2 d \mu^{\sigma} = \frac{1}{\rho} K.
\end{align}
where~$\rho>0$ is a constant in Poincar\'e inequality independent of~$m$ and~$K$, which proves the upper bound. \\

Proof of the lower bound relies on~\cite[Lemma~9]{Me11}. Let~$W$ be a random variable distributed according to the distribution
\begin{align}
\nu(dz)= \frac{1}{Z} \exp\left( -\psi(z) -tz\right)dz.
\end{align}
In~\cite[Lemma~9]{Me11}, it was shown that there is a constant~$0<C<\infty$ such that for all~$t \in \mathbb{R}$,
\begin{align} \label{e_me11_var_estimate}
\frac{1}{C} \leq \var_{\nu}\left(W\right) \leq C.
\end{align}
Observe that the conditional random variable~$ X_i \mid X_j : j \ne i  $ has the Lebesgue density
\begin{align}\label{e_conditional lebesgue density}
\mu^{\sigma}\left( dx_i | \bar{x_i} \right) = \frac{1}{Z} \exp\left( -\psi(x_i ) - \left(- Jx_{i-1} - Jx_{i+1} +s_i- \sigma\right) x_i \right)dx_i.
\end{align}
Let~$\bar{\mu}^{\sigma}(d\bar{x_i})$ be the marginal measure defined by the following property
\begin{align}
\mu^{\sigma}(dx)=\mu^{\sigma} \left( dx_i | \bar{x_i} \right) \bar{\mu }^{\sigma}(d\bar{x_i}).
\end{align}
Then it follows that
\begin{align}
\var_{\mu^{\sigma}}\left(X_i\right) & = \int \var_{\mu^{\sigma}(dx_i \mid \bar{x_i})}\left(X_i\right) \bar{\mu}^{\sigma}(d\bar{x_i}) + \var_{\bar{\mu}^{\sigma}(d\bar{x_i})}\left( \int x_i \mu^{\sigma}(dx_i \mid \bar{x_i }) \right) \\
& \geq \int \var_{\mu^{\sigma}(dx_i \mid \bar{x_i})}\left(X_i\right) \bar{\mu}^{\sigma}(d\bar{x_i}) \\
& \geq \int \frac{1}{C}\bar{\mu}^{\sigma}(d\bar{x_i}) =\frac{1}{C}. \label{pl32_1}
\end{align}
Moreover, Menz and Nittka (see~\cite[Lemma 2.1]{MeNi14} for example) proved that for~$J\geq0$ we have
\begin{align}
\cov\left(X_i, X_j \right) \geq 0. \label{pl32_2}
\end{align}
Then it follows using~\eqref{pl32_1} and~\eqref{pl32_2} that
\begin{align}
\var \left(\sum_{i=1}^{K} X_i \right) &= \sum_{i=1}^{K} \var\left(X_i \right) + \sum_{i \ne j} \cov\left(X_i, X_j \right)  \geq \frac{1}{C}K.
\end{align}
This finishes the proof of Lemma~\ref{l_variance estimates}.
\qed

\medskip

Now we are ready to give proofs of Lemma~\ref{p_convexity_free_energy_gce} and Lemma~\ref{p_strict_convexity_Legendre_transform}. \\

\textsc{Proof of Lemma~\ref{p_convexity_free_energy_gce}}. \ It is a direct consequence of Lemma~\ref{l_variance estimates}. Indeed, we have
\begin{align}
\frac{d}{d\sigma} A_{gce}(\sigma) &= \frac{d}{d\sigma}\left(\frac{1}{K} \ln \int_{\mathbb{R}^K} \exp\left( \sigma \sum_{i=1}^K x_i - H(x)\right) dx \right) \\
& = \frac{1}{K} \frac{ \int_{\mathbb{R}^K} \sum_{i=1}^{K} x_i \exp\left( \sigma \sum_{i=1}^K x_i - H(x)\right) dx }{\int_{\mathbb{R}^K}  \exp\left( \sigma \sum_{i=1}^K x_i - H(x)\right) dx } = \frac{1}{K} \mathbb{E} \left[ \sum_{i=1}^K X_i \right].
\end{align}
Taking the derivative with respect to $\sigma$ again, we obtain
\begin{align}
\frac{d^2}{d\sigma^2} A_{gce}(\sigma) & = \frac{1}{K} \int_{\mathbb{R}^K} \sum_{i=1}^K x_i \left(\sum_{j=1}^K \left(x_j -\mathbb{E} \left[ X_j \right] \right) \right) d\mu^{\sigma} (dx) \\
& = \frac{1}{K} \mathbb{E} \left[ \sum_{i=1}^K X_i \sum_{j=1}^K \left(X_j -\mathbb{E}\left[X_j\right] \right) \right] = \frac{1}{K} \var\left(\sum_{i=1}^{K} X_i \right).
\end{align}
Therefore, we conclude from Lemma~\ref{l_variance estimates} that there is a constant~$C>0$ with
\begin{align}
\frac{1}{C} \leq \frac{d^2}{d \sigma^2} A_{gce}(\sigma) \leq C.
\end{align}
\qed
\medskip

\textsc{Proof of Lemma~\ref{p_strict_convexity_Legendre_transform}}. \
Let us introduce an auxiliary notation~$m_i$ defined by
\begin{align} \label{d_def_mi}
m_i : = \mathbb{E} \left[ X_i \right] 
\end{align}
Since~$\mathcal{H}_K (m)$ is the Legendre transform of the strict convex function~$\widehat{ \mathcal{H}}_K (\sigma)$, there exists a unique~$\sigma=\sigma(m)$ such that
\begin{align}
\mathcal{H}_K (m) = \sigma(m)m - \widehat{ \mathcal{H}}_K (\sigma(m)).
\end{align}
Moreover, for each~$m$,~$\sigma(m)$ satisfies
\begin{align}
\frac{d}{d \sigma} \left( \sigma m - \widehat{\mathcal{H}}_K (\sigma) \right) = 0,
\end{align}
which is equivalent to
\begin{align}
m&= \frac{d}{d \sigma} \widehat{\mathcal{H}}_K (\sigma)  \\
&= \frac{1}{K}\frac{ \int_{\mathbb{R}^k} \left(\sum_{i=1}^{K} x_i\right) \exp \left( \sum_{i=1}^K \sigma   x_i - H(x) \right) dx}{\int_{\mathbb{R}^k} \exp \left( \sum_{i=1}^K \sigma   x_i - H(x) \right) dx} \\
& = \frac{1}{K} \mathbb{E} \left[ \sum_{i=1}^K X_i \right] \overset{\eqref{d_def_mi}}{=} \frac{1}{K} \sum_{i=1}^{K} m_i. \label{e_m=1/Kmi}
\end{align}
Then it follows that
\begin{align}
\frac{d}{dm} \mathcal{H}_K(m) &= \frac{d}{dm}  \left(\sigma(m)m - \widehat{ \mathcal{H}}_K (\sigma(m))\right) \\
&= \frac{d \sigma(m)}{dm}m + \sigma(m)- \frac{d}{d\sigma}\widehat{ \mathcal{H}}_K (\sigma) \cdot \frac{d\sigma}{dm}\\
&= \frac{d \sigma}{dm}m + \sigma- \frac{1}{K} \mathbb{E}\left[ \sum_{i=1}^K X_i \right]\cdot \frac{d\sigma}{dm} =  \sigma.
\end{align}
In Lemma~\ref{p_convexity_free_energy_gce} we proved
\begin{align} \label{e_dmdsigma}
\frac{d}{d\sigma} m = \frac{d}{d\sigma}\left(\frac{1}{K} \sum_{i=1}^K \mathbb{E} \left[ X_i \right] \right) = \frac{1}{K} \var\left(\sum_{i=1}^K X_i \right).
\end{align}
Thus Lemma~\ref{l_variance estimates} implies there exists a constant~$C>0$ with
\begin{align}
\frac{1}{C} \leq \frac{d^2}{d \sigma^2} \mathcal{H}_K (m) = \frac{d\sigma}{dm} = \left(\frac{dm}{d\sigma}\right)^{-1} \leq C. 
\end{align}
\qed

Let us now turn to the proof of Theorem~\ref{p_equivalence_of_ensembles}. We need some more auxiliary results. The first one is Cram\'er's trick of exponential shift of measures.

\begin{lemma}\label{l_cramers_trick}
It holds that
\begin{align}
g_{K, m} (0) = \exp\left( K A_{ce}(\sigma) - K  A_{gce} (\sigma) \right)\label{cramers representation1}\\
\overset{\eqref{e_rewriting_difference_free_energies}}{=} \exp\left( K \mathcal H_K (m)  - K\bar H (m)\right).\label{cramers representation2}
\end{align}
Here,~$g_{K, m} $ denotes the distribution of 
\begin{align}
\frac{1}{\sqrt{K}} \sum_{i=1}^K \left(X_i -m \right).
\end{align}
\end{lemma}

\textsc{Proof of Lemma~\ref{l_cramers_trick}}. \ The lemma follows from a direct computation:
\begin{align}
&K\mathcal H_K (m) - K \bar H (m) \\
&\qquad = K \sigma(m)m - \ln \int_{\mathbb{R}^K} \exp \left( \sum_{i=1}^K \sigma(m)   x_i - H(x) \right) dx \\
& \qquad \quad +\ln \int_{\left\{ \frac{1}{K} \sum_{i=1}^K  x_i = m  \right\}} \exp \left(- H(x)\right) \mathcal{L}^{K-1} (dx) \\
& \qquad = \ln \int_{\left\{ \frac{1}{K} \sum_{i=1}^K  x_i = m  \right\}} \exp \left(K \sigma(m)m- H(x)\right) \mathcal{L}^{K-1} (dx)\\
& \qquad  \quad -\ln \int_{\mathbb{R}^K} \exp \left( \sum_{i=1}^K \sigma(m)   x_i - H(x) \right) dx \\
& \qquad  = \ln \frac{ \int_{\left\{ \frac{1}{\sqrt{K}} \sum_{i=1}^K  \left(x_i -m\right)=0  \right\}} \exp \left( \sigma(m)\sum_{i=1}^K x_i - H(x)\right) \mathcal{L}^{K-1} (dx) }{\int_{\mathbb{R}^K} \exp \left( \sum_{i=1}^K \sigma(m)   x_i - H(x) \right) dx} \\
& \qquad  = \ln g_{K, m} (0).
\end{align}
Taking the exponential function and using~\eqref{e_rewriting_difference_free_energies} yields the lemma as desired.
\qed
\medskip

Next, we need the following direct consequence of Lemma~\ref{l_variance estimates}.
\begin{lemma}\label{c_amgm}
Assume that the single-site potential~$\psi$ satisfies~\eqref{e_single_site_potential} and~\eqref{e_nonconvexity_bounds_on_perturbation}. Recall the definition~\eqref{d_def_mi} of~$m_i$. Then for any finite set~$A_i \subset \{1, 2, \cdots, K\}$ and~$k \in \mathbb{N}$, we have
\begin{align} \label{e_amgm}
\left|\mathbb{E} \left[ \sum_{i_1 \in A_{1}}\cdots \sum_{i_k \in A_{k}} \left(X_{i_1}-m_{i_1} \right) \cdots \left(X_{i_k}- m_{i_k} \right) \right]\right| \lesssim \left|A_{1}\right|\cdots\left|A_k \right|,
\end{align}
where the constant only depends on~$k \in \mathbb{N}$.
\end{lemma}

\textsc{Proof of Lemma~\ref{c_amgm}}. \ As noted in the proof of Lemma~\ref{l_variance estimates}, we know that gce~$\mu^\sigma$ satisfies a uniform Poincar\'e inequality (PI). Let us first prove that for each~$n \in \mathbb{N}$, there exists a uniform constant~$C(2n) \in (0, \infty)$ with
\begin{align} \label{e_uniform_bound_of_2n_power}
\mathbb{E} \left[ \left| X_i -m_i  \right|^{2n} \right] \leq C(2n).
\end{align}
Note that~\eqref{e_uniform_bound_of_2n_power} is true for~$n=1$ by Lemma~\ref{l_variance estimates}. Then the following observation and induction on~$n$ imply~\eqref{e_uniform_bound_of_2n_power} for all~$n\in\mathbb{N}$:
\begin{align}
\mathbb{E}\left[ \left| X_i -m_i  \right|^{2n+2} \right] &= \var \left( \left| X_i -m_i \right|^{n+1} \right) + \mathbb{E} \left[ \left| X_i -m_i \right|^{n+1} \right]^2 \\
& \overset{PI}{\leq} \frac{1}{\rho} \int (n+1)^2 \left| X_i -m_i \right|^{2n} d\mu^{\sigma} + \mathbb{E} \left[ \left| X_i -m_i \right|^{n+1} \right]^2 \\
& \lesssim \mathbb{E}  \left[ \left| X_i -m_i \right|^{2n} \right]+ \mathbb{E}  \left[ \left| X_i -m_i \right|^{n+1} \right]^2.
\end{align}
Note also that a combination of~\eqref{e_uniform_bound_of_2n_power} and H\"older's inequality implies for each~$n\in \mathbb{N}$, there is a constant~$C(n)$ with
\begin{align}\label{e_uniform_bound_of_n_power}
\mathbb{E} \left[ \left| X_i -m_i  \right|^{n} \right] \leq C(n).
\end{align}
Then the lemma now follows from arithmetic-geometric mean inequality:
\begin{align}
&\left|\mathbb{E} \left[ \sum_{i_1 \in A_{1}}\cdots \sum_{i_k \in A_{k}} \left(X_{i_1}-m_{i_1} \right) \cdots \left(X_{i_k}- m_{i_k} \right) \right]\right| \\
& \qquad  \leq \sum_{i_1 \in A_{1}}\cdots \sum_{i_k \in A_{k}} \mathbb{E} \left[ \frac{1}{k}\left|\left(X_{i_1}-m_{i_1} \right)\right|^k +  \cdots + \frac{1}{k}\left|\left(X_{i_k}- m_{i_k} \right)\right|^k \right] \\
& \qquad \overset{\eqref{e_uniform_bound_of_n_power}}{\lesssim} \left|A_{1}\right|\cdots\left|A_k \right|.
\end{align}
\qed
\medskip

The next auxiliary lemma states that on a one dimensional lattice with nearest-neighbor interaction, the gce has uniform exponential decay of correlations.
\begin{lemma}\label{l_exponential decay of correlations}
For a function~$f:\mathbb{R}^K \rightarrow \mathbb{R}$, denote~$\text{supp} \ f = \{i_1, i_2, \cdots, i_k\}$ by the minimal subset of~$\{1, 2, \cdots, K\}$ with $f(x) = f(x_{i_1}, \cdots x_{i_k})$. Then for any~$f, g :\mathbb{R}^K \rightarrow \mathbb{R}$,
\begin{align}
&\left|\cov\left(f(X), g(X) \right) \right| \\
& \qquad \lesssim \left( \int | \nabla f |^2 d \mu^{\sigma}  \right)^{\frac{1}{2}}  \left( \int |\nabla g |^2 d \mu^\sigma \right)^{\frac{1}{2}} \exp\left(-C{\text{dist}\left(\text{supp} \ f, \text{supp}\ g\right)   }\right)  .\label{e_exp_decay}
\end{align}
\end{lemma}

\textsc{Proof of Lemma~\ref{l_exponential decay of correlations}}. \ As noted in the proof of Lemma~\ref{l_variance estimates}, gce~$\mu^{\sigma}$ satisfies a uniform logarithmic Sobolev inequality. Then~\cite[Theorem 2.1]{Yos01} implies there exist constants~$C>0$ and~$C(f,g)>0$ with
\begin{align} \label{e_exp_decay_yos}
\left| \cov \left( f(X), g(X) \right) \right| \leq C(f, g) \exp\left(-C{\text{dist}\left(\text{supp} \ f, \text{supp}\ g\right)   }\right),
\end{align}
where~$C(f,g)$ only depends on~$\|f\|_{\infty}, \|g \|_{\infty}, \text{supp} \ f$ and~$\text{supp} \ g$ (see also~\cite[Theorem 2.9]{HeMe16}). On the right hand side of~\eqref{e_exp_decay_yos}, the constant~$C(f,g)$ depends on~$\|f\|_{\infty}, \|g \|_{\infty}$. Here, with a small change in the proof, one can easily deduce the desired inequality~\eqref{e_exp_decay}. To convince the reader, there is a similar result for algebraic decaying interactions contained in~\cite{HeMe16}. 
\qed
\medskip

Now, we get to the core estimates needed for the proof of Theorem~\ref{p_equivalence_of_ensembles} and of Theorem~\ref{p_local_cramer}. 
\begin{proposition}\label{p_main computation}
For each~$\alpha>0$ and~$\beta>\frac{1}{2}$, there exists a uniform constant~$0<C<\infty$ and an integer~$K_0\in \mathbb{N}$ such that for all~$K \geq K_0$ and all~$\sigma \in \mathbb{R}$
\begin{align}
\frac{1}{C} \leq g_{K, m} (0) \leq C, \label{0th derivative}\\
\left| \frac{d}{d \sigma}g_{K, m} (0) \right| \lesssim K^{\alpha},\label{1st derivative}\\
\left| \frac{d^2}{d \sigma^2}g_{K, m} (0) \right| \lesssim K^{\beta}.\label{2nd derivative}
\end{align}
\end{proposition}

The statement of Proposition~\ref{p_main computation} should be compared to Proposition 31 in~\cite{GrOtViWe09} or Proposition 3.1 in~\cite{MeOt13}. The main difference is that in our situation the random variables~$X_1,\ldots, X_K$ are dependent. This also makes the proof of Proposition~\ref{p_main computation} a lot harder.\\

The estimates of Proposition~\ref{p_main computation} are motivated from deducing a quantitative local central limit theorem for the properly normalized sum of the random variables $X_1, \ldots, X_K$. For example, if the random variables~$X_1, \ldots, X_K$ are iid, the estimate~\eqref{0th derivative} is a weaker version of the quantitative local clt estimate 
\begin{align}
\left| g_{K, m} (0) - \frac{1}{\sqrt{2 \pi}} \right| \lesssim \frac{1}{\sqrt{K}}.
\end{align}
The last inequality states that the density of the normalized sum at point~$0$ converges to the density of the normal distribution. As we mentioned in the introduction, we believe that one could strengthen the estimates of Proposition~\ref{p_main computation} to get a local central limit theorem for dependent random variables. However, we choose to derive weaker bounds instead because they are sufficiently strong for deducing our main results (see Theorem~\ref{p_equivalence_of_ensembles} and Theorem~\ref{p_local_cramer}). Deducing those weaker estimates is already quite subtle and challenging.\\

We deduce Proposition~\ref{p_main computation} in Section~\ref{s_main_computation}. There, we also comment on how to overcome the problem of considering dependent random variables and not independent ones. Now, we are prepared for the proof of Theorem~\ref{p_equivalence_of_ensembles}.\\

\textsc{Proof of Theorem~\ref{p_equivalence_of_ensembles}}. \ Let us begin with the estimate~\eqref{e_eoe_C_0}. From~\eqref{cramers representation1}, we have
\begin{align} \label{e_cramer_rep_log_form}
A_{ce}(\sigma)-A_{gce}(\sigma) = \frac{1}{K} \ln g_{K,m}(0).
\end{align}
Then a combination of~\eqref{0th derivative} and~\eqref{e_cramer_rep_log_form} yields, as desired,
\begin{align}
\left| A_{ce}(\sigma)-A_{gce}(\sigma) \right| \lesssim \frac{1}{K}.
\end{align}
Let us turn to the estimate~\eqref{e_eoe_C_1}. Taking the derivative with respect to~$\sigma$ in~\eqref{e_cramer_rep_log_form} yields
\begin{align}\label{e_cramer_rep_1st_der}
\frac{d}{d\sigma} A_{ce}(\sigma)- \frac{d}{d\sigma} A_{gce}(\sigma) = \frac{1}{K} \frac{1}{g_{K,m}(0)} \frac{dg_{K,m}(0)}{d\sigma}.
\end{align}
Let us choose~$\alpha=\varepsilon$. Then a combination of~\eqref{0th derivative},~\eqref{1st derivative} and~\eqref{e_cramer_rep_1st_der} implies
\begin{align}
 \left| \frac{d}{d\sigma} A_{ce}(\sigma)- \frac{d}{d\sigma} A_{gce}(\sigma) \right| \lesssim \frac{1}{K}K^{\alpha}= \frac{1}{K^{1-\varepsilon}}.
\end{align}
Let us turn to the estimate~\eqref{e_eoe_C_2}. Differentiating~\eqref{e_cramer_rep_1st_der} again, we get
\begin{align}
\frac{d^2}{d\sigma ^2} A_{ce}(\sigma)- \frac{d^2}{d\sigma ^2} A_{gce}(\sigma) &= -\frac{1}{K} \frac{1}{\left(g_{K,m}(0) \right)^2} \left(\frac{dg_{K,m}(0)}{d\sigma} \right)^2 + \frac{1}{K} \frac{1}{g_{K,m}(0)} \frac{d^2 g_{K,m}(0)}{d\sigma^2}.
\end{align}
Then after choosing~$\beta= \frac{1}{2} + \varepsilon$, a combination of~\eqref{0th derivative},~\eqref{1st derivative} and~\eqref{2nd derivative} yields
\begin{align}
\left|\frac{d^2}{d\sigma ^2} A_{ce}(\sigma)- \frac{d^2}{d\sigma ^2} A_{gce}(\sigma) \right| \lesssim \frac{1}{K^{1-2\alpha}} + \frac{1}{K^{1- \beta}}  \lesssim \frac{1}{K^{\frac{1}{2}-\varepsilon}}.
\end{align}
\qed

Let us proceed to the proof of Theorem~\ref{p_local_cramer}. \\

\textsc{Proof of Theorem~\ref{p_local_cramer}}. \ 
Recall the difference of the free energies~\eqref{e_rewriting_difference_free_energies} of $A_{gce}$ and~$A_{ce}$
\begin{align}
A_{gce}(\sigma) - A_{ce}(\sigma) = \bar{H}(m)-\mathcal{H}_K (m). \end{align}
Then the first desired estimate~\eqref{e_theorem1 1} follows from a combination of~\eqref{e_rewriting_difference_free_energies} and~\eqref{e_eoe_C_0} in Theorem~\ref{p_equivalence_of_ensembles}. Let us turn to the estimate~\eqref{e_theorem1 2}. A direct computation yields
\begin{align}
\frac{d}{dm} \left(\mathcal H_K (m) - \bar H (m)\right)
&= \frac{d}{dm}\left( A_{ce}(\sigma) - A_{gce}(\sigma) \right) \\
& = \frac{d}{d\sigma}\left( A_{ce}(\sigma) - A_{gce}(\sigma) \right) \frac{d\sigma}{dm} \\
& = \frac{d}{d\sigma}\left( A_{ce}(\sigma) - A_{gce}(\sigma) \right) \left(\frac{dm}{d\sigma}\right)^{-1} \label{e_dhdsigma}.
\end{align}
Then~\eqref{e_eoe_C_1},~\eqref{e_dmdsigma} and Lemma~\ref{l_variance estimates} implies, as desired,
\begin{align}
\left|\frac{d}{dm} \left(\mathcal H_K (m) - \bar H (m)\right) \right| \lesssim  \frac{1}{K^{1-\varepsilon}}.
\end{align}

Before we move on to the estimate~\eqref{e_theorem1 3}, let us deduce some auxiliary results. A direct calculation yields
\begin{align}
\frac{d}{d \sigma} \mathbb{E}\left[f(X)\right] & = \frac{d}{d\sigma} \int f(x) \frac{ \exp \left( \sum_{i=1}^K \sigma x_i -H(x)\right)}{\int\exp \left( \sum_{i=1}^K \sigma x_i -H(x)\right) dx } dx \\
& = \int \frac{df}{d\sigma} (x) \frac{ \exp \left( \sum_{i=1}^K \sigma x_i -H(x)\right)}{\int\exp \left( \sum_{i=1}^K \sigma x_i -H(x)\right) dx } dx \\
& \quad + \int f(x) \left( \sum_{i=1}^{K} \left(x_i - m_i\right) \right) \frac{ \exp \left( \sum_{i=1}^K \sigma x_i -H(x)\right)}{\int\exp \left( \sum_{i=1}^K \sigma x_i -H(x)\right) dx } dx \\
& = \mathbb{E}\left[\frac{df}{d\sigma}(X) \right] + \mathbb{E} \left[ f(X) \left( \sum_{i=1}^{K}\left( X_i -m_i\right) \right) \right]. \label{e_derivative_formula}
\end{align}
In particular we have
\begin{align} 
\frac{d}{d\sigma}m_k &=\frac{d}{d\sigma} \mathbb{E} \left[ X_k \right]= \mathbb{E} \left[X_k \left( \sum_{i=1}^K \left(X_i -m_i \right) \right)\right] = \mathbb{E} \left[\left(X_k -m_k\right) \left( \sum_{i=1}^K \left(X_i -m_i \right) \right)\right]\label{e_derivative_of_mk}
\end{align}
and
\begin{align}
&\frac{d}{d\sigma} \left( \mathbb{E} \left[ \left( \sum_{i=1}^K \left(X_i -m_i \right) \right) ^2 \right] \right) \\
&\overset{\eqref{e_derivative_formula}}{=} \mathbb{E} \left[ \frac{d}{d\sigma} \left( \sum_{i=1}^K \left(X_i -m_i \right) \right) ^2\right] + \mathbb{E} \left[ \left( \sum_{i=1}^K \left(X_i -m_i \right) \right) ^3\right] \\
& \overset{\eqref{e_derivative_of_mk}}{=} \mathbb{E} \left[ 2 \left(\sum_{i=1}^K \left(X_i -m_i \right) \right)\left(- \sum_{i=1}^K \mathbb{E} \left[\left( X_i -m_i \right)\sum_{k=1}^K \left(X_k -m_k \right) \right] \right) \right] + \mathbb{E} \left[ \left( \sum_{i=1}^K \left(X_i -m_i \right) \right) ^3\right] \\
&  =-2 \mathbb{E} \left[\left( \sum_{k=1}^K \left(X_k -m_k \right) \right)^2 \right]  \mathbb{E} \left[ \sum_{i=1}^K \left(X_i -m_i \right) \right] + \mathbb{E} \left[ \left( \sum_{i=1}^K \left(X_i -m_i \right) \right) ^3\right] \\
& \overset{\eqref{d_def_mi}}{=} 0+ \mathbb{E} \left[ \left( \sum_{i=1}^K \left(X_i -m_i \right) \right) ^3\right]  = \mathbb{E} \left[ \left( \sum_{i=1}^K \left(X_i -m_i \right) \right) ^3\right]. \label{e_derivative_second_power}
\end{align}
Now we claim that
\begin{align}
\left|\frac{d}{d\sigma}\left(\frac{d \sigma}{dm}\right)\right| \lesssim 1.
\end{align}
To begin with, it follows from the calculation from above that
\begin{align}
\frac{d}{d\sigma}\left(\frac{d \sigma}{dm}\right) & \overset{\eqref{e_dmdsigma}}{=} \frac{d}{d\sigma}\left( \frac{K}{\mathbb{E} \left[  \left(\sum_{i=1}^{K} \left(X_i -m_i \right) \right)^2 \right]}\right)\\
& = - \frac{K}{\mathbb{E} \left[  \left(\sum_{i=1}^{K} \left(X_i -m_i \right) \right)^2 \right]^2} \frac{d}{d\sigma}\left(\mathbb{E} \left[  \left(\sum_{i=1}^{K} \left(X_i -m_i \right) \right)^2 \right]\right) \\
&  \overset{\eqref{e_derivative_second_power}}{=} - \frac{K}{\left( \var\left(\sum_{i=1}^{K} X_i \right)\right)^2} \mathbb{E} \left[  \left(\sum_{i=1}^{K} \left(X_i -m_i \right) \right)^3 \right]. \label{e_dsigmadm_estimation}
\end{align}
Note that a direct computation yields
\begin{align}
\left|\mathbb{E} \left[ \left(\sum_{i=1}^{K} \left(X_i -m_i \right) \right)^3 \right]\right| & \lesssim \sum_{i\leq j \leq k} \left|\mathbb{E} \left[ \left(X_i -m_i \right) \left(X_j -m_j \right) \left(X_k -m_k \right)\right]\right|\\
&  \lesssim \sum_{j=1}^{K} \sum_{s=0}^{j} \sum_{t=0}^{K-j} \left| \mathbb{E}\left[ \left(X_{j-s} -m_{j-s} \right) \left(X_j -m_j \right) \left(X_{j+t} -m_{j+t} \right)\right]\right|\\
& \leq \sum_{j=1}^{K} \sum_{s=0}^{K} \sum_{t=0}^{s} \left| \mathbb{E}\left[ \left(X_{j-s} -m_{j-s} \right) \left(X_j -m_j \right) \left(X_{j+t} -m_{j+t} \right)\right]\right| \label{e_3rd moment1}\\
&  \quad + \sum_{j=1}^{K} \sum_{t=0}^{K} \sum_{s=0}^{t} \left| \mathbb{E}\left[ \left(X_{j-s} -m_{j-s} \right) \left(X_j -m_j \right) \left(X_{j+t} -m_{j+t} \right)\right]\right|.\label{e_3rd moment2}
\end{align}
Then Lemma~\ref{c_amgm} and Lemma~\ref{l_exponential decay of correlations} imply
\begin{align}
T_{\eqref{e_3rd moment1}} &= \sum_{j=1}^{K} \sum_{s=0}^{K} \sum_{t=0}^{s} \left| \cov\left( \left(X_j -m_j \right)\left( X_{j+t}-m_{j+t} \right), X_{j-s}-m_{j-s} \right)\right| \\
& \lesssim \sum_{j=1}^{K} \sum_{s=0}^{K} \sum_{t=0}^{s} \exp\left(-Cs\right)  = K \sum_{s=0}^{K} s \exp\left(-Cs\right) \lesssim K. \label{e_estimate_e_3rd moment1}
\end{align}
A similar argument gives
\begin{align}T_{\eqref{e_3rd moment2}} \lesssim K \label{e_estimate_e_3rd moment2}
\end{align}
Therefore a combination of~\eqref{e_dsigmadm_estimation},~\eqref{e_estimate_e_3rd moment1},~\eqref{e_estimate_e_3rd moment2}, and Lemma~\ref{l_variance estimates} yields, as desired,
\begin{align}
\left|\frac{d}{d\sigma}\left(\frac{d \sigma}{dm}\right)\right|&= \left| \frac{K}{\left(\var\left(\sum_{i=1}^K X_i\right)\right)^2} \right| \left|\mathbb{E}\left[\left( \sum_{i=1}^K \left(X_i -m_i \right) \right)^3 \right] \right| \lesssim \frac{K}{K^2}K =1.\label{e_ddsigmadsigmadm}
\end{align} 

Let us turn to the estimate~\eqref{e_theorem1 3}. We differentiate~\eqref{e_dhdsigma} to obtain
\begin{align}
&\frac{d^2}{dm^2} \left(\mathcal H_K (m) - \bar H (m)\right)\\
&\qquad = \frac{d}{dm} \left( \frac{d}{d\sigma}\left( A_{ce}(\sigma) - A_{gce}(\sigma) \right) \frac{d\sigma}{dm} \right) \\
& \qquad = \frac{d}{d\sigma}\left( \frac{d}{d\sigma}\left( A_{ce}(\sigma) - A_{gce}(\sigma) \right) \frac{d\sigma}{dm} \right)  \frac{d\sigma}{dm} \\
&\qquad = \frac{d^2}{d\sigma^2} \left( A_{ce}(\sigma) - A_{gce}(\sigma) \right) \left( \frac{d\sigma}{dm}\right)^2 + \frac{d}{d\sigma}\left(A_{ce}(\sigma) - A_{gce}(\sigma) \right) \frac{d}{d\sigma}\left(\frac{d\sigma}{dm}\right) \frac{d\sigma}{dm}.
\end{align}
Then a combination of~\eqref{e_dmdsigma},~\eqref{e_ddsigmadsigmadm}, Theorem~\ref{p_equivalence_of_ensembles} and Lemma~\ref{l_variance estimates} yields
\begin{align}
\left|\frac{d^2}{dm^2}\left(\mathcal H_K (m) - \bar H (m)\right)\right|\lesssim \frac{1}{K^{\frac{1}{2}-\varepsilon}},
\end{align}
and this finishes the proof of Theorem~\ref{p_local_cramer}.
\qed

\section{Proof of Proposition~\ref{p_main computation}}\label{s_main_computation}

The proof of Proposition~\ref{p_main computation} represents the core of our argument. Before turning to the precise argument let us motivate and explain our approach in more detail. We will especially emphasize on how the problem of considering dependent and not independent random variables~$X_l$ is solved.\\

As we mentioned before, the argument is inspired from deducing local central limit theorems via the Fourier inversion method (see~\cite{Fel71} or~\cite{MeOt13}). The main idea of this method is to write the the density of the random variable~$$Z= \frac{1}{\sqrt{K}} \sum_{l=1}^K X_l$$ by Fourier inversion as an integral involving the characteristic function (see~\eqref{Inverse Fourier representation} below)~$$\varphi_Z(\xi) = \mathbb{E}\left[ \exp (i \xi Z)\right] .$$
The next step is to split up the integral into an inner integral over the interval~$|\xi| \leq \delta \sqrt{K}$ and an outer integral over the interval~$|\xi| > \delta \sqrt{K}$. The outer integral usually is an error term and the main contribution comes from the inner integral.\\

The big advantage of considering independent random variables~$X_l$ is that the characteristic function~$\varphi_Z$ becomes a product of the characteristic functions~$\varphi_{X_l}$ i.e.
\begin{align}
 \varphi_Z(\xi) & = \mathbb{E}\left[ \exp \left( i \xi \frac{1}{\sqrt{K}} \sum_{l=1}^K X_l  \right)\right]  = \prod_{l=1}^K \varphi_{X_l} \left( \frac{\xi}{\sqrt{K}} \right). \label{e_product_characteristic_function} 
\end{align}
Then the outer integral is small because each characteristic function~$\varphi_{X_l}<1$ is small and decays at least of the order~$|\xi|^{-1}$. For the inner integral, one applies a Taylor expansion onto the functions~$\ln \varphi_{X_l} $ and gets the correct contribution due to the normalization of the random variables. This strategy would yield the desired estimate~\eqref{0th derivative} in the case of independent random variables (see also Section~3 in~\cite{MeOt13}).\\

For deducing the estimates~\eqref{1st derivative} and~\eqref{2nd derivative} in the case of independent random variables one proceeds in a similar way. The obtained integral representation is split up into an inner and an outer integral. One shows that the outer integral is small by using decay of the characteristic functions. The inner integral is estimated again by Taylor expansion. However, the situation becomes more subtle when considering dependent random variables. The obtained integral representation involves several new terms that look like covariances i.e.~they are covariances if~$\xi=0$. Therefore, we have to be a lot more careful when applying this strategy. \\

The following observation helps a lot when considering dependent random variables: Because we only consider nearest-neighbor interaction, the odd random variables~$X_{odd}$ become independent if we condition on the values of the even random variables~$X_{even}$. Additionally, because our Hamiltonian is quadratic the variances of the conditioned random variables~$X_{odd}| X_{even}$ are uniformly bounded from above and from below (see proof of Lemma~\ref{l_exchanging exponential terms} ).\\

Using this observation the outer integral can be estimated in a straight-forward manner. We condition on the even random variables~$X_{even}$. By conditional independence we get that the conditional characteristic function becomes a product i.e.
\begin{align}
 \varphi_Z(\xi) & = \mathbb{E}\left[ \exp \left( i \xi \frac{1}{\sqrt{K}} \sum_{l=1}^K X_l  \right)\right]  \\
 & = \mathbb{E}\left[ \exp \left( i \xi \frac{1}{\sqrt{K}} \sum_{j : even} X_j  \right)  \mathbb{E}\left[\exp \left( i \xi \frac{1}{\sqrt{K}} \sum_{i : odd} X_i  \right)  | X_{j}, j : even \right]\right] \\
 & = \mathbb{E}\left[ \exp \left( i \xi \frac{1}{\sqrt{K}} \sum_{j : even} X_j  \right) \prod_{i : odd}  \mathbb{E}\left[\exp \left( i \xi \frac{1}{\sqrt{K}} X_i  \right)  | X_{j}, j : even\right]\right].
\end{align}
Because the variances of the conditional random variables~$X_{odd} |X_{even}$ are controlled uniformly in the conditioned values~$X_{even}$ we have that the conditional characteristic functions
\begin{align}
\mathbb{E}\left[\exp \left( i \xi \frac{1}{\sqrt{K}} X_k  \right)  | X_{j}, j : even\right]
\end{align}
decay uniformly (see Lemma~\ref{l_exchanging exponential terms} below). Over-simplifying the argument, this yields the correct bounds on the outer integrals.\\

The situation for the inner integrals is more tricky and one has to proceed differently for the estimate~\eqref{0th derivative} and for the estimates~\eqref{1st derivative} and~\eqref{2nd derivative}. Let us first consider the argument for~\eqref{0th derivative}. In the inner integral, we condition on the even random variables~$X_{even}$. We use the conditional independence and the control on the conditional variances to do a Taylor expansion just for the characteristic functions of the conditional random variables~$X_{odd}|X_{even}$. Then we show that this suffices to get the desired estimate of~\eqref{0th derivative}.\\

Let us turn to~\eqref{1st derivative} and~\eqref{2nd derivative} and explain how the inner integrals are estimated there. As mentioned above, the Taylor expansion becomes a lot more tricky than for~\eqref{0th derivative}. For each additional derivative, the argument becomes more and more elaborate. The reason is that whenever calculating the inner and outer integral one ends up with more and more error terms. The first step of the argument is to carefully group those terms such that certain terms cancel and other terms become covariance-like. This means that those terms are a covariance if~$\xi=0$. However, if~$\xi \neq 0$ they are not covariances and cannot be estimated by the decay of correlations. We are able to estimate those terms using the following idea: For each error term, we partition the sites of the lattice system into blocks (see Figure~\ref{f_sets_E} and Figure~\ref{f_sets_F} below). Then we carry out a multivariate Taylor expansion. Let's say we expand the function~$F(\xi_1, \xi_2)$ and after expanding we set~$\xi_1=\xi$ and~$\xi_2=\xi$. The variable~$\xi_1$ corresponds to the sites within the block and the variable~$\xi_2$ corresponds to terms outside of the block. We carry out the Taylor expansion with respect to~$\xi_1$. The resulting terms are controlled either by the help of decay of correlations or by the size of the blocks. \\

The proof is organized in the following way. In Section~\ref{s_est_char_fct} we deduce auxiliary estimates for the conditional characteristic functions.
In Section~\ref{s_0th_derivative} we deduce the estimate~\eqref{0th derivative}. In Section~\ref{s_2nd_derivative} we verify the estimate~\eqref{2nd derivative}. The estimate~\eqref{1st derivative} can be derived by similar arguments. In Section~\ref{s_proof_of_aux_lemmas_inner} we deduce auxiliary lemmas used in Section~\ref{s_2nd_derivative}.

\subsection{Auxiliary estimates}\label{s_est_char_fct} In this section we provide some auxiliary estimates that are needed in the proof of Proposition~\ref{p_main computation}. 
Let us introduce the auxiliary sets (cf.~Figure~\ref{f_sets_E})
\begin{align}
&E_1 ^l : = \{ k\ | \ |k-l| \leq L\} \label{e_E1l} \\
& E_2 ^l : = \{ k\  | \ |k-l|>L \} \label{e_E2l}
\end{align}
and
(cf.~Figure~\ref{f_sets_F})
\begin{align}
&F_1 ^{n, l} : = \{ k\ | \ |k-n| \leq L \text{ or } |k-l| \leq L \} \label{e_F1nl}\\
& F_2 ^{n, l} : = \{ k\ | \ |k-n| > L \text{ and } |k-l| > L \}, \label{e_F2nl}
\end{align}
where~$L\ll K$ is a positive integer that will be chosen later.\\

Let us also introduce the auxiliary notations. Let us denote~$m_{i,2}$ and~$s_{i,2}^2$ to be
\begin{align}
m_{i, 2} &: = \mathbb{E}\left[ X_i \mid X_j, j : even \right], \label{e_def_of_ml2} \\
s_{i,2}^2 &: = \mathbb{E} \left[ \left( X_i -m_{i, 2} \right)^2 \mid X_j, j : even \right]. \label{e_def_si2}
\end{align}
Define the function~$\hat{e} : \mathbb{R} \rightarrow \mathbb{C}$ by
\begin{align} \label{e_def_of_hat_e}
\hat{e}\left(\xi\right) := \exp \left( i \frac{1}{\sqrt{K}}\sum_{k : even} \left( X_k - m_k \right) \xi + i \frac{1}{\sqrt{K}} \sum_{i : odd} \left( m_{i, 2} -m_i \right) \xi\right).
\end{align}
Then we have the following lemma:

\begin{figure}[t]
\centering
\begin{tikzpicture}[xscale=1.3]

\draw[fill] (0,0) circle [radius=0.05];
\draw[fill] (.25,0) circle [radius=0.05];
\draw[fill] (.5,0) circle [radius=0.05];
\draw[fill] (.75,0) circle [radius=0.05];
\draw[fill] (1,0) circle [radius=0.05];

\draw[fill] (1.475,0) circle [radius=0.02];
\draw[fill] (1.625,0) circle [radius=0.02];
\draw[fill] (1.775,0) circle [radius=0.02];

\draw[fill] (2.25,0) circle [radius=0.05];
\draw[fill] (2.5,0) circle [radius=0.05];
\draw[fill] (2.75,0) circle [radius=0.05];
\draw[fill] (3,0) circle [radius=0.05];
\draw[fill] (3.25,0) circle [radius=0.05];

\draw (3.5,0) circle[radius=0.05];
\draw (3.75,0) circle[radius=0.05];
\draw (4,0) circle[radius=0.05];
\draw (4.25,0) circle[radius=0.05];
\draw (4.5,0) circle[radius=0.05];
\draw (4.75,0) circle[radius=0.05];
\draw (5,0) circle[radius=0.05];
\draw (5.25,0) circle[radius=0.05];
\draw (5.5,0) circle[radius=0.05];

\draw[fill] (5.75,0) circle [radius=0.05];
\draw[fill] (6,0) circle [radius=0.05];
\draw[fill] (6.25,0) circle [radius=0.05];
\draw[fill] (6.5,0) circle [radius=0.05];
\draw[fill] (6.75,0) circle [radius=0.05];

\draw[fill] (7.225,0) circle [radius=0.02];
\draw[fill] (7.375,0) circle [radius=0.02];
\draw[fill] (7.525,0) circle [radius=0.02];

\draw[fill] (8,0) circle [radius=0.05];
\draw[fill] (8.25,0) circle [radius=0.05];
\draw[fill] (8.5,0) circle [radius=0.05];
\draw[fill] (8.75,0) circle [radius=0.05];
\draw[fill] (9,0) circle [radius=0.05];
\node[align=center, below] at (0,-.05) {$1$};
\node[align=center, below] at (4.5,-.05) {$l$};
\node[align=center, below] at (9,-.05) {$K$};
\node[align=center, below] at (4.5,-1.2) {};

\draw[decorate,decoration={brace,mirror},thick] (0,-.5) -- node[below]{$E_{2}^{l}$} (3.25,-.5);
\draw[decorate,decoration={brace,mirror},thick] (3.5,-.5) -- node[below]{$E_{1}^{l}$} (5.5,-.5);
\draw[decorate,decoration={brace,mirror},thick] (5.75,-.5) -- node[below]{$E_{2}^{l}$} (9,-.5);

\draw[thick,<->] (3.5,.5) -- (4.45,.5);
\node[align=center, above] at (4,.5) {$L$};
\draw[thick,<->] (4.55,.5) -- (5.5,.5);
\node[align=center, above] at (5,.5) {$L$};
\end{tikzpicture}
\caption{The sets $E_1^l$ and $E_2 ^l$ given by~\eqref{e_E1l} and~\eqref{e_E2l}}\label{f_sets_E}
\end{figure}
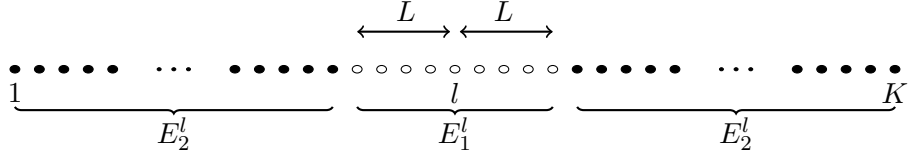

\begin{lemma}\label{l_exchanging exponential terms}
For large enough~$K$ and~$\delta>0$ small enough, there exists a positive constant~$C>0$ such that the following inequalities hold for all~$\xi \in \mathbb{R}$ with~$ \frac{\left|\xi\right|}{\sqrt{K}}  \leq \delta$.
\begin{align}
&\left|\mathbb{E}\left[\exp\left( i  \sum_{k =1}^K \left(X_k -m_k \right) \frac{\xi}{\sqrt{K}}\right) \right] \right. \\
&\qquad  \qquad \left.- \mathbb{E} \left[ \hat{e}\left(\xi\right) \exp\left( - \sum_{i : odd} \frac{ s_{i, 2}^2}{2K}\xi^2 \right) \right]\right|  \lesssim \frac{1}{\sqrt{K}} \left|\xi\right|^3 \exp\left(-C\xi^2 \right),  \label{e_lemma exchanging 0th der}  \\
& \left|\mathbb{E}\left[\exp\left( i \sum_{k \in E_2^l } \left( X_k -m_k \right) \frac{\xi}{\sqrt{K}}\right) \mid \mathscr{F}_l \right] \right| \lesssim \left( 1 + \xi^2 \right)\exp\left(-C\xi^2 \right),
\label{e_lemma exchanging 1st der} \\
& \left|\mathbb{E}\left[\exp\left( i \sum_{k \in F_2^{n,l} } \left( X_k -m_k \right) \frac{\xi}{\sqrt{K}}\right) \mid \mathscr{G}_{n,l} \right] \right| \lesssim \left( 1 + \xi^2 \right)\exp\left(-C\xi^2 \right) .
\label{e_lemma exchanging 2nd der}
\end{align}
where~$\mathscr{F}_l$,~$\mathscr{G}_{n,l}$ denote the sigma algebras defined by
\begin{align}
\mathscr{F}_l : = \sigma\left(X_k , \ k \in E_1 ^l \right) \qquad \text{and} \qquad \mathscr{G }_{n,l} : = \sigma \left( X_k , \ k \in F_1 ^{n,l} \right).
\end{align} 
\end{lemma}
Lemma~\ref{l_exchanging exponential terms} will be used in the estimation of the inner integrals when deriving~\eqref{0th derivative},~\eqref{1st derivative} and~\eqref{2nd derivative}. 

\medskip

\textsc{Proof of Lemma~\ref{l_exchanging exponential terms}}. \
We first deduce~\eqref{e_lemma exchanging 0th der}. Let us consider the conditional expectation with respect to~$\{X_j \mid j : even\}$. In the case of nearest-neighbor interaction, the conditional Lebesgue density~
$\mu^{\sigma}(dx_1 dx_3 \cdots \mid x_j, \ j : even)$ can be written as
\begin{align}
&\mu^{\sigma}(dx_1 dx_3 \cdots \mid x_j, \ j : even) \\
&\qquad  = \frac{1}{Z}\exp\left( \sum_{i : odd} -\psi(x_i) - (-Jx_{i-1}-Jx_{i+1}+s_i  - \sigma )x_i \right) \\
& \qquad = \prod_{i : odd} \frac{1}{Z} \exp\left( -\psi(x_i) - (-Jx_{i-1}-Jx_{i+1}+s_i  - \sigma )x_i \right), \label{e_conditional_density_odd_variables}
\end{align}
which implies that~$\{\mathbb{E}\left[ X_i \mid X_j , \ j : even\right] \mid i : odd \}$ are independent and in particular,
\begin{align}
&\mu^{\sigma}\left(dx_i | x_2, x_4, \cdots \right)= \frac{1}{Z} \exp \left( -\psi(x_i) - (-Jx_{i-1}-Jx_{i+1}+s_i  - \sigma )x_i  \right). \label{e_conditional leb for odd}
\end{align}
Note that
\begin{align}
\frac{d}{d\sigma} \mathbb{E} \left[  X_i  \mid X_j, j : even \right] &= \frac{d}{d\sigma} \int x_i \cdot \frac{1}{Z} \exp \left( -\psi(x_i) - (-Jx_{i-1}-Jx_{i+1}+s_i  - \sigma )x_i  \right) dx\\
& = \int x_i \left( x_i - m_{i,2} \right) \cdot \frac{1}{Z} \exp \left( -\psi(x_i) - (-Jx_{i-1}-Jx_{i+1}+s_i  - \sigma )x_i  \right) dx \\
& =\mathbb{E} \left[ \left( X_i -m_{i, 2} \right)^2 \mid X_j, j : even \right] = s_{i,2}^2, \label{e_def_derivative_si2}
\end{align}
and a similar computation yields
\begin{align}
\frac{d^2}{d\sigma^2} \mathbb{E} \left[  X_i  \mid X_j, j : even \right] = \mathbb{E} \left[ \left( X_i -m_{i, 2} \right)^3 \mid X_j, j : even \right] = : t_{i,2}. \label{e_def_ti2}
\end{align}
From the observations~\eqref{e_conditional_density_odd_variables} and~\eqref{e_conditional leb for odd} from above, we get the product structure of conditional expectation
\begin{align}
&\mathbb{E}\left[ \exp\left( i \frac{1}{\sqrt{K}}\sum_{k=1}^K \left( X_k - m_k \right) \xi \right) \right] \\
& \qquad = \mathbb{E}\left[ \exp\left( i \frac{1}{\sqrt{K}}\sum_{k : even} \left( X_k - m_k \right) \xi \right)\right. \\
&\qquad \quad \quad \quad \quad \quad \quad \left.\times \mathbb{E}\left[ \exp\left(i \frac{1}{\sqrt{K}} \sum_{i : odd} \left(X_i -m_i \right)\xi\right) \mid X_j, j : even \right] \right] \\
& \qquad = \mathbb{E} \left[ \exp\left( i \frac{1}{\sqrt{K}}\sum_{k : even} \left( X_k - m_k \right) \xi  + i \frac{1}{\sqrt{K}} \sum_{i : odd} \left( m_{i, 2} -m_i \right) \xi\right)\right. \\
&\qquad \quad  \quad \quad \quad \quad \quad \left.\times \mathbb{E} \left[\exp\left( i \frac{1}{\sqrt{K}} \sum_{ i : odd} \left( X_i -m_{i,2} \right) \xi \right) \mid X_j, j : even \right] \right] \\
&\qquad  \overset{\eqref{e_def_of_hat_e}}{=} \mathbb{E} \left[ \hat{e}\left(\xi\right) \prod_{ i : odd} \mathbb{E} \left[\exp\left( i \frac{1}{\sqrt{K}}  \left( X_i -m_{i,2} \right) \xi \right) \mid X_j, j : even \right] \right]. \label{e_proof for exchanging 0th der}
\end{align}
Let~$h_i $ denote the complex valued function 
\begin{align}
h_i (\xi) : = -\ln \left( \mathbb{E} \left[\exp\left( i  \left( X_i -m_{i,2} \right) \xi \right) \mid X_j, j : even \right]\right).
\end{align}
Equivalently, denote
\begin{align}
F_i (\xi) : &= \exp\left( - h_i (\xi) \right)  =  \mathbb{E} \left[\exp\left( i  \left( X_i -m_{i,2} \right) \xi \right) \mid X_j, j : even \right].
\end{align}
Differentiating both sides. we have
\begin{align}
F_i ' (\xi) & = -h_i ' (\xi) \exp\left( - h_i (\xi) \right)\\
& = i \mathbb{E} \left[  \left( X_i -m_{i,2} \right) \exp\left( i  \left( X_i -m_{i,2} \right) \xi \right) \mid X_j, j : even \right] ,
\end{align}
\begin{align}
F_i ''(\xi) &= -h_i '' (\xi) \exp\left( - h_i (\xi) \right) + h_i '(\xi)^2 \exp\left( - h_i (\xi) \right)\\
& = - \mathbb{E} \left[  \left( X_i -m_{i,2} \right)^2 \exp\left( i  \left( X_i -m_{i,2} \right) \xi \right) \mid X_j, j : even \right],
\end{align}
\begin{align}
F_i '''(\xi) & = -h_i '''( \xi) \exp\left( -h_i (\xi)\right) + 3h_i '' (\xi) h_i '(\xi) \exp\left( -h_i (\xi)\right) \\
& \quad - \left( h_i '(\xi)\right)^3 \exp\left(-h_i (\xi)\right) \\
& = -i\mathbb{E} \left[  \left( X_i -m_{i,2} \right)^3 \exp\left( i  \left( X_i -m_{i,2} \right) \xi \right) \mid X_j, j : even \right], \label{e_0th der Taylor third}
\end{align}
which implies~$h_i (0) = h_i ' (0) =0 $ and 
\begin{align}
h_i ''(0) = \mathbb{E} \left[  \left( X_i -m_{i,2} \right)^2 \mid X_j, j : even \right] \overset{\eqref{e_def_si2}}{=} s_{i,2}^2.
\end{align} 
Before we proceed, let us note that analogue of the equation~\eqref{e_uniform_bound_of_n_power} holds for the conditional expectation. More precisely, because the conditional measures~$\mu^{\sigma}\left(dx_i | x_2, x_4, \cdots \right)$ are bounded perturbations of a one dimensional Gaussian measure, it holds that for each~$n \in \mathbb{N}$ there is a constant~$C(n)$ with
\begin{align} \label{e_analogue of proving prop lemma}
\mathbb{E} \left[ \left| X_i -m_{i, 2} \right|^n \mid X_j, j : even \right] \leq C(n).
\end{align}
Note also that~\eqref{e_analogue of proving prop lemma} implies
\begin{align}
\left| F_i ' (\xi) \right| &\leq \mathbb{E} \left[ \left| X_i -m_{i, 2} \right| \mid X_j, j : even \right] \lesssim 1.
\end{align}
Combined with the fact that~$F_i (0) = 1$ we have for~$\left| \xi \right|$ small enough
\begin{align}
\frac{1}{2} \leq \left| F_i (\xi) \right|= \left| \exp\left(h_i (\xi) \right)\right| \leq \frac{3}{2}.
\end{align}
Inserting this into~\eqref{e_0th der Taylor third}, we obtain
\begin{align}
&\left| h_i ''' (\xi) \right| \\
& \qquad = \Big| 3h_i '' (\xi) h_i ' (\xi) - \left( h_i ' (\xi) \right) ^3  i \mathbb{E} \left[ \left( X_i - m_{i, 2} \right)^3 \exp\left( i \left( X_i -m_{i, 2} \right) \xi \right) \mid X_j , j : even \right] \Big| \\
& \qquad \quad \quad \times  \left| \exp\left(h_i (\xi) \right) \right| \\
& \qquad  \lesssim \mathbb{E} \left[ \left| X_i -m_{i, 2} \right| ^2 \mid X_j, j : even \right] \cdot \mathbb{E} \left[ \left| X_i -m_{i, 2} \right| \mid X_j, j : even \right]  \\
& \qquad \quad + \left( \mathbb{E} \left[\left| X_i -m_{i, 2}\right| \mid X_j, j : even \right] \right)^3 + \mathbb{E} \left[ \left| X_i -m_{i, 2} \right| ^3 \mid X_j, j : even \right]  \\
& \qquad  \overset{\eqref{e_analogue of proving prop lemma}}{\lesssim} 1.
\end{align}
We thus have for~$\left|\xi \right|$ small,
\begin{align}
\left| h_i ( \xi)  - \frac{1}{2}s_{i, 2}^2  \xi^2 \right| \lesssim \left|\xi \right|^3,
\end{align}
Summing up for all odd~$i$'s, we have
\begin{align}
\left| \sum_{ i : odd} h_i \left( \frac{\xi}{\sqrt{K}} \right) - \sum_{ i : odd} \frac{1}{2K} s_{i, 2} ^2 \xi^2 \right| \lesssim \frac{1}{\sqrt{K}} \left|\xi \right|^3. \label{e_taylor_variance}
\end{align}
Note that for odd~$i$'s,~\eqref{e_conditional_density_odd_variables},~\eqref{e_conditional leb for odd} and~\cite[Lemma 9]{Me11} imply~$s_{i,2}^2$ is uniformly bounded above and below. That is, there exists a uniform constant~$C_{\eqref{e_uniform_bdd_si2}}>0$ such that
\begin{align}
\frac{1}{C_{\eqref{e_uniform_bdd_si2}}} \leq s_{i,2}^2 \leq C_{\eqref{e_uniform_bdd_si2}}. \label{e_uniform_bdd_si2}
\end{align}
As a consequence, there is a constant~$C>0$ with
\begin{align}\label{e_uniform boundedness of sl2}
\frac{1}{C} \leq \sum_{ i : odd} \frac{1}{2K} s_{i, 2}^2 \leq C .
\end{align}
In particular for~$K$ large enough and~$\left|\frac{\xi}{\sqrt{K}}\right|\leq \delta$, the estimate~\eqref{e_taylor_variance} yields
\begin{align}
\text{Re} \left( \sum_{ i : odd} h_i \left( \frac{\xi}{\sqrt{K}}\right) \right) \geq \sum_{ i : odd} \frac{ s_{i, 2}^2 }{4K}\xi^2.
\end{align}
Furthermore, Lipschitz continuity of complex function~$y \mapsto \exp\left(y\right) \in \mathbb{C}$ on $\text{Re} y \leq - \sum_{ i : odd} \frac{ s_{i, 2}^2 }{4K}\xi^2$ yields
\begin{align}
&\left| \exp \left( - \sum_{i : odd} h_i \left( \frac{\xi}{\sqrt{K}} \right) \right)- \exp \left(- \sum_{ i : odd} \frac{ s_{i, 2}^2 }{2K} \xi^2 \right)\right| \\
& \qquad \qquad \qquad \qquad \qquad \qquad \qquad  \lesssim \frac{1}{\sqrt{K}} \left|\xi\right|^3 \exp \left( - \sum_{ i : odd} \frac{ s_{i, 2}^2 }{4K} \xi^2 \right). \label{changing exponential term}
\end{align}
Then a combination of~\eqref{e_proof for exchanging 0th der},~\eqref{e_uniform boundedness of sl2} and~\eqref{changing exponential term} implies the desired estimate~\eqref{e_lemma exchanging 0th der}. \\

Let us now address the estimate~\eqref{e_lemma exchanging 1st der}. The proof for this case is almost identical to~\eqref{e_lemma exchanging 0th der}. Let us denote~$\mathscr{F}_l ^2 : = \sigma \left( X_k , k \in E_1 ^l \text{ or } k : even \right)$. Then we have
\begin{align}
&\mathbb{E}\left[\exp\left( i \sum_{k \in E_2 ^l } \left( X_k -m_k \right) \frac{\xi}{\sqrt{K}}\right) \mid \mathscr{F}_l \right] \\
&\qquad =\mathbb{E}\left[\exp\left( i \sum_{\substack{k \in E_2 ^l \\ k : even} } \left( X_k -m_k \right) \frac{\xi}{\sqrt{K}} + i \sum_{\substack{i \in E_2 ^l\\ i : odd}} \left(\mathbb{E}\left[X_i \mid \mathscr{F}_l ^2 \right] -m_i \right) \frac{\xi}{\sqrt{K}}  \right)\right. \\
&\qquad  \quad \quad \quad  \quad \quad \times \left. \mathbb{E}\left[\exp\left( i \sum_{\substack{i \in E_2 ^l \\ i : odd} } \left( X_i -\mathbb{E}\left[X_i \mid \mathscr{F}_l ^2 \right] \right) \frac{\xi}{\sqrt{K}}\right) \mid \mathscr{F}_l ^2 \right] \mid \mathscr{F}_l \right].
\label{trick 2}
\end{align}
Let us also denote~$\tilde{s}_{i, 2}^2 := \mathbb{E}\left[ \left( X_i - \mathbb{E}\left[X_i \mid \mathscr{F}_l ^2 \right] \right)^2 \mid \mathscr{F}_l ^2 \right] $. Then one can easily prove the analogue of~\eqref{changing exponential term}. That is, for $\left| \frac{\xi}{\sqrt{K}} \right| \leq \delta$ and~$L \ll K$, it holds that
\begin{align} 
&\left| \mathbb{E}\left[\exp\left( i \sum_{\substack{i \in E_2 ^l \\ i : odd} } \left( X_i -\mathbb{E}\left[X_i \mid \mathscr{F}_l ^2 \right] \right) \frac{\xi}{\sqrt{K}}\right) \mid \mathscr{F}_l ^2 \right] - \exp \left(- \sum_{\substack{ i \in E_2 ^l \\ i : odd}} \frac{ \tilde{s}_{i, 2}^2 }{2K} \xi^2 \right)\right| \\
& \qquad  \lesssim \frac{1}{\sqrt{K}} \left|\xi\right|^3 \exp \left( - \sum_{ \substack{i \in E_2 ^l\\ i : odd}} \frac{ \tilde{s}_{i, 2}^2 }{4K} \xi^2 \right).  \label{analogue of changing exponential term}
\end{align}
Moreover, for~$L \ll K$ and~$K$ large enough, there exists a uniform constant~$C>0$ with
\begin{align} 
\frac{1}{C} \leq \sum_{ \substack{i \in E_2 ^l\\ i : odd}} \frac{ \tilde{s}_{i, 2}^2 }{4K} \leq C.
\end{align}
This implies, as desired,
\begin{align}
\left| \mathbb{E}\left[\exp\left( i \sum_{\substack{i \in E_2 ^l\\ i : odd }} \left( X_i -\mathbb{E}\left[X_i \mid \mathscr{F}_l ^2 \right] \right) \frac{\xi}{\sqrt{K}}\right) \mid \mathscr{F}_l ^2 \right] \right| &\lesssim \left(1+ \frac{\left|\xi\right|^3}{\sqrt{K}} \right)\exp \left( - \sum_{\substack{i \in E_2 ^l\\ i : odd }} \frac{ \tilde{s}_{i, 2}^2 }{4K} \xi^2 \right) \\
&\lesssim \left(1+ \frac{\left|\xi\right|^3}{\sqrt{K}} \right) \exp\left(-C \xi^2 \right). \label{proof for unknown2}
\end{align}
In particular with~\eqref{trick 2},
\begin{align}
\left|\mathbb{E}\left[\exp\left( i \sum_{k \in E_2^l } \left( X_k -m_k \right) \frac{\xi}{\sqrt{K}}\right) \mid \mathscr{F}_l \right] \right|&\lesssim \left( 1 + \frac{ \left|\xi\right|^3}{\sqrt{K}} \right)\exp\left(-C\xi^2 \right)  \lesssim \left(1+ \xi^2 \right) \exp\left(-C\xi^2 \right),
\end{align}
which proves~\eqref{e_lemma exchanging 1st der}. The inequality~\eqref{e_lemma exchanging 2nd der} is deduced by same type of argument.
\qed

\subsection{Proof of ~\eqref{0th derivative} in Proposition ~\ref{p_main computation}}\label{s_0th_derivative} \ An application of the inverse Fourier transform with~\eqref{e_m=1/Kmi} yields the representation
\begin{align}
2 \pi g_{K, m} (0) &= \int_{\mathbb{R}} \mathbb{E}\left[ \exp\left( i \frac{1}{\sqrt{K}}\sum_{k=1}^K \left( X_k - m \right) \xi \right) \right] d \xi \\
& \overset{\eqref{e_m=1/Kmi}}{=} \int_{\mathbb{R}} \mathbb{E}\left[ \exp\left( i \frac{1}{\sqrt{K}}\sum_{k=1}^K \left( X_k - m_k \right) \xi \right) \right] d \xi. \label{Inverse Fourier representation}
\end{align}
For~$\delta>0$, let us divide the integral~\eqref{Inverse Fourier representation} into two parts
\begin{align}
&\int_{\mathbb{R}} \mathbb{E}\left[ \exp\left( i \frac{1}{\sqrt{K}}\sum_{k=1}^K \left( X_k - m_k \right) \xi \right) \right] d \xi \\
&\qquad  = \int_{ \{ \left| \left(1/ \sqrt{K}\right)\xi\right|\leq\delta\}} \mathbb{E}\left[ \exp\left( i \frac{1}{\sqrt{K}}\sum_{k=1}^K \left( X_k - m_k \right) \xi \right) \right] d \xi \label{e_0th der (a)} \\
&\qquad  \quad + \int_{ \{ \left| \left(1/ \sqrt{K}\right)\xi\right|>\delta\}} \mathbb{E}\left[ \exp\left( i \frac{1}{\sqrt{K}}\sum_{k=1}^K \left( X_k - m_k \right) \xi \right) \right] d \xi  \label{e_0th der (b)}.
\end{align}
\medskip

\textbf{Argument for~\eqref{e_0th der (a)} : Estimation of the inner integral}. \ Recall the definitions~\eqref{e_def_of_ml2}, \eqref{e_def_si2} and~\eqref{e_def_of_hat_e} of~$m_{i,2}$,~$s_{i,2}^2$ and~$\hat{e}\left(\xi\right)$. Choosing~$\delta>0$ small enough and for~$K$ large enough,~\eqref{e_lemma exchanging 0th der} gives
\begin{align}
&\left| T_{\eqref{e_0th der (a)}} - \int_{ \{ \left| \left(1/ \sqrt{K}\right)\xi\right|\leq\delta\}}\mathbb{E} \left[ \hat{e}\left(\xi\right) \exp\left( - \sum_{i : odd} \frac{ s_{i, 2}^2}{2K}\xi^2 \right) \right] d\xi \right|\\
& \qquad  \lesssim \int_{ \{ \left| \left(1/ \sqrt{K}\right)\xi\right|\leq\delta\}} \frac{1}{ \sqrt{K}} \left| \xi \right|^3 \exp\left( - C\xi^2 \right) d\xi  \lesssim \frac{1}{\sqrt{K}}. \label{e_estimation of (a) in 0th der}
\end{align}
Note that we have by Fubini's theorem
\begin{align}
&\int_{ \{ \left| \left(1/ \sqrt{K}\right)\xi\right|\leq\delta\}}\mathbb{E} \left[ \hat{e}\left(\xi\right) \exp\left( - \sum_{i : odd} \frac{ s_{i, 2}^2}{2K}\xi^2 \right) \right] d\xi \\
&\qquad  = \mathbb{E}\left[ \int_{ \{ \left| \left(1/ \sqrt{K}\right)\xi\right|\leq\delta\}}\hat{e}\left(\xi\right)\exp\left( - \sum_{i : odd} \frac{ s_{i, 2}^2}{2K}\xi^2 \right)d\xi\right].
\label{exchanged inner integral}
\end{align}
We claim that~\eqref{exchanged inner integral} is uniformly bounded above and below by positive constants. To prove this, we divide the integral as follows

\begin{align}
&\mathbb{E}\left[ \int_{ \{ \left| \left(1/ \sqrt{K}\right)\xi\right|\leq\delta\}}\hat{e}\left(\xi\right)\exp\left( - \sum_{i : odd} \frac{ s_{i, 2}^2}{2K}\xi^2 \right)d\xi\right]\\
&\qquad  =\mathbb{E}\left[ \int_{ \mathbb{R}}\hat{e}\left(\xi\right)\exp\left( - \sum_{i : odd} \frac{ s_{i, 2}^2}{2K}\xi^2 \right)d\xi\right] \label{e_inner=whole-outer, whole}\\
&\qquad  \quad - \mathbb{E}\left[ \int_{ \{ \left| \left(1/ \sqrt{K}\right)\xi\right|>\delta\}}\hat{e}\left(\xi\right)\exp\left( - \sum_{i : odd} \frac{ s_{i, 2}^2}{2K}\xi^2 \right)d\xi\right]. \label{e_inner=whole-outer, outer}
\end{align}
Our aim is to prove that the whole integral given by~\eqref{e_inner=whole-outer, whole} is uniformly bounded above and below while the outer integral~\eqref{e_inner=whole-outer, outer} is relatively small. Let us begin with the estimation of~\eqref{e_inner=whole-outer, whole}. The upper bound of~\eqref{e_inner=whole-outer, whole} follows from~\eqref{e_uniform boundedness of sl2} that
\begin{align}
\left|T_{~\eqref{e_inner=whole-outer, whole}} \right| \leq \mathbb{E}\left[\int_{ \mathbb{R}} \exp\left( - \sum_{i : odd} \frac{ s_{i, 2}^2}{2K}\xi^2 \right)d\xi \right] \leq \mathbb{E}\left[\int_{ \mathbb{R}} \exp\left( - C \xi^2 \right)d\xi \right] \lesssim 1.
\end{align}
To deduce the lower bound of~\eqref{e_inner=whole-outer, whole}, we compute the integral inside the expectation directly. It holds that
\begin{align}
& \int_{ \mathbb{R}}\hat{e}\left(\xi\right)\exp\left( - \sum_{i : odd} \frac{ s_{i, 2}^2}{2K}\xi^2 \right)d\xi \\
&  =\int_{ \mathbb{R}} \exp\left( - \sum_{i : odd} \frac{s_{i,2}^2}{2K} \left(\xi- i \frac{ \sqrt{K}}{\sum_{i : odd} s_{i,2}^2 } \left( \sum_{k : even} \left(X_k -m_k \right)+ \sum_{i : odd} \left(m_{i,2} -m_i \right) \right) \right)^2 \right) \\
&  \quad \quad \times \exp\left( - \frac{ \left( \sum_{k : even} \left(X_k -m_k \right)+ \sum_{i : odd} \left(m_{i,2} -m_i \right)\right)^2 }{ 2\sum_{ i : odd} s_{i,2}^2 } \right) d\xi \\
&  = \int_{ \mathbb{R}} \exp\left( - \sum_{i : odd} \frac{s_{i,2}^2}{2K} \left(\xi- i \frac{ \sqrt{K}}{\sum_{i : odd} s_{i,2}^2 } \left( \sum_{k : even} \left(X_k -m_k \right)+ \sum_{i : odd} \left(m_{i,2} -m_i \right) \right) \right)^2 \right)d\xi \\
&  \quad \quad \times \exp\left( - \frac{ \left( \sum_{k : even} \left(X_k -m_k \right)+ \sum_{i : odd} \left(m_{i,2} -m_i \right)\right)^2 }{ 2\sum_{ i : odd} s_{i,2}^2 } \right). \label{e_complex contour integration}
\end{align}
Note that~$\frac{ \sqrt{K}}{\sum_{i : odd} s_{i,2}^2 } \left( \sum_{k : even} \left(X_k -m_k \right)+ \sum_{i : odd} \left(m_{i,2} -m_i \right) \right) \in \mathbb{R}$. Then complex contour integration implies
\begin{align}
T_{\eqref{e_complex contour integration}}& =  \int_{\mathbb{R} } \exp\left( - \sum_{ i : odd} \frac{ s_{i, 2}^2}{2K} \xi^2\right) d\xi \\
& \quad \quad \times \exp\left( - \frac{ \left( \sum_{k : even} \left(X_k -m_k \right)+ \sum_{i : odd} \left(m_{i,2} -m_i \right)\right)^2 }{ 2\sum_{ i : odd} s_{i,2}^2 } \right)\\
& =\sqrt{ \frac{2K \pi} { \sum_{ i : odd} s_{i,2}^2}} \exp\left( - \frac{ \left( \sum_{k : even} \left(X_k -m_k \right)+ \sum_{i : odd} \left(m_{i,2} -m_i \right)\right)^2 }{ 2\sum_{ i : odd} s_{i,2}^2 } \right)   \\
& \overset{\eqref{e_uniform boundedness of sl2}}{\gtrsim} \exp\left( - \frac{ \left( \sum_{k : even} \left(X_k -m_k \right)+ \sum_{i : odd} \left(m_{i,2} -m_i \right)\right)^2 }{ 2\sum_{ i : odd} s_{i,2}^2 } \right).  \label{whole integral estimate}
\end{align}
By taking the expectation and applying Jensen's inequality, it holds that
\begin{align}
&\mathbb{E}\left[ \exp\left( - \frac{ \left( \sum_{k : even} \left(X_k -m_k \right)+ \sum_{i : odd} \left(m_{i,2} -m_i \right)\right)^2 }{ 2\sum_{ i : odd} s_{i,2}^2 } \right)  \right] \\
&\qquad \geq \exp\left( - \mathbb{E}\left[  \frac{ \left( \sum_{k : even} \left(X_k -m_k \right)+ \sum_{i : odd} \left(m_{i,2} -m_i \right)\right)^2 }{ 2\sum_{ i : odd} s_{i,2}^2 }  \right]\right)  \\
&\qquad  \overset{\eqref{e_uniform boundedness of sl2}}{\geq} \exp\left( -  \mathbb{E}\left[ \frac{C}{K}  \left( \sum_{k : even} \left(X_k -m_k \right)+ \sum_{i : odd} \left(m_{i,2} -m_i \right)\right)^2  \right]\right) \\
&\qquad  \geq \exp\left( -  \mathbb{E}\left[ \frac{2C}{K}  \left( \sum_{k : even} \left(X_k -m_k \right) \right)^2+ \frac{2C}{K} \left(\sum_{i : odd} \left(m_{i,2} -m_l \right)\right)^2  \right]\right).
\end{align}
Now it follows from~\cite[Lemma 9]{Me11} and Lemma ~\ref{l_variance estimates} that as desired,
\begin{align}
&\exp\left( -  \mathbb{E}\left[ \frac{2C}{K}  \left( \sum_{k : even} \left(X_k -m_k \right) \right)^2+ \frac{2C}{K}\left(\sum_{i : odd} \left(m_{i,2} -m_i \right)\right)^2  \right]\right) \\
&\qquad = \exp\left( - \frac{2C}{K} \var\left(\sum_{k : even} X_k \right) -\frac{2C}{K} \var\left( \sum_{i : odd} m_{i, 2} \right) \right)\\
&\qquad = \exp\left(- \frac{2C}{K} \var\left(\sum_{k : even} X_k \right) - \frac{2C}{K} \sum_{i : odd} \var\left(m_{i,2}\right) \right) \geq C. \label{exchanged whole estimate}
\end{align}
Let us turn to the estimation of~\eqref{e_inner=whole-outer, outer}. Using~\eqref{e_uniform boundedness of sl2}, we have
\begin{align}
\left|T_{\eqref{e_inner=whole-outer, outer}} \right| 
&\leq \mathbb{E}\left[ \int_{ \{ \left| \left(1/ \sqrt{K}\right)\xi\right|>\delta\}}\exp\left( - \sum_{i : odd} \frac{ s_{i, 2}^2}{2K}\xi^2 \right)d\xi\right] \\
&\leq \mathbb{E}\left[ \int_{ \{ \left| \left(1/ \sqrt{K}\right)\xi\right|>\delta\}}\exp\left( - C\xi^2 \right)d\xi\right] \lesssim \frac{1}{\sqrt{K}}. \label{exchanged outer estimate}
\end{align}
Therefore a combination of ~\eqref{exchanged whole estimate} and ~\eqref{exchanged outer estimate} implies that there exists a positive constant~$C>0$ with
\begin{align}
\mathbb{E}\left[ \int_{ \{ \left| \left(1/ \sqrt{K}\right)\xi\right|\leq\delta\}}\hat{e}\left(\xi\right)\exp\left( - \sum_{i : odd} \frac{ s_{i, 2}^2}{2K}\xi^2 \right)d\xi\right] \in \left(\frac{1}{C}, C\right)
\end{align}
for~$K$ large enough. Combined with~\eqref{e_estimation of (a) in 0th der}, it holds that for~$K$ large enough,
\begin{align} \label{e_(a) estimate final}
\frac{1}{C} \leq T_{\eqref{e_0th der (a)}} \leq C .
\end{align}

\textbf{Argument for~\eqref{e_0th der (b)} : Estimation of the outer integral}. \ Let us recall the observations~\eqref{e_conditional_density_odd_variables} and~\eqref{e_conditional leb for odd}. We know the conditional random variable~$ X_i \mid X_j, j : even$ has the conditional Lebesgue density given by~\eqref{e_conditional leb for odd}. We observe that this conditional density has the same form as the one dimensional measure considered in~\cite[Lemma 3.4]{MeOt13} and~\cite[(47)]{MeOt13}. Changing~$i$ or the conditioned spins~$X_j$ only changes the linear term in the Hamiltonian of~\eqref{e_conditional leb for odd}. Because the estimates of~\cite[Lemma 3.4]{MeOt13} and~\cite[(47)]{MeOt13} are uniform in the linear term of the Hamiltonian, an application to the random variable $X_i \mid X_j, j : even$ and its distribution $\mu^{\sigma}\left(dx_i \mid x_j , j : even\right)$ yields the following: For any~$\hat \delta>0$ and odd~$i$'s, there exists~$\lambda<1$ with
\begin{align}
\left| \mathbb{E}\left[ \exp\left( i  \frac{ X_i -m_{i,2} }{s_{i,2}} \xi \right) \mid X_j, j : even \right] \right|  \leq \lambda \quad \text{for all } \left|\xi\right|\geq \hat \delta \label{e_outer less lambda_prev}
\end{align}
and
\begin{align}
\left| \mathbb{E}\left[ \exp\left( i  \frac{ X_i -m_{i,2} }{s_{i,2}} \xi \right) \mid X_j, j : even \right] \right|  \lesssim \left|\xi \right|^{-1}. \label{e_outer less xi_prev}
\end{align}
Recall that~$s_{i,2}^2$ is uniformly bounded from above and below (see~\eqref{e_uniform_bdd_si2}). By setting~$\xi = \frac{\hat{\xi}}{s_{i,2}}$ and~$\delta= \sqrt{C_{\eqref{e_uniform_bdd_si2}}}                                                                                                                                                                                                                                                                                                                                                                                                                                                                                                                                                                                                                         \hat \delta $ the estimates~\eqref{e_outer less lambda_prev} and~\eqref{e_outer less xi_prev} yield
\begin{align}
\left| \mathbb{E}\left[ \exp\left( i  \left( X_i -m_{i,2} \right) \xi \right) \mid X_j, j : even \right] \right|  \leq \lambda \quad \text{for all } \left|\xi\right| \geq \delta\label{e_outer less lambda}
\end{align}
and
\begin{align}
\left| \mathbb{E}\left[ \exp\left( i  \left( X_i -m_{i,2} \right) \xi \right) \mid X_j, j : even \right] \right|  \lesssim \frac{1}{s_{i,2} \left|\xi\right|} \overset{~\eqref{e_uniform_bdd_si2}}{\lesssim} \frac{1}{\left|\xi\right|} \label{e_outer less xi}.
\end{align}
Now consider the conditional expectation with respect to~$\{ X_j \mid j : even \}$:
\begin{align}
&\left|\int_{ \{ \left| \left(1/ \sqrt{K}\right)\xi\right|>\delta\}} \mathbb{E}\left[ \exp\left( i \frac{1}{\sqrt{K}}\sum_{k=1}^K \left( X_k - m_k \right) \xi \right) \right] d \xi \right| \\
&\qquad  =\left|\int_{ \{ \left| \left(1/ \sqrt{K}\right)\xi\right|>\delta\}} \mathbb{E}\left[ \exp\left( i \frac{1}{\sqrt{K}}\sum_{k : even} \left( X_k - m_k \right) \xi \right)\right.  \right. \\
& \quad \quad \quad \quad \times\left. \left. \mathbb{E}\left[  \exp\left( i \frac{1}{\sqrt{K}}\sum_{i : odd} \left( X_i - m_i \right) \xi \right) \mid X_j, j : even \right]\right] d \xi\right|  \\
& \qquad  \leq \int_{ \{ \left| \left(1/ \sqrt{K}\right)\xi\right|>\delta\}} \mathbb{E}\left[\left| \mathbb{E}\left[  \exp\left( i \frac{1}{\sqrt{K}}\sum_{i : odd} \left( X_i - m_i \right) \xi \right) \mid X_j, j : even \right]\right|\right] d \xi . \label{0th der outer estimate 1}
\end{align}
We apply~\eqref{e_outer less lambda} (on~$\frac{K}{2}-2$ of~$\frac{K}{2}$ factors) and~\eqref{e_outer less xi} (on the remaining 2 factors) to obtain
\begin{align}
\left|T_{\eqref{0th der outer estimate 1}} \right|
&   \lesssim \int_{ \{ \left| \left(1/ \sqrt{K}\right)\xi\right|>\delta\}} \lambda^{\frac{K}{2}-2} \left( \frac{1}{1+ \left(1/ \sqrt{K}\right)\left|\xi\right|}\right)^2 d\xi \\
&  \lesssim \int_{ \{ \left| \left(1/ \sqrt{K}\right)\xi\right|>\delta\}} K \lambda^{\frac{K}{2} -2 } \frac{1}{K + \xi^2} d\xi\\
&   \lesssim \int_{ \{ \left| \left(1/ \sqrt{K}\right)\xi\right|>\delta\}} K \lambda^{ \frac{K}{2} -2} \frac{1}{1+ \xi^2} d\xi\\
&   \leq  K \lambda^{ \frac{K}{2} -2} \int_{ \mathbb{R}} \frac{1}{1+ \xi^2} d\xi   \lesssim K \lambda^{\frac{K}{2}-2} \overset{\lambda<1}{\lesssim} \frac{1}{\sqrt{K}}.
\end{align}
Therefore we conclude that
\begin{align} \label{e_(b) estimate final}
\left|T_{\eqref{e_0th der (b)}}\right| \lesssim \frac{1}{\sqrt{K}}.
\end{align} 
Choosing~$K>0$ large enough,~\eqref{e_(a) estimate final} and~\eqref{e_(b) estimate final} imply the desired estimate~\eqref{0th derivative}.\qed

\subsection{Proof of~\eqref{2nd derivative} in Proposition~\ref{p_main computation}}\label{s_2nd_derivative} \ 
Let us address the inequality~\eqref{2nd derivative}. As before, we start with dividing the integral into inner and outer part.
\begin{align}
2\pi \frac{d^2}{d \sigma^2} g_{K, m} (0) & =\frac{d^2}{d \sigma^2}\int_{ \mathbb{R}} \mathbb{E}\left[ \exp\left( i \frac{1}{\sqrt{K}}\sum_{k=1}^K \left( X_k - m_k \right) \xi \right) \right] d \xi\\
&  = \frac{d^2}{d \sigma^2}\int_{ \{ \left| \left(1/ \sqrt{K}\right)\xi\right|\leq\delta\}} \mathbb{E}\left[ \exp\left( i \frac{1}{\sqrt{K}}\sum_{k=1}^K \left( X_k - m_k \right) \xi \right) \right] d \xi \label{2nd der inner + outer-inner}\\
&  \ + \frac{d^2}{d \sigma^2}\int_{ \{ \left| \left(1/ \sqrt{K}\right)\xi\right|>\delta\}} \mathbb{E}\left[ \exp\left( i \frac{1}{\sqrt{K}}\sum_{k=1}^K \left( X_k - m_k \right) \xi \right) \right] d \xi. \label{2nd der inner + outer-outer}
\end{align}
Recall the definition~\eqref{e_F1nl} and~\eqref{e_F2nl} of the auxiliary sets~$F_1 ^{n, l}$ and~$F_2 ^{n, l}$ (cf. Figure~\ref{f_sets_F}). We choose $L= K^{\varepsilon} \ll K$.\\

\textbf{Argument for~\eqref{2nd der inner + outer-inner} : Estimation of the inner integral}. We introduce a trick that allows to apply Taylor expansions for blocks. We consider the law of the random vector~$X$. Instead of considering for the every site~$i$ the homogeneous linear term~$\sigma x_i$, we artificially introduce the heterogeneous linear term~$\sigma_i x_i$. More precisely, it holds that
\begin{align}\label{e_introducing_sigma_i}
\mathbb{E} \left[ f(X) \right] &=  \int f(x) \frac{ \exp \left( \sum_{i=1}^K \sigma_i x_i -H(x) \right)}{\int \exp \left( \sum_{i=1}^K \sigma_i x_i -H(x) \right) dx } dx,
\end{align}
where we introduced on every site~$i$ the variable~$\sigma_i$ and set~$\sigma_i=\sigma$. Introducing the heterogeneous field~$\sigma_i$ gives more flexibility for Taylor expansions. We group sites together in blocks and then take advantage of the decay of correlations. Let us introduce an auxiliary notation
\begin{align} \label{e_def_yk}
Y_k : = X_k -m_k ,
\end{align}
and define a function~$G_{n,l} : \mathbb{R}^2 \rightarrow \mathbb{C}$ by
\begin{align}\label{e_def_gnl}
G_{n, l}(\xi_1, \xi_2) : = \mathbb{E} \left[\exp \left( i \sum_{i \in F_1 ^{n,l}} Y_i \xi_1 + i \sum_{j \in F_2 ^{n,l}} Y_j \xi_2  \right) \right] .
\end{align}
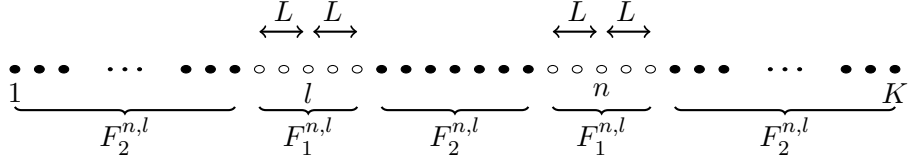
\begin{figure}[t]
\centering
\begin{tikzpicture}[xscale=1.3]

\draw[fill] (0,0) circle [radius=0.05];
\draw[fill] (.25,0) circle [radius=0.05];
\draw[fill] (.5,0) circle [radius=0.05];

\draw[fill] (0.975,0) circle [radius=0.02];
\draw[fill] (1.125,0) circle [radius=0.02];
\draw[fill] (1.275,0) circle [radius=0.02];

\draw[fill] (1.75,0) circle [radius=0.05];
\draw[fill] (2,0) circle [radius=0.05];
\draw[fill] (2.25,0) circle [radius=0.05];

\draw (2.5,0) circle[radius=0.05];
\draw (2.75,0) circle[radius=0.05];
\draw (3,0) circle[radius=0.05];
\draw (3.25,0) circle[radius=0.05];
\draw (3.5,0) circle[radius=0.05];

\draw[fill] (3.75,0) circle [radius=0.05];
\draw[fill] (4,0) circle [radius=0.05];
\draw[fill] (4.25,0) circle [radius=0.05];
\draw[fill] (4.5,0) circle [radius=0.05];
\draw[fill] (4.75,0) circle [radius=0.05];
\draw[fill] (5,0) circle [radius=0.05];
\draw[fill] (5.25,0) circle [radius=0.05];

\draw (5.5,0) circle[radius=0.05];
\draw (5.75,0) circle[radius=0.05];
\draw (6,0) circle[radius=0.05];
\draw (6.25,0) circle[radius=0.05];
\draw (6.5,0) circle[radius=0.05];

\draw[fill] (6.75,0) circle [radius=0.05];
\draw[fill] (7,0) circle [radius=0.05];
\draw[fill] (7.25,0) circle [radius=0.05];

\draw[fill] (7.725,0) circle [radius=0.02];
\draw[fill] (7.875,0) circle [radius=0.02];
\draw[fill] (8.025,0) circle [radius=0.02];

\draw[fill] (8.5,0) circle [radius=0.05];
\draw[fill] (8.75,0) circle [radius=0.05];
\draw[fill] (9,0) circle [radius=0.05];

\node[align=center, below] at (0,-.05) {$1$};
\node[align=center, below] at (3,-.05) {$l$};
\node[align=center, below] at (6,-.05) {$n$};
\node[align=center, below] at (9,-.05) {$K$};

\draw[decorate,decoration={brace,mirror},thick] (0,-.5) -- node[below]{$F_{2}^{n,l}$} (2.25,-.5);
\draw[decorate,decoration={brace,mirror},thick] (2.5,-.5) -- node[below]{$F_{1}^{n,l}$} (3.5,-.5);
\draw[decorate,decoration={brace,mirror},thick] (3.75,-.5) -- node[below]{$F_{2}^{n,l}$} (5.25,-.5);
\draw[decorate,decoration={brace,mirror},thick] (5.5,-.5) -- node[below]{$F_{1}^{n,l}$} (6.5,-.5);
\draw[decorate,decoration={brace,mirror},thick] (6.75,-.5) -- node[below]{$F_{2}^{n,l}$} (9,-.5);

\draw[thick,<->] (2.5,.5) -- (2.95,.5);
\node[align=center, above] at (2.75,.5) {$L$};
\draw[thick,<->] (3.05,.5) -- (3.5,.5);
\node[align=center, above] at (3.25,.5) {$L$};
\draw[thick,<->] (5.5,.5) -- (5.95,.5);
\node[align=center, above] at (5.75,.5) {$L$};
\draw[thick,<->] (6.05,.5) -- (6.5,.5);
\node[align=center, above] at (6.25,.5) {$L$};
\end{tikzpicture}
\caption{The sets $F_1^{n,l}$ and $F_2 ^{n,l}$ given by~\eqref{e_F1nl} and~\eqref{e_F2nl}}\label{f_sets_F}
\end{figure}
Note that
\begin{align}
&\frac{d^2}{d\sigma^2}\int_{ \{ \left| \left(1/ \sqrt{K}\right)\xi\right|\leq\delta\}} \mathbb{E}\left[ \exp\left( i \frac{1}{\sqrt{K}}\sum_{k=1}^K Y_k \xi \right) \right] d \xi  \\
&\qquad  = \frac{d^2}{d\sigma^2} \int_{ \{ \left| \left(1/ \sqrt{K}\right)\xi\right|\leq\delta\}}G_{n, l}\left(\frac{\xi}{\sqrt{K}}, \frac{\xi}{\sqrt{K}}\right) d \xi \\
&\qquad  =  \sum_{l=1}^{K} \sum_{n=1}^{K} \int_{ \{ \left| \left(1/ \sqrt{K}\right)\xi\right|\leq\delta\}}\frac{d}{d\sigma_n} \frac{d}{d\sigma_l} G_{n, l}\left(\frac{\xi}{\sqrt{K}}, \frac{\xi}{\sqrt{K}}\right) \Bigg|_{\sigma_l =\sigma_n  = \sigma} d \xi.
\end{align}
For convenience, we will write~$\frac{d}{d\sigma_n} f$ to denote~$\left.\frac{d}{d\sigma_n} f \right|_{\sigma_n =\sigma}$ for any index~$n$. For deducing a correct estimate for~\eqref{2nd der inner + outer-inner} it is sufficient to show that for given~$\beta>\frac{1}{2}$ and~$K$ large enough
\begin{align}\label{e_aim for 2nd inner}
\sum_{l=1}^K \sum_{n=1}^K \left|\int_{ \{ \left| \left(1/ \sqrt{K}\right)\xi\right|\leq\delta\}}\frac{d}{d\sigma_n} \frac{d}{d\sigma_l} G_{n, l}\left(\frac{\xi}{\sqrt{K}}, \frac{\xi}{\sqrt{K}}\right) d \xi\right| \lesssim K^{\beta}.
\end{align}
To establish this, we take the 2nd order Taylor expansion
\begin{align}
&\int_{\{ \left| \left(1/ \sqrt{K}\right)\xi\right|\leq\delta\}} \frac{d}{d \sigma_n }\frac{d}{d \sigma_l }  G_{n,l} \left(\frac{\xi}{\sqrt{K}}, \frac{\xi}{\sqrt{K}}\right)d\xi \\
& \qquad  = \int_{\{ \left| \left(1/ \sqrt{K}\right)\xi\right|\leq\delta\}} \frac{d}{d \sigma_n }\frac{d}{d \sigma_l }  G_{n,l} \left(0 , \frac{\xi}{\sqrt{K}}\right) d\xi \label{e_2nd derivative split form (a)}\\
& \qquad  \quad + \int_{\{ \left| \left(1/ \sqrt{K}\right)\xi\right|\leq\delta\}}\frac{d}{d \xi_1}  \frac{d}{d \sigma_n }\frac{d}{d \sigma_l } G_{n,l} \left(0, \frac{\xi}{\sqrt{K}}\right) \frac{\xi}{\sqrt{K}} d\xi \label{e_2nd derivative split form (b)}\\
&\qquad  \quad +\frac{1}{2}\int_{\{ \left| \left(1/ \sqrt{K}\right)\xi\right|\leq\delta\}}\frac{d^2}{d \xi_1^2}  \frac{d}{d \sigma_n }\frac{d}{d \sigma_l } G_{n,l} \left(\tilde{ \frac{\xi}{\sqrt{K}}}, \frac{\xi}{\sqrt{K}}\right) \left(\frac{\xi}{\sqrt{K}}\right)^2 d\xi \label{e_2nd derivative split form (c)},
\end{align}
where~$\frac{\tilde{\xi}}{\sqrt{K}}$ is a real number between~$0$ and~$\frac{\xi}{\sqrt{K}}$. In particular we have~$\left|\frac{\tilde{\xi}}{\sqrt{K}}\right| \leq \left|\frac{\xi}{\sqrt{K}}\right|$. The proof of~\eqref{e_aim for 2nd inner} is based on the following three lemmas. The first lemma is an estimation of~\eqref{e_2nd derivative split form (a)}. 
\begin{lemma}\label{l_proof_of_2nd_der_(a)}
Under the assumptions of Proposition~\ref{p_main computation}, it holds that for~$\frac{\left|\xi\right|}{\sqrt{K}} \leq \delta \leq 1$ and any~$n, l \in \{1, 2, \cdots K\}$,
\begin{align}
\left|\frac{d}{d \sigma_n }\frac{d}{d \sigma_l }  G_{n,l} \left(0 , \frac{\xi}{\sqrt{K}}\right) \right| \lesssim \left( 1 + \xi^2 \right)\exp \left(-CL\right).
\end{align}
\end{lemma}
The second lemma provides a bound on~\eqref{e_2nd derivative split form (b)}. 
\begin{lemma}\label{l_proof_of_2nd_der_(b)}
Under the assumptions of Proposition~\ref{p_main computation}, it holds that for~$\frac{\left|\xi\right|}{\sqrt{K}} \leq \delta \leq 1$ and any~$n, l \in \{1, 2, \cdots K\}$,
\begin{align}
\left|\frac{d}{d \xi_1}  \frac{d}{d \sigma_n }\frac{d}{d \sigma_l } G_{n,l} \left(0, \frac{\xi}{\sqrt{K}}\right) \right|  \lesssim  KL^2\left(1+ \xi^2 \right) \exp\left(-CL\right).
\end{align}
\end{lemma}
The last lemma is an estimation of~\eqref{e_2nd derivative split form (c)}. 
\begin{lemma}\label{l_proof_of_2nd_der_(c)}
Under the assumptions of Proposition~\ref{p_main computation}, it holds that for $\frac{\left|\tilde{\xi}\right|}{\sqrt{K}}\leq \frac{\left|\xi\right|}{\sqrt{K}} \leq \delta \leq 1$ and~$n,l$ with~$|n-l|>2L$,
\begin{align}
&\left|\frac{d^2}{d \xi_1^2}  \frac{d}{d \sigma_n }\frac{d}{d \sigma_l } G_{n,l} \left(\tilde{ \frac{\xi}{\sqrt{K}}}, \frac{\xi}{\sqrt{K}}\right)  \right|  \\
&\qquad \lesssim L^3 \left(1+ \xi^2 \right) \exp\left(-CL\right) +    L^4 \frac{\left|\xi\right|}{\sqrt{K}} \left(1+\xi^2\right) \exp\left(-C \xi^2 \right), \label{e_estimate_2nd_c_large}
\end{align}
and for~$n,l$ with~$|n-l|\leq 2L$,
\begin{align}
&\left|\frac{d^2}{d \xi_1^2}  \frac{d}{d \sigma_n }\frac{d}{d \sigma_l } G_{n,l} \left(\tilde{ \frac{\xi}{\sqrt{K}}}, \frac{\xi}{\sqrt{K}}\right)  \right|  \\
&\qquad \lesssim L^3 \left(1+ \xi^2 \right) \exp\left(-CL\right) +    L^4  \left(1+\xi^2\right) \exp\left(-C \xi^2 \right). \label{e_estimate_2nd_c_small}
\end{align}
\end{lemma}

We give the proof of the lemmas from above in Section~\ref{s_proof_of_aux_lemmas_inner}. We shall now see how it can be used to prove~\eqref{e_aim for 2nd inner} and therefore~\eqref{2nd der inner + outer-inner}. \\

\textbf{Estimation of~\eqref{e_aim for 2nd inner}:} \ Using the lemmas from above and recalling that~$\beta>\frac{1}{2}$ yield
\begin{align}
\left|T_{\eqref{e_2nd derivative split form (a)}} \right| \overset{Lemma~\ref{l_proof_of_2nd_der_(a)}}{\lesssim} \int_{ \{ \left| \left(1/ \sqrt{K}\right)\xi\right|\leq \delta\}} \left( 1+\xi^2 \right)\exp\left(-CL\right) d\xi  \overset{L=K^{\varepsilon}}{\lesssim} \frac{1}{K^{2-\beta}} \label{e_2nd_(a)_bound},
\end{align}
\begin{align}
\left|T_{\eqref{e_2nd derivative split form (b)}} \right| &\overset{Lemma~\ref{l_proof_of_2nd_der_(b)}}{\lesssim} \int_{ \{ \left| \left(1/ \sqrt{K}\right)\xi\right|\leq \delta\}}  KL^2\left(1+ \xi^2 \right) \exp\left(-CL\right)  \frac{\left|\xi \right|}{\sqrt{K}} d\xi \overset{L=K^{\varepsilon}}{\lesssim} \frac{1}{K^{2-\beta}}, \label{e_2nd_(b)_bound}
\end{align}
Summing over all pairs~$(n,l)$, it holds that
\begin{align}
\sum_{n=1}^K \sum_{l=1}^K \left| \int_{\{ \left| \left(1/ \sqrt{K}\right)\xi\right|\leq\delta\}} \frac{d}{d \sigma_n }\frac{d}{d \sigma_l }  G_{n,l} \left(0 , \frac{\xi}{\sqrt{K}}\right) d\xi\right|  \lesssim K^{\beta}, \\
\sum_{n=1}^K \sum_{l=1}^K \left| \int_{\{ \left| \left(1/ \sqrt{K}\right)\xi\right|\leq\delta\}}\frac{d}{d \xi_1}  \frac{d}{d \sigma_n }\frac{d}{d \sigma_l } G_{n,l} \left(0, \frac{\xi}{\sqrt{K}}\right) \frac{\xi}{\sqrt{K}} d\xi\right| \lesssim K^{\beta}.
\end{align}
Let us consider~\eqref{e_2nd derivative split form (c)}. For a pair~$\left(n, l \right)$ with~$|n-l|> 2L$, it holds that
\begin{align}
\left|T_{\eqref{e_2nd derivative split form (c)}} \right| &\overset{\eqref{e_estimate_2nd_c_large}}{\lesssim} \int_{ \{ \left| \left(1/ \sqrt{K}\right)\xi\right|\leq \delta\}} L^3 \left(1+ \xi^2 \right) \exp\left(-CL\right) \left( \frac{\left|\xi \right|}{\sqrt{K}} \right)^2 d\xi \\
& \quad \quad + \int_{ \{ \left| \left(1/ \sqrt{K}\right)\xi\right|\leq \delta\}} L^4 \frac{\left|\xi\right|}{\sqrt{K}} \left(1+\xi^2\right) \exp\left(-C \xi^2 \right) \left( \frac{\left|\xi \right|}{\sqrt{K}} \right)^2 d\xi \\
& \quad \lesssim L^3 \left(\sqrt{K}+K^{3/2}\right) \exp\left(-CL\right) + \frac{L^4}{K^{3/2}} . \label{e_2nd_(c)_bound1}
\end{align}
For~$\left(n, l\right)$ with~$|n-l|\leq 2L$,
\begin{align}
\left|T_{\eqref{e_2nd derivative split form (c)}} \right| &\overset{\eqref{e_estimate_2nd_c_small}}{\lesssim} \int_{ \{ \left| \left(1/ \sqrt{K}\right)\xi\right|\leq \delta\}} L^3 \left(1+ \xi^2 \right) \exp\left(-CL\right) \left( \frac{\left|\xi \right|}{\sqrt{K}} \right)^2 d\xi \\
&  \quad + \int_{ \{ \left| \left(1/ \sqrt{K}\right)\xi\right|\leq \delta\}} L^4  \left(1+\xi^2\right) \exp\left(-C \xi^2 \right) \left( \frac{\left|\xi \right|}{\sqrt{K}} \right)^2 d\xi \\
& \lesssim L^3 \left(\sqrt{K}+K^{3/2}\right) \exp\left(-CL\right) + \frac{L^4}{K}. \label{e_2nd_(c)_bound2}
\end{align}
A combination of~\eqref{e_2nd_(c)_bound1} and~\eqref{e_2nd_(c)_bound2} yields
\begin{align}
&\sum_{n=1}^K \sum_{l=1}^K \left|\int_{\{ \left| \left(1/ \sqrt{K}\right)\xi\right|\leq\delta\}}\frac{d^2}{d \xi_1^2}  \frac{d}{d \sigma_n }\frac{d}{d \sigma_l } G_{n,l} \left(\tilde{ \frac{\xi}{\sqrt{K}}}, \frac{\xi}{\sqrt{K}}\right) \left(\frac{\xi}{\sqrt{K}}\right)^2 d\xi \right| \\
& \qquad = \sum_{n=1}^K \sum_{l :  |n-l|>2L} \left|\int_{\{ \left| \left(1/ \sqrt{K}\right)\xi\right|\leq\delta\}}\frac{d^2}{d \xi_1^2}  \frac{d}{d \sigma_n }\frac{d}{d \sigma_l } G_{n,l} \left(\tilde{ \frac{\xi}{\sqrt{K}}}, \frac{\xi}{\sqrt{K}}\right) \left(\frac{\xi}{\sqrt{K}}\right)^2 d\xi \right| \\
& \qquad \quad + \sum_{n=1}^K \sum_{l :  |n-l|\leq 2L} \left|\int_{\{ \left| \left(1/ \sqrt{K}\right)\xi\right|\leq\delta\}}\frac{d^2}{d \xi_1^2}  \frac{d}{d \sigma_n }\frac{d}{d \sigma_l } G_{n,l} \left(\tilde{ \frac{\xi}{\sqrt{K}}}, \frac{\xi}{\sqrt{K}}\right) \left(\frac{\xi}{\sqrt{K}}\right)^2 d\xi \right| \\
& \qquad \lesssim K^2 \left(L^3 \left(\sqrt{K}+K^{3/2}\right) \exp\left(-CL\right) + \frac{L^4}{K^{3/2}} \right) + KL \left(L^3 \left(\sqrt{K}+K^{3/2}\right) \exp\left(-CL\right) + \frac{L^4}{K} \right) \\
& \quad \quad \overset{L=K^{\varepsilon}}{\lesssim} K^{\beta}.
\end{align}
This gives the desired estimate~\eqref{e_aim for 2nd inner}.
\qed  

\medskip

\textbf{Argument for~\eqref{2nd der inner + outer-outer} : Estimation of the outer integral}. \ Let us start with providing auxiliary ingredients. The first lemma provides auxiliary estimates: 
\begin{lemma} \label{l_proof_dervatives_outer}
Recall the definition~\eqref{e_def_of_ml2} of~$m_{i,2}$. It holds that
\begin{align}
\left|\frac{d}{d\sigma} \mathbb{E}\left[ \exp\left( i  \left( X_i -m_{i,2} \right) \xi \right) \mid X_j, j : even \right] \right| \lesssim 1+ \left|\xi\right|,  \label{e_outer first der} \\
\left|\frac{d^2}{d\sigma ^2} \mathbb{E}\left[ \exp\left( i  \left(X_i -m_{i,2} \right) \xi \right) \mid X_j, j : even \right] \right| \lesssim 1+ \left|\xi\right|^2. \label{e_outer second der}
\end{align}
\end{lemma}

\textsc{Proof of Lemma~\ref{l_proof_dervatives_outer}}. \ 
Let us begin with providing an auxiliary computation. Recalling~\eqref{e_conditional leb for odd}, one can easily prove that (cf. see also~\eqref{e_def_derivative_si2},~\eqref{e_def_ti2} and~\eqref{e_derivative_formula})
\begin{align}
&\frac{d}{d\sigma} \mathbb{E} \left[ f(X_i) \mid X_j, j : even \right] \\
&\qquad = \mathbb{E} \left[ \left( X_i -m_{i,2} \right) f(X_i) \mid X_j, j : even \right] + \mathbb{E} \left[ \frac{d}{d\sigma}f(X_i) \mid X_j, j : even  \right]. \label{e_derivative_formula_for_conditional_exp}
\end{align}
This implies in particular
\begin{align}
&\frac{d}{d\sigma} \mathbb{E}\left[ \exp\left( i  \left( X_i -m_{i,2} \right) \xi \right) \mid X_j, j : even \right]  \\
& \qquad \overset{\eqref{e_def_derivative_si2}}{=}  \mathbb{E} \left[ \left( X_i -m_{i,2} \right)  \exp\left( i  \left( X_i -m_{i,2} \right) \xi \right) \mid X_j, j : even   \right] \label{e_outer_first_der_aa}  \\
& \qquad \quad + \mathbb{E} \left[ -i \xi s_{i,2}^2 \exp\left( i  \left( X_i -m_{i,2} \right) \xi \right) \mid X_j, j : even  \right]. \label{e_outer_first_der_ab}
\end{align}
Note that~\eqref{e_analogue of proving prop lemma} and~\eqref{e_uniform_bdd_si2} yields
\begin{align}
\left| T_{\eqref{e_outer_first_der_aa}}\right| = \left| \mathbb{E} \left[ \left( X_i -m_{i,2} \right)  \exp\left( i  \left( X_i -m_{i,2} \right) \xi \right) \mid X_j, j : even   \right] \right| \overset{\eqref{e_analogue of proving prop lemma}}{ \lesssim} 1, \label{e_outer_first_der_a} \\
\left| T_{\eqref{e_outer_first_der_ab}} \right| =\left| \mathbb{E} \left[ -i \xi s_{i,2}^2 \exp\left( i  \left( X_i -m_{i,2} \right) \xi \right) \mid X_j, j : even  \right]\right| \overset{\eqref{e_uniform_bdd_si2}}{\lesssim} \left|\xi\right|. \label{e_outer_first_der_b}
\end{align}
As desired, the estimates~\eqref{e_outer_first_der_a} and~\eqref{e_outer_first_der_b} yield
\begin{align}
\left| \frac{d}{d\sigma} \mathbb{E}\left[ \exp\left( i  \left( X_i -m_{i,2} \right) \xi \right) \mid X_j, j : even \right]  \right| \leq \left|T_{\eqref{e_outer_first_der_aa}} \right| + \left|T_{\eqref{e_outer_first_der_ab}} \right| \lesssim 1 + \left| \xi \right|.
\end{align}

Let us turn to the derivation of~\eqref{e_outer second der}. The calculation from above yields that
\begin{align}
&\frac{d^2}{d\sigma^2} \mathbb{E}\left[ \exp\left( i  \left( X_i -m_{i,2} \right) \xi \right) \mid X_j, j : even \right]  \\
& \qquad {=}  \frac{d}{d\sigma}\mathbb{E} \left[ \left( X_i -m_{i,2} \right)  \exp\left( i  \left( X_i -m_{i,2} \right) \xi \right) \mid X_j, j : even   \right]  \label{e_outer_second_der_aa} \\
& \qquad \quad + \frac{d}{d\sigma} \mathbb{E} \left[ -i \xi s_{i,2}^2 \exp\left( i  \left( X_i -m_{i,2} \right) \xi \right) \mid X_j, j : even  \right]. \label{e_outer_second_der_ab}
\end{align}
Let us consider~\eqref{e_outer_second_der_aa}: A direct computation using~\eqref{e_derivative_formula_for_conditional_exp} and~\eqref{e_def_derivative_si2} implies
\begin{align}
& T_{\eqref{e_outer_second_der_aa}} = \mathbb{E} \left[ \left( X_i -m_{i,2} \right)^2  \exp\left( i  \left( X_i -m_{i,2} \right) \xi \right) \mid X_j, j : even \right] \label{e_outer_second_der_a1} \\
& \qquad \quad + \mathbb{E} \left[ -s_{i,2}^2  \exp\left( i  \left( X_i -m_{i,2} \right) \xi \right) \mid X_j, j : even \right] \label{e_outer_second_der_a2}\\
& \qquad \quad + \mathbb{E} \left[ \left( X_i -m_{i,2} \right) \left(-i \xi s_{i,2}^2 \right) \exp\left( i  \left( X_i -m_{i,2} \right) \xi \right) \mid X_j, j : even  \right]. \label{e_outer_second_der_a3}
\end{align}
The terms~\eqref{e_outer_second_der_a1},~\eqref{e_outer_second_der_a2} and~\eqref{e_outer_second_der_a3} can be estimated similarly. Indeed, it holds that
\begin{align}
\left| T_{\eqref{e_outer_second_der_a1}}\right|, \left| T_{\eqref{e_outer_second_der_a2}}\right| \lesssim 1 \qquad \text{ and } \qquad
\left| T_{\eqref{e_outer_second_der_a3}}\right| \lesssim \left|\xi\right|.
\end{align}
As a consequence, we have
\begin{align}
\left|T_{\eqref{e_outer_second_der_aa}}\right| &\leq \left|T_{\eqref{e_outer_second_der_a1}} \right|+ \left|T_{\eqref{e_outer_second_der_a2}} \right|+ \left|T_{\eqref{e_outer_second_der_a3}} \right| \lesssim 1+ \left|\xi\right| \lesssim 1+ \left|\xi\right|^2.
\end{align}
Let us consider~\eqref{e_outer_second_der_ab}: Recall the definition~\eqref{e_def_ti2} of~$t_{i,2}$. Then it holds from~\eqref{e_derivative_formula_for_conditional_exp} that
\begin{align}
T_{\eqref{e_outer_second_der_ab}} &= \mathbb{E} \left[ \left( X_i -m_{i,2} \right) \left(-i \xi s_{i,2}^2 \right) \exp\left( i  \left( X_i -m_{i,2} \right) \xi \right) \mid X_j, j : even  \right] \label{e_outer_second_der_b1}  \\
&\quad + \mathbb{E} \left[ \left(-i \xi t_{i,2} \right) \exp\left( i  \left( X_i -m_{i,2} \right) \xi \right) \mid X_j, j : even  \right] \label{e_outer_second_der_b2} \\
& \quad + \mathbb{E} \left[ \left(-i \xi s_{i,2}^2 \right)^2 \exp\left( i  \left( X_i -m_{i,2} \right) \xi \right) \mid X_j, j : even  \right]. \label{e_outer_second_der_b3}
\end{align}
Now~\eqref{e_analogue of proving prop lemma} and~\eqref{e_uniform_bdd_si2} imply
\begin{align}
\left|T_{\eqref{e_outer_second_der_b1}} \right|, \left|T_{\eqref{e_outer_second_der_b2}} \right| \lesssim \left| \xi \right| \qquad \text{ and } \qquad 
\left|T_{\eqref{e_outer_second_der_b3}} \right| \lesssim \left|\xi \right|^2.
\end{align}
Therefore we obtain
\begin{align}
\left|T_{\eqref{e_outer_second_der_ab}} \right| \leq \left|T_{\eqref{e_outer_second_der_b1}} \right|+ \left|T_{\eqref{e_outer_second_der_b2}} \right|+ \left|T_{\eqref{e_outer_second_der_b3}} \right| \lesssim  \left| \xi \right| + \left| \xi \right|^2 \lesssim 1+ \left|\xi\right|^2.
\end{align}
To conclude, we sum up the estimates from above. As desired, we have
\begin{align}
\left| \frac{d^2}{d\sigma^2} \mathbb{E}\left[ \exp\left( i  \left( X_i -m_{i,2} \right) \xi \right) \mid X_j, j : even \right]  \right| 
& \leq \left|T_{\eqref{e_outer_second_der_aa}} \right|+\left|T_{\eqref{e_outer_second_der_ab}} \right|  \lesssim 1+ \left|\xi\right|^2.
\end{align}
\qed
\medskip

The second lemma is an auxiliary computation:
\begin{lemma} \label{l_derivatives_of_hat_e}
Recall the definition~\eqref{e_def_of_hat_e} and~\eqref{e_def_ti2} of~$\hat{e}$ and~$t_{i,2}$, respectively. It holds that
\begin{align} \label{e_1st_derivative_hat_e}
\frac{d}{d\sigma} \hat{e}(\xi) &= \hat{e}(\xi) \left( - i \frac{1}{\sqrt{K}}\xi  \mathbb{E} \left[ \left(\sum_{n=1}^K\left(X_n -m_n \right) \right)^2 \right] + i \frac{1}{\sqrt{K}} \xi \sum_{i : odd} s_{i,2}^2 \right), \\
\frac{d^2}{d\sigma ^2} \hat{e}(\xi) &=\hat{e}(\xi) \left( - i \frac{1}{\sqrt{K}}\xi \mathbb{E} \left[ \left(\sum_{n=1}^K\left(X_n -m_n \right) \right) ^2 \right] + i \frac{1}{\sqrt{K}} \xi \sum_{i : odd} s_{i,2}^2 \right)^2 \\
& \quad + \hat{e}\left(\xi\right) \left(-i \frac{1}{\sqrt{K}} \xi  \mathbb{E} \left[\left( \sum_{n=1}^{K} \left(X_n -m_n \right) \right)^3 \right] +  i \frac{1}{\sqrt{K}} \xi \sum_{i : odd} t_{i,2} \right). \label{e_2nd_derivative_hat_e}
\end{align}
\end{lemma}

\medskip

\textsc{Proof of Lemma~\ref{l_derivatives_of_hat_e}}. \ 
We start with deducing~\eqref{e_1st_derivative_hat_e}. Let us recall the equations~\eqref{e_derivative_of_mk} and~\eqref{e_def_derivative_si2}. It follows from a direct computation that
\begin{align}
\frac{d}{d\sigma} \hat{e}(\xi) & =  \left( -i \frac{1}{\sqrt{K}}\xi \sum_{k : even} \mathbb{E} \left[\left( \sum_{n=1}^K \left(X_n -m_n \right) \right) X_k   \right]  \right. \\
& \qquad \quad \left. + i \frac{1}{\sqrt{K}}\xi \sum_{i : odd} \left(s_{i,2}^2 - \mathbb{E} \left[\left( \sum_{n=1}^K \left(X_n -m_n \right) \right) X_k   \right] \right) \right) \\
& \quad \times \exp \left( i \frac{1}{\sqrt{K}}\sum_{k : even} \left( X_k - m_k \right) \xi + i \frac{1}{\sqrt{K}} \sum_{i : odd} \left( m_{i, 2} -m_i \right) \xi\right) \\
& = \left( -i \frac{1}{\sqrt{K}}\xi \sum_{k=1}^{K} \mathbb{E} \left[\left( \sum_{n=1}^K \left(X_n -m_n \right) \right) X_k   \right] + i \frac{1}{\sqrt{K}}\xi \sum_{i : odd} s_{i,2}^2 \right) \hat{e}(\xi) \\
& = \left( - i \frac{1}{\sqrt{K}} \xi \mathbb{E} \left[ \left(\sum_{n=1}^K\left(X_n -m_n \right) \right)^2 \right] + i \frac{1}{\sqrt{K}} \xi \sum_{i : odd} s_{i,2}^2 \right) \hat{e}(\xi).
\end{align}
The equation~\eqref{e_2nd_derivative_hat_e} can be derived with a similar argument using~\eqref{e_def_ti2}. \qed 

\medskip

Now let us move on to the estimation of the outer integral~\eqref{2nd der inner + outer-outer}. Recalling~\eqref{e_derivative_formula}, we get that the integrand in~\eqref{2nd der inner + outer-outer} is
\begin{align} \label{second derivative outer integral split form}
&\frac{d^2}{d \sigma^2} \mathbb{E}\left[ \exp\left(i \frac{1}{\sqrt{K}} \sum_{k=1}^{K} \left(X_k -m_k \right) \xi \right) \right]\\
&\ =\frac{d^2}{d \sigma^2} \mathbb{E} \left[ \hat{e}\left(\xi\right) \prod_{ i : odd} \mathbb{E} \left[\exp\left( i \frac{1}{\sqrt{K}}  \left( X_i -m_{i,2} \right) \xi \right) \mid X_j, j : even \right] \right]  \\
&\  = \frac{d}{d\sigma}\mathbb{E} \left[ \left(\sum_{n=1}^K \left(X_n -m_n \right) \right) \hat{e}\left(\xi\right)\prod_{ i : odd} \mathbb{E} \left[\exp\left( i \frac{1}{\sqrt{K}}  \left( X_i -m_{i,2} \right) \xi \right) \mid X_j, j : even \right] \right] \label{e_sec outer 1}\\
&\ \quad + \frac{d}{d\sigma}\mathbb{E} \left[ \frac{d}{d\sigma} \left(\hat{e}\left(\xi\right) \right)  \prod_{ i : odd} \mathbb{E} \left[\exp\left( i \frac{1}{\sqrt{K}}  \left( X_i -m_{i,2} \right) \xi \right) \mid X_j, j : even \right] \right] \label{e_sec outer 2}\\
&\  \quad +\frac{d}{d\sigma}\mathbb{E} \left[ \hat{e}\left(\xi\right) \frac{d}{d\sigma} \left( \prod_{ i : odd} \mathbb{E} \left[\exp\left( i \frac{1}{\sqrt{K}}  \left( X_i -m_{l,2} \right) \xi \right) \mid X_j, j : even \right] \right) \right].\label{e_sec outer 3}
\end{align}
We estimate each term on the right hand side separately. We begin with estimating~\eqref{e_sec outer 1}. Applying~\eqref{e_derivative_formula}, we get
\begin{align}
T_{\eqref{e_sec outer 1}}
& = \mathbb{E} \left[ \left(\sum_{n=1}^K \left(X_n -m_n \right) \right)^2 \hat{e}\left(\xi\right) \prod_{ i : odd} \mathbb{E} \left[\exp\left( i \frac{1}{\sqrt{K}}  \left( X_i -m_{i,2} \right) \xi \right) \mid X_j, j : even \right] \right]  \label{e_sec outer 11}\\
&\quad + \mathbb{E} \left[ \frac{d}{d\sigma}\left(\sum_{n=1}^K \left(X_n -m_n \right) \right) \hat{e}\left(\xi\right) \right. \\
&\quad \quad \quad \quad \quad \quad \left. \times \prod_{ i : odd} \mathbb{E} \left[\exp\left( i \frac{1}{\sqrt{K}}  \left( X_i -m_{i,2} \right) \xi \right) \mid X_j, j : even \right] \right] \label{e_sec outer 12}\\
& \quad + \mathbb{E} \left[ \left(\sum_{n=1}^K \left(X_n -m_n \right) \right) \frac{d}{d\sigma} \left(\hat{e}\left(\xi\right) \right)\right. \\
& \quad \quad \quad \quad \quad  \quad \left. \times \prod_{ i : odd} \mathbb{E} \left[\exp\left( i \frac{1}{\sqrt{K}}  \left( X_i -m_{i,2} \right) \xi \right) \mid X_j, j : even \right] \right] \label{e_sec outer 13}\\
& \quad + \mathbb{E} \left[\left(\sum_{n=1}^K \left(X_n -m_n \right) \right) \hat{e}\left(\xi\right) \right.\\
&\quad \quad \quad \quad \quad \quad \left. \times \frac{d}{d\sigma} \left( \prod_{ i : odd} \mathbb{E} \left[\exp\left( i \frac{1}{\sqrt{K}}  \left( X_i -m_{i,2} \right) \xi \right) \mid X_j, j : even \right] \right) \right]. \label{e_sec outer 14}
\end{align}
Let us consider~\eqref{e_sec outer 11}: Applying~\eqref{e_outer less lambda} to~$\frac{K}{2}-2$ many factors and~\eqref{e_outer less xi} on the remaining 2 factors leads to
\begin{align}
\left|T_{\eqref{e_sec outer 11}}\right| &\leq \mathbb{E}\left[\left(\sum_{n=1}^K \left(X_n -m_n \right) \right)^2 \lambda^{\frac{K}{2}-2} \left( \frac{1}{ 1+ (1/\sqrt{K})\left|\xi\right| } \right)^2 \right].
\end{align}
Then Lemma~\ref{c_amgm} implies
\begin{align}
\left|T_{\eqref{e_sec outer 11}}\right| &\lesssim K^2 \lambda^{\frac{K}{2}-2} \left(\frac{1}{ 1+ (1/\sqrt{K})\left|\xi\right| } \right)^2 \lesssim K^3 \lambda^{\frac{K}{2}-2} \frac{1}{K + \xi^2} \lesssim K^3 \lambda^{\frac{K}{2}-2} \frac{1}{1+\xi^2}.
\end{align}
Let us consider~\eqref{e_sec outer 12}: A similar estimation using~\eqref{e_derivative_of_mk} yields that
\begin{align}
\left|T_{\eqref{e_sec outer 12}}\right| & = \left| \mathbb{E} \left[- \mathbb{E}\left[\left(\sum_{n=1}^K \left(X_n -m_n \right) \right)^2\right] \hat{e}\left(\xi\right)\right.\right.\\
&\quad \quad \quad \quad \quad \quad \left.\left.\times \vphantom{\left(\sum_{n=1}^K \left(X_n -m_n \right) \right)^2} \prod_{ i : odd} \mathbb{E} \left[\exp\left( i \frac{1}{\sqrt{K}}  \left( X_i -m_{i,2} \right) \xi \right) \mid X_j, j : even \right] \right] \right| \\
& \lesssim  K^3 \lambda^{\frac{K}{2}-2} \frac{1}{ 1+ \xi^2 }.
\end{align}
Let us consider~\eqref{e_sec outer 13}: We first expand the term using~\eqref{e_1st_derivative_hat_e}. Then applying~\eqref{e_outer less lambda} to~$\frac{K}{2}-3$ many factors and~\eqref{e_outer less xi} to~$3$ factors yields
\begin{align}
\left|T_{\eqref{e_sec outer 13}}\right|& \overset{\eqref{e_1st_derivative_hat_e}}{=} \left| \mathbb{E} \left[ \left(\sum_{n=1}^K \left(X_n -m_n \right) \right)\hat{e}\left(\xi\right) \right.\right.\\
& \quad \quad \quad \times \left( -i \frac{1}{\sqrt{K}} \xi \mathbb{E} \left[\left( \sum_{n=1}^{K} \left(X_n -m_n \right) \right)^2 \right] +  i \frac{1}{\sqrt{K}} \xi \sum_{i : odd} s_{i,2}^2 \right) \\
&\quad \quad \quad \quad \quad \quad \quad \left.\left. \times \prod_{ i : odd} \mathbb{E} \left[\exp\left( i \frac{1}{\sqrt{K}}  \left( X_i -m_{i,2} \right) \xi \right) \mid X_j, j : even \right] \right] \right| \\
& \overset{\eqref{e_amgm}, \eqref{e_uniform boundedness of sl2}, \eqref{e_outer less lambda}, \eqref{e_outer less xi}}{\lesssim} K \left(\frac{\left|\xi\right|}{\sqrt{K}}K^2 +\frac{\left|\xi\right|}{\sqrt{K}}K\right) \lambda^{\frac{K}{2}-3} \left( \frac{1}{ 1+ (1/\sqrt{K})\left|\xi\right| } \right)^3\\
& \qquad  \ \lesssim K^4 \lambda^{\frac{K}{2}-3} \left|\xi\right| \left(\frac{1}{\sqrt{K} + \left|\xi\right|}\right)^3 \lesssim K^4 \lambda^{\frac{K}{2}-3} \frac{1}{1+\xi^2 }.
\end{align}
Let us consider~\eqref{e_sec outer 14}: A combination of~\eqref{e_amgm},~\eqref{e_outer less lambda},~\eqref{e_outer less xi} and~\eqref{e_outer first der} yields
\begin{align}
\left|T_{\eqref{e_sec outer 14}}\right|& \lesssim  \frac{K}{2} \cdot K \left( 1+ \left| \frac{\xi}{\sqrt{K}} \right| \right)   \lambda^{\frac{K}{2} - 4} \left( \frac{1}{1 + \left(1/\sqrt{K}\right) \left|\xi\right| } \right)^3  \lesssim  K^3 \lambda^{ \frac{K}{2} -4} \frac{1}{ 1+ \xi^2}.
\end{align}
Hence we have
\begin{align}
\left|T_{\eqref{e_sec outer 1}}\right| & \leq \left|T_{\eqref{e_sec outer 11}}\right|+\left|T_{\eqref{e_sec outer 12}}\right|+\left|T_{\eqref{e_sec outer 13}}\right|+\left|T_{\eqref{e_sec outer 14}} \right| \\
& \lesssim K^3 \lambda^{\frac{K}{2}-2} \frac{1}{ 1+ \xi^2 }+ K^4 \lambda^{\frac{K}{2}-3} \frac{1}{1+\xi^2 }+ K^3 \lambda^{ \frac{K}{2} -4} \frac{1}{ 1+ \xi^2}  \lesssim K^4 \lambda^{ \frac{K}{2} -4} \frac{1}{ 1+ \xi^2} .
\end{align}
Let us turn to the estimation of~\eqref{e_sec outer 2}. A direct computation using~\eqref{e_derivative_formula} yields
\begin{align}
T_{\eqref{e_sec outer 2}} &= \mathbb{E} \left[\left(\sum_{n=1}^K \left(X_n -m_n \right) \right) \frac{d}{d\sigma} \left(\hat{e}\left(\xi\right) \right)\right. \\
& \quad \quad \quad \quad \left. \times  \prod_{ i : odd} \mathbb{E} \left[\exp\left( i \frac{1}{\sqrt{K}}  \left( X_i -m_{i,2} \right) \xi \right) \mid X_j, j : even \right] \right] \label{e_sec outer 21}\\
& \quad + \mathbb{E} \left[ \frac{d^2}{d\sigma ^2} \left(\hat{e}\left(\xi\right) \right) \right.\\
& \quad \quad \quad \quad \left. \times  \prod_{ i : odd} \mathbb{E} \left[\exp\left( i \frac{1}{\sqrt{K}}  \left( X_i -m_{i,2} \right) \xi \right) \mid X_j, j : even \right] \right] \label{e_sec outer 22}\\
& \quad + \mathbb{E} \left[ \frac{d}{d\sigma} \left(\hat{e}\left(\xi\right) \right)\right. \\
& \quad \quad \quad \quad \left. \times  \frac{d}{d\sigma}\left(\prod_{ i : odd} \mathbb{E} \left[\exp\left( i \frac{1}{\sqrt{K}}  \left( X_i -m_{i,2} \right) \xi \right) \mid X_j, j : even \right]\right) \right]. \label{e_sec outer 23}
\end{align}
Let us estimate each term on right hand side step by step. We start with estimating~\eqref{e_sec outer 21}. It holds that
\begin{align}
\left|T_{\eqref{e_sec outer 21}}\right|= \left|T_{\eqref{e_sec outer 13}}\right| \lesssim K^4 \lambda^{\frac{K}{2}-3} \frac{1}{1+\xi^2 }.
\end{align}
Let us consider~\eqref{e_sec outer 22}: An application of Lemma~\ref{c_amgm} to~\eqref{e_2nd_derivative_hat_e} yields
\begin{align}
\left|\frac{d^2}{d\sigma^2} \left(\hat{e}\left(\xi\right) \right)\right| &\overset{\eqref{e_amgm}}{\lesssim} \left(\frac{\left|\xi\right|}{\sqrt{K}}K^2 + \frac{\left|\xi\right|}{\sqrt{K}}K \right)^2  + \left( \frac{\left|\xi\right|}{\sqrt{K}}K^3 + \frac{\left|\xi\right|}{\sqrt{K}}K \right)\\
& \lesssim K^3  \xi^2 + K^{5/2}\left|\xi\right| \lesssim K^3 \left(1+ \xi^2 \right). \label{e_d2sigma2 hate}
\end{align}
Inserting~\eqref{e_d2sigma2 hate} into~\eqref{e_sec outer 22} gives
\begin{align}
\left|T_{\eqref{e_sec outer 22}}\right| 
& \lesssim K^3 \left(1+ \xi^2 \right)\mathbb{E} \left[ \prod_{ i : odd}\left| \mathbb{E} \left[\exp\left( i \frac{1}{\sqrt{K}}  \left( X_i -m_{i,2} \right) \xi \right) \mid X_j, j : even \right] \right|\right]\\
& \overset{\eqref{e_outer less lambda}, \eqref{e_outer less xi} }{\lesssim K^3} \left(1+ \xi^2 \right) \lambda^{\frac{K}{2} - 4}\left( \frac{1}{ 1+ (1/\sqrt{K})\left|\xi\right| } \right)^4   \lesssim K^5 \lambda^{\frac{K}{2} - 4} \frac{1}{1+\xi^2}.
\end{align}
Let us consider~\eqref{e_sec outer 23}: Using~\eqref{e_1st_derivative_hat_e},~\eqref{e_outer less lambda},~\eqref{e_outer less xi} and~\eqref{e_outer first der} yields
\begin{align}
\left|T_{\eqref{e_sec outer 23}}\right| & \lesssim \frac{K}{2} \cdot \left(\frac{\left|\xi\right|}{\sqrt{K}}K^2 +\frac{\left|\xi\right|}{\sqrt{K}}K\right) \left( 1+ \left| \frac{\xi}{\sqrt{K}} \right| \right)   \lambda^{\frac{K}{2} - 5} \left( \frac{1}{1 + \left(1/\sqrt{K}\right) \left|\xi\right| } \right)^4 \\
&\lesssim K^4 \lambda^{\frac{K}{2} - 5}  \frac{1}{1 + \xi^2} .
\end{align}
Therefore it holds that
\begin{align}
\left|T_{\eqref{e_sec outer 2}} \right| &\leq \left|T_{\eqref{e_sec outer 21}}\right|+\left|T_{\eqref{e_sec outer 22}}\right|+\left|T_{\eqref{e_sec outer 23}}\right|\\
& \lesssim K^4 \lambda^{\frac{K}{2}-3} \frac{1}{1+\xi^2 }+ K^5 \lambda^{\frac{K}{2} - 4} \frac{1}{1+\xi^2} +  K^4 \lambda^{\frac{K}{2} - 5}  \frac{1}{1 + \xi^2}  \lesssim K^5 \lambda^{\frac{K}{2} - 5}  \frac{1}{1 + \xi^2}.
\end{align}
Let us turn to the estimation of~\eqref{e_sec outer 3}, which is deduced by the same type of argument. More precisely, we get
\begin{align}
T_{\eqref{e_sec outer 3}} &= \mathbb{E} \left[\left(\sum_{n=1}^{K} \left(X_n -m_n \right) \right) \hat{e}\left(\xi\right)\right. \\
& \quad \quad \quad \left.\times \frac{d}{d\sigma} \left( \prod_{ i : odd} \mathbb{E} \left[\exp\left( i \frac{1}{\sqrt{K}}  \left( X_i -m_{i,2} \right) \xi \right) \mid X_j, j : even \right] \right) \right] \label{e_sec outer 31}\\
& \quad + \mathbb{E} \left[ \frac{d}{d\sigma}\left(\hat{e}\left(\xi\right)\right) \frac{d}{d\sigma} \left( \prod_{ i : odd} \mathbb{E} \left[\exp\left( i \frac{1}{\sqrt{K}}  \left( X_i -m_{i,2} \right) \xi \right) \mid X_j, j : even \right] \right) \right] \label{e_sec outer 32} \\
& \quad + \mathbb{E} \left[ \hat{e}\left(\xi\right) \frac{d^2}{d\sigma ^2} \left( \prod_{ i : odd} \mathbb{E} \left[\exp\left( i \frac{1}{\sqrt{K}}  \left( X_i -m_{i,2} \right) \xi \right) \mid X_j, j : even \right] \right) \right]. \label{e_sec outer 33}
\end{align}
Note that~$T_{\eqref{e_sec outer 31}}=T_{\eqref{e_sec outer 14}}$ and~$T_{\eqref{e_sec outer 32}}=T_{\eqref{e_sec outer 23}}$.
Let us estimate~\eqref{e_sec outer 33}. We have by using~\eqref{e_outer first der},~\eqref{e_outer second der},~\eqref{e_outer less lambda} and~\eqref{e_outer less xi} that
\begin{align}
&\left|\frac{d^2}{d\sigma ^2} \left( \prod_{ i : odd} \mathbb{E} \left[\exp\left( i \frac{1}{\sqrt{K}}  \left( X_i -m_{i,2} \right) \xi \right) \mid X_j, j : even \right] \right) \right| \\
& \qquad \lesssim \frac{K}{2} \cdot \left( 1+ \left| \frac{\xi}{\sqrt{K}} \right|^2 \right)  \lambda^{\frac{K}{2} - 5} \left( \frac{1}{1 + \left(1/\sqrt{K}\right) \left|\xi\right| } \right)^4\\
& \qquad \quad + \frac{K}{2}\left(\frac{K}{2}-1\right)\left( 1+ \left| \frac{\xi}{\sqrt{K}} \right| \right) ^2 \lambda^{\frac{K}{2} - 6} \left( \frac{1}{1 + \left(1/\sqrt{K}\right) \left|\xi\right| } \right)^{4} \\
& \qquad \lesssim K^2 \lambda^{\frac{K}{2} - 5} \frac{1}{1+ \xi^2}+K^3 \lambda^{\frac{K}{2} - 6} \frac{1}{1+ \xi^2}.
\end{align}
Thus we conclude
\begin{align}
\left|T_{\eqref{e_sec outer 33}}\right| \lesssim K^3 \lambda^{\frac{K}{2} - 6}\frac{1}{1+ \xi^2}.
\end{align}
It follows that
\begin{align}
\left|T_{\eqref{e_sec outer 3}} \right| & \leq \left|T_{\eqref{e_sec outer 31}}\right|+\left|T_{\eqref{e_sec outer 32}}\right|+\left|T_{\eqref{e_sec outer 33}}\right|  \\
& \lesssim K^3 \lambda^{ \frac{K}{2} -4} \frac{1}{ 1+ \xi^2} + K^4 \lambda^{\frac{K}{2} - 5}  \frac{1}{1 + \xi^2} + K^3 \lambda^{\frac{K}{2} - 6}\frac{1}{1+ \xi^2} \lesssim K^4 \lambda^{ \frac{K}{2} -6}\frac{1}{1+ \xi^2} .
\end{align}
\medskip
Overall, we have proven that
\begin{align}
&\left|\frac{d^2}{d \sigma^2} \mathbb{E}\left[ \exp\left(i \frac{1}{\sqrt{K}} \sum_{k=1}^{K} \left(X_k -m_k \right) \xi \right) \right]\right|
\\
& \qquad \leq \left|T_{\eqref{e_sec outer 1}}\right| +\left|T_{\eqref{e_sec outer 2}}\right| +\left|T_{\eqref{e_sec outer 3}} \right| \\
& \qquad \lesssim K^4 \lambda^{ \frac{K}{2} -4} \frac{1}{ 1+ \xi^2} + K^5 \lambda^{\frac{K}{2} - 5}  \frac{1}{1 + \xi^2} + K^4 \lambda^{ \frac{K}{2} -6}\frac{1}{1+ \xi^2}\\
& \qquad \lesssim K^5 \lambda^{\frac{K}{2}-6} \frac{1}{1+ \xi^2} \overset{\lambda <1}{\lesssim} K^{\beta} \frac{1}{1+ \xi^2}.
\end{align}
This implies, as desired,
\begin{align}
&\left|\frac{d^2}{d \sigma^2}\int_{\{ \left| \left(1/ \sqrt{K}\right)\xi\right|>\delta\}} \mathbb{E}\left[ \exp\left(i \frac{1}{\sqrt{K}} \sum_{k=1}^{K} \left(X_k -m_k \right) \xi \right) \right]d\xi\right| \\
& \qquad = \left|\int_{\{ \left| \left(1/ \sqrt{K}\right)\xi\right|>\delta\}}\frac{d^2}{d \sigma^2} \mathbb{E}\left[ \exp\left(i \frac{1}{\sqrt{K}} \sum_{k=1}^{K} \left(X_k -m_k \right) \xi \right) \right]d\xi\right|\\
&\qquad \lesssim K^{\beta} \int_{\{ \left| \left(1/ \sqrt{K}\right)\xi\right|>\delta\}} \frac{1}{1+ \xi^2} d\xi \\
&\qquad \lesssim K^{\beta}.
\end{align}
Combined with the argument for~\eqref{2nd der inner + outer-inner}, we obtain~\eqref{2nd derivative} in Proposition~\ref{p_main computation}.
\qed

\subsection{Proof of Lemma~\ref{l_proof_of_2nd_der_(a)}, Lemma~\ref{l_proof_of_2nd_der_(b)} and Lemma~\ref{l_proof_of_2nd_der_(c)}.}\label{s_proof_of_aux_lemmas_inner}

We shall only present the proof of Lemma~\ref{l_proof_of_2nd_der_(a)} and Lemma~\ref{l_proof_of_2nd_der_(c)}. One can derive Lemma~\ref{l_proof_of_2nd_der_(b)} by using similar arguments. The main ingredients are Lemma~\ref{c_amgm}, Lemma~\ref{l_exponential decay of correlations} and Lemma~\ref{l_exchanging exponential terms}. Additionally, we will need the following auxiliary observation. Recalling that~$Y_l = X_l -\mathbb{E}\left[ X_l \right]$, a direct calculation using~\eqref{e_introducing_sigma_i} yields
\begin{align}\label{e_calculation_d_dsigma_f}
\frac{d}{d\sigma_l} \mathbb{E} \left[ f(X) \right] & = \mathbb{E} \left[Y_l f(X) + \frac{d}{d\sigma_l} f(X)\right].
\end{align}
Hence, we have by additionally using that by definition~$\mathbb{E} \left[ Y_k \right] = 0$ for any~$k$
\begin{align} \label{e_par_der_Y_k}
 \frac{d}{d\sigma_l} Y_k = - \mathbb{E}\left[ Y_l X_k \right] = -\mathbb{E} \left[ Y_l Y_k \right], \ \mbox{and} \
\frac{d}{d\sigma_n} \mathbb{E} \left[ Y_l Y_k \right] = \mathbb{E} \left[ Y_n Y_l Y_k \right].
\end{align}
The last two formulas will be used many times in the upcoming arguments.
\medskip

\textsc{Proof of Lemma~\ref{l_proof_of_2nd_der_(a)}}. \ Let~$e(\xi_1, \xi_2)$ indicate
\begin{align}
e(\xi_1, \xi_2) := \exp \left( i \sum_{i \in F_1 ^{n,l}} Y_i \xi_1 + i \sum_{j \in F_2 ^{n,l}} Y_j \xi_2  \right).
\end{align}
Recall the definition~\eqref{e_def_gnl} of~$G_{n,l}\left(\xi_1, \xi_2\right)$. A direct calculation using~\eqref{e_calculation_d_dsigma_f} yields 
\newpage
\begin{align}
\frac{d}{d \sigma_n }\frac{d}{d \sigma_l } G_{n,l} \left( \xi_1 , \xi_2 \right)   & =\frac{d}{d \sigma_n }\mathbb{E}\left[ Y_l e(\xi_1, \xi_2) \right] \label{e_dsigman a1}  \\
&\quad +  \frac{d}{d \sigma_n }\mathbb{E} \left[ -i \xi_1 \mathbb{E}\left[ Y_l  \sum_{ i \in F_1 ^{n,l} } Y_i  \right] e(\xi_1, \xi_2) \right] \label{e_dsigman a2}\\
&  \quad + \frac{d}{d \sigma_n }\mathbb{E} \left[ -i \xi_2 \mathbb{E} \left[ Y_l  \sum_{j \in F_2 ^{n,l} } Y_j \right] e(\xi_1, \xi_2) \right]. \label{e_dsigman a3}
\end{align}
Let us consider~\eqref{e_dsigman a1}. A combination of~\eqref{e_calculation_d_dsigma_f} and~\eqref{e_par_der_Y_k} yields
\begin{align}
\frac{d}{d \sigma_n }\mathbb{E}\left[ Y_l e(\xi_1, \xi_2) \right] &  =\mathbb{E}\left[Y_n Y_l e(\xi_1, \xi_2) \right] \label{e_second der a11}  \\
& \quad + \mathbb{E} \left[ -\mathbb{E} \left[ Y_n Y_l \right] e(\xi_1, \xi_2) \right]\label{e_second der a12} \\
&  \quad + \mathbb{E} \left[ Y_l \left(-i \xi_1\right) \mathbb{E}\left[ Y_n \sum_{i \in F_1 ^{n, l} }  Y_i \right] e(\xi_1, \xi_2) \right] \label{e_second der a13}\\
& \quad + \mathbb{E} \left[ Y_l \left(-i \xi_2\right) \mathbb{E}\left[ Y_n\sum_{j \in F_2 ^{n, l} }  Y_j \right] e(\xi_1, \xi_2) \right] \label{e_second der a14}.
\end{align}
Setting~$(\xi_1 , \xi_2 )=\left(0, \frac{\xi}{\sqrt{K}}\right)$ and using Lemma~\ref{c_amgm} and Lemma~\ref{l_exponential decay of correlations}, we obtain
\begin{align}
& \left|T_{\eqref{e_second der a11}}+T_{\eqref{e_second der a12}}\right|  = \left| \cov \left( Y_n Y_l , e\left(0, \frac{\xi}{\sqrt{K}} \right) \right) \right| \lesssim \left|\xi\right|\exp\left(-CL\right), \\
& T_{\eqref{e_second der a13}} = 0, \qquad \mbox{and} \\
& \left|T_{\eqref{e_second der a14}} \right| \leq \left| \frac{\xi}{\sqrt{K}} \right| \left|\cov \left(Y_n, \sum_{j \in F_2 ^{n, l} }  Y_j \right) \right| \mathbb{E}\left|Y_l \right| \lesssim  \left|\xi\right| \exp\left(-CL\right). 
\end{align}
Therefore we conclude that for~$(\xi_1 , \xi_2 )=\left(0, \frac{\xi}{\sqrt{K}}\right)$,
\begin{align}
\left|T_{\eqref{e_dsigman a1}}  \right| &= \left| T_{\eqref{e_second der a11}}+T_{\eqref{e_second der a12}}+T_{\eqref{e_second der a13}}+T_{\eqref{e_second der a14}}\right| \\
&\leq \left|T_{\eqref{e_second der a11}}+T_{\eqref{e_second der a12}} \right|+\left|T_{\eqref{e_second der a13}}\right| +  \left|T_{\eqref{e_second der a14}}\right| \\
&\lesssim \left|\xi\right|\exp\left(-CL\right) \lesssim \left(1+ \xi^2 \right)\exp\left(-CL\right). \label{e_second der a1 estimate}
\end{align}
We observe that when~$(\xi_1, \xi_2 ) = \left(0, \frac{\xi}{\sqrt{K}}\right)$ we have
\begin{align}
T_{\eqref{e_dsigman a2}} = 0. \label{e_second der a2 estimate}
\end{align}
Let us turn to the estimation of~\eqref{e_dsigman a3}. A direct computation using~\eqref{e_calculation_d_dsigma_f} and~\eqref{e_par_der_Y_k} yields
\begin{align}
&\frac{d}{d \sigma_n }\mathbb{E} \left[ -i \xi_2 \mathbb{E} \left[ Y_l  \sum_{j \in F_2  ^{n,l}} Y_j \right] e\left(\xi_1, \xi_2\right) \right]\\
&\qquad   = -i \xi_2 \frac{d}{d \sigma_n } \left( \mathbb{E} \left[ Y_l \left( \sum_{j \in F_2 ^{n,l}} Y_j \right)\right]\right) \mathbb{E} \left[ e\left(\xi_1, \xi_2\right) \right]  \\
&\qquad   \quad -i \xi_2  \mathbb{E} \left[ Y_l  \sum_{j \in F_2 ^{n,l}} Y_j \right]\frac{d}{d \sigma_n } \left(\mathbb{E} \left[ e\left(\xi_1, \xi_2\right) \right] \right)  \\
& \qquad = - i \xi_2 \mathbb{E} \left[ Y_n Y_l  \sum_{j \in F_2 ^{n,l}} Y_j \right]\mathbb{E}\left[e\left(\xi_1, \xi_2\right) \right] \label{e_second der a31}\\
& \qquad \quad  - i \xi_2 \mathbb{E} \left[ Y_l \sum_{j \in F_2 ^{n,l}} Y_j \right]\mathbb{E}\left[ Y_n e\left(\xi_1, \xi_2\right) \right] \label{e_second der a32}\\
& \qquad \quad - i \xi_2 \mathbb{E} \left[ Y_l \sum_{j \in F_2 ^{n,l}} Y_j \right] \mathbb{E}\left[ \left(-i \xi_1 \right)\mathbb{E} \left[ Y_n \sum_{i \in F_1 ^{n,l}} Y_i \right]e\left(\xi_1, \xi_2\right) \right] \label{e_second der a33}\\
& \qquad \quad - i \xi_2 \mathbb{E} \left[ Y_l \sum_{j_1 \in F_2 ^{n,l}} Y_{j_1} \right] \mathbb{E}\left[ \left(-i \xi_2 \right)\mathbb{E} \left[ Y_n \sum_{j_2 \in F_2 ^{n,l}} Y_{j_2} \right]e\left(\xi_1, \xi_2\right) \right] \label{e_second der a34}.
\end{align}
Then plugging in~$(\xi_1, \xi_2 ) = \left(0, \frac{\xi}{\sqrt{K}}\right)$ and applying Lemma~\ref{c_amgm} and Lemma~\ref{l_exponential decay of correlations} yield
\begin{align}
\left|T_{\eqref{e_second der a31}}\right| &\leq \left| \frac{\xi}{\sqrt{K}} \right| \left| \cov\left( Y_n Y_l , \sum_{j \in F_2 ^{n,l} } Y_j \right)\right|   \lesssim \left|\xi\right| \exp\left(-CL\right),
\end{align}
\begin{align}
\left|T_{\eqref{e_second der a32}}\right| \leq \left| \frac{\xi}{\sqrt{K}} \right| \left|\cov\left( Y_l , \sum_{j \in F_2 ^{n,l}} Y_j \right) \right| \mathbb{E} \left|Y_n \right| \lesssim \left|\xi\right|\exp\left(-CL\right),
\end{align}
\begin{align}
T_{\eqref{e_second der a33}}=0, \qquad \mbox{and} 
\end{align}
\begin{align}
\left|T_{\eqref{e_second der a34}}\right|& \leq \left(\frac{\xi}{\sqrt{K}}\right)^2 \left|\cov\left(Y_l , \sum_{j_1 \in F_2 ^{n,l}} Y_{j_1} \right)\right| \left|\cov\left(Y_n , \sum_{j_2 \in F_2 ^{n,l}} Y_{j_2} \right)\right| \lesssim \xi^2 \exp\left(-2CL\right).
\end{align}
Thus a combination of the bounds from above yields
\begin{align}
\left| T_{\eqref{e_dsigman a3}} \right| &= \left|T_{\eqref{e_second der a31}}+T_{\eqref{e_second der a32}}+T_{\eqref{e_second der a33}}+T_{\eqref{e_second der a34}} \right| \\
& \leq \left|T_{\eqref{e_second der a31}}\right|+\left|T_{\eqref{e_second der a32}}\right|+\left|T_{\eqref{e_second der a33}}\right|+\left|T_{\eqref{e_second der a34}} \right| \\
&\lesssim \left|\xi\right| \exp\left(-CL\right) + \xi^2 \exp\left(-2CL\right)\\
&\lesssim \left(  \left|\xi \right| + \xi^2 \right) \exp\left(-CL\right) \\
& \lesssim \left(1+ \xi^2 \right)\exp\left(-CL\right). \label{e_second der a3 estimate}
\end{align}
We then sum up the bounds of~\eqref{e_dsigman a1},~\eqref{e_dsigman a2} and~\eqref{e_dsigman a3} given by~\eqref{e_second der a1 estimate},~\eqref{e_second der a2 estimate} and~\eqref{e_second der a3 estimate} respectively, and this finishes the proof of Lemma~\ref{l_proof_of_2nd_der_(a)}. \qed 
\medskip

The proof of Lemma~\ref{l_proof_of_2nd_der_(c)} is similar, but it needs more careful estimation. \\

\textsc{Proof of Lemma~\ref{l_proof_of_2nd_der_(c)}}. \ A direct computation using~\eqref{e_calculation_d_dsigma_f} and~\eqref{e_par_der_Y_k} yields
\begin{align}
&\frac{d^2}{d \xi_1^2} \frac{d}{d \sigma_n}\frac{d}{d \sigma_l} G_{n,l} \left(\xi_1, \xi_2\right) \\
&\qquad  = \frac{d}{d \sigma_n}\frac{d}{d \sigma_l}\frac{d^2}{d \xi_1^2} G_{n,l} \left(\xi_1, \xi_2\right) \\
& \qquad  = \frac{d}{d \sigma_n}\frac{d}{d \sigma_l} \mathbb{E} \left[ \left( i \sum_{i \in F_1 ^{n,l}} Y_i \right)^2 e\left(\xi_1, \xi_2\right) \right] \\
& \qquad  = \frac{d}{d \sigma_n}\mathbb{E} \left[ Y_l \left( i \sum_{i \in F_1 ^{n,l}} Y_i  \right)^2 e\left(\xi_1, \xi_2\right) \right]\label{e_dsigman c1} \\
& \qquad  \quad + \frac{d}{d \sigma_n}\mathbb{E} \left[ 2 \left( i \sum_{ i_1 \in F_1 ^{n,l}} Y_{i_1}  \right) \left( - i \right) \mathbb{E}\left[ Y_l \sum_{ i_2 \in F_1 ^{n,l}} Y_{i_2} \right] e\left(\xi_1, \xi_2\right) \right]\label{e_dsigman c2} \\
&\qquad   \quad + \frac{d}{d \sigma_n}\mathbb{E} \left[ \left( i \sum_{ i_1 \in F_1 ^{n,l}} Y_{i_1}   \right)^2 \left( -i\xi_1 \right) \mathbb{E} \left[Y_l \sum_{ i_2 \in F_1 ^{n,l}} Y_{i_2} \right] e\left(\xi_1, \xi_2\right) \right]\label{e_dsigman c3}
\\
&\qquad   \quad + \frac{d}{d \sigma_n}\mathbb{E} \left[ \left( i \sum_{ i \in F_1 ^{n,l}} Y_i  \right)^2 \left( -i  \xi_2 \right) \mathbb{E} \left[Y_l \sum_{ j \in F_2 ^{n,l}} Y_j \right] e\left(\xi_1, \xi_2\right) \right].\label{e_dsigman c4} 
\end{align}
Estimating the terms~\eqref{e_dsigman c3} and~\eqref{e_dsigman c4} is a lot easier than estimating~\eqref{e_dsigman c1} and~\eqref{e_dsigman c2}. The argument consists of a straightforward calculation and application of Lemma~\ref{c_amgm}, Lemma~\ref{l_exponential decay of correlations} and Lemma~\ref{l_exchanging exponential terms}. More precisely, for any~$n, l \in \{1,2, \cdots ,K\}$ and~$\left(\xi_1, \xi_2 \right) = \left( \frac{\tilde{\xi}}{\sqrt{K}}, \frac{\xi}{\sqrt{K}}\right)$  such that~$\frac{\left|\tilde{\xi}\right|}{\sqrt{K}} \leq \frac{\left|\xi\right|}{\sqrt{K}} \leq \delta \leq 1$, it holds that
\begin{align}
\left|T_{~\eqref{e_dsigman c3}} \right| &\lesssim L^3 \left(1+ \xi^2\right) \exp\left(-CL\right) + L^4 \frac{ \left|\xi\right|}{\sqrt{K}} \left( 1+ \xi^2 \right) \exp\left(-C \xi^2 \right),
\end{align}
and
\begin{align}
\left|T_{\eqref{e_dsigman c4}}\right| & \lesssim L^3 \left(1+ \xi^2 \right) \exp\left(-CL\right) .
\end{align}
We leave the details to the reader and turn to the more subtle estimation of~\eqref{e_dsigman c1} and~\eqref{e_dsigman c2}. The argument is more evolved because both terms have to be considered together. Only then, one sees covariances and can take advantage of the decay of correlations (cf.~Lemma~\ref{l_exponential decay of correlations}). Let us now turn to the details. Let us begin with expanding~\eqref{e_dsigman c1} and~\eqref{e_dsigman c2}. A direct computation using~\eqref{e_calculation_d_dsigma_f} and~\eqref{e_par_der_Y_k} yields that~\eqref{e_dsigman c1} is
\begin{align}
&\frac{d}{d\sigma_n}\mathbb{E} \left[ Y_l \left( i \sum_{i \in F_1 ^{n,l}}Y_i  \right)^2 e\left(\xi_1, \xi_2\right) \right]\\
&\qquad  =- \mathbb{E} \left[ Y_n Y_l \left(  \sum_{ i \in F_1 ^{n,l}} Y_i  \right)^2 e\left(\xi_1, \xi_2\right) \right] \label{e_second der c11}\\
&\qquad  \quad -\mathbb{E} \left[ -\mathbb{E} \left[ Y_n Y_l \right] \left(  \sum_{ i \in F_1 ^{n,l}} Y_i  \right)^2 e\left(\xi_1, \xi_2\right) \right] \label{e_second der c12}\\
&\qquad   \quad - \mathbb{E} \left[ Y_l 2\left(  \sum_{ i_1 \in F_1 ^{n,l}} Y_{i_1}  \right) \left(-\mathbb{E}\left[ Y_n \sum_{ i_2 \in F_1 ^{n,l}} Y_{i_2} \right] \right) e\left(\xi_1, \xi_2\right) \right] \label{e_second der c13}\\
&\qquad   \quad - \mathbb{E} \left[ Y_l \left( \sum_{ i_1 \in F_1 ^{n,l}} Y_{i_1}  \right)^2 \left( -i \xi_1 \right)  \mathbb{E}\left[Y_n \sum_{i_2 \in F_1 ^{n,l}} Y_{i_2} \right] e\left(\xi_1, \xi_2\right) \right] \label{e_second der c14}\\
&\qquad  \quad - \mathbb{E} \left[ Y_l \left( \sum_{ i \in F_1 ^{n,l}} Y_i  \right)^2 \left( -i \xi_2 \right)  \mathbb{E}\left[ Y_n \sum_{j \in F_2 ^{n,l}} Y_j \right] e\left(\xi_1, \xi_2\right) \right].  \label{e_second der c15} 
\end{align}
and~\eqref{e_dsigman c2} is 
\begin{align}
&\frac{d}{d \sigma_n} \mathbb{E} \left[ 2 \left( i \sum_{ i_1 \in F_1 ^{n,l}}Y_{i_1} \right)\left( - i \right) \mathbb{E}\left[ Y_l \sum_{ i_2 \in F_1 ^{n,l}} Y_{i_2} \right] e\left(\xi_1, \xi_2\right) \right] \\
&\qquad  = 2\frac{d}{d\sigma_n} \left( \mathbb{E}\left[ Y_l \sum_{ i_1 \in F_1 ^{n,l}} Y_{i_1} \right] \right)\mathbb{E} \left[  \left( \sum_{ i_2 \in F_1 ^{n,l}} Y_{i_2}  \right)  e\left(\xi_1, \xi_2\right) \right]  \\
& \qquad  \quad +  2\mathbb{E}\left[ Y_l \sum_{ i_1 \in F_1 ^{n,l}} Y_{i_1} \right] \frac{d}{d\sigma_n} \left( \mathbb{E} \left[  \left( \sum_{ i_2 \in F_1 ^{n,l}} Y_{i_2}  \right)  e\left(\xi_1, \xi_2\right) \right]\right)  \\
&\qquad =2\mathbb{E}\left[ Y_n Y_l \sum_{i_1 \in F_1 ^{n,l}} Y_{i_1} \right]\mathbb{E} \left[ \left( \sum_{i_2 \in F_1 ^{n,l}} Y_{i_2}  \right)  e\left(\xi_1, \xi_2\right) \right] \label{e_second der c21}\\
&\qquad  \quad + 2\mathbb{E}\left[ Y_l \sum_{ i_1 \in F_1 ^{n,l}} Y_{i_1} \right]  \mathbb{E}\left[ Y_n \left( \sum_{ i_2 \in F_1 ^{n,l}} Y_{i_2} \right)  e\left(\xi_1, \xi_2\right) \right] \label{e_second der c22}\\
& \qquad \quad +  2\mathbb{E}\left[ Y_l \sum_{ i_1 \in F_1 ^{n,l}} Y_{i_1} \right]   \mathbb{E} \left[ - \mathbb{E} \left[ Y_n  \sum_{i_2 \in F_1 ^{n,l}} Y_{i_2} \right] e\left(\xi_1, \xi_2\right) \right] \label{e_second der c23}\\
& \qquad \quad + 2\mathbb{E}\left[ Y_l \sum_{ i_1 \in F_1 ^{n,l}} Y_{i_1} \right]  \mathbb{E}\left[  \left( \sum_{i_2 \in F_1 ^{n,l}} Y_{i_2}  \right)\left( -i \xi_1 \right)  \mathbb{E}\left[ Y_n \sum_{i_3 \in F_1 ^{n,l}} Y_{i_3} \right] e\left(\xi_1, \xi_2\right)\right] \label{e_second der c24}\\
& \qquad \quad + 2\mathbb{E}\left[ Y_l \sum_{i_1 \in F_1 ^{n,l}} Y_{i_1} \right]  \mathbb{E}\left[  \left( \sum_{ i_2 \in F_1 ^{n,l}} Y_{i_2}  \right) \left( -i \xi_2 \right)  \mathbb{E}\left[ Y_n \sum_{j \in F_2 ^{n,l}} Y_j \right] e\left(\xi_1, \xi_2\right)\right]. \label{e_second der c25}
\end{align}
We will divide into two cases. The first case is when~$|n-l|>2L$ and the second case is when~$|n-l|\leq 2L$. In the case when~$|n-l|>2L$, we have at least~$K\left(K-2L\right)$ pairs of~$\left(n,l\right)$. Hence in order to get the right estimate we have to apply the decay of correlations (cf.~Lemma~\ref{l_exponential decay of correlations}).
The case when~$|n-l|\leq 2L$ is much easier because there are at most~$2KL$ pairs of~$\left(n,l\right)$. In this case, we just provide a rough estimate using Lemma~\ref{c_amgm}, Lemma~\ref{l_exponential decay of correlations} and Lemma~\ref{l_exchanging exponential terms}, which is still good enough for deducing Lemma~\ref{l_proof_of_2nd_der_(c)}.\\

\noindent \textbf{Case 1}. \ $|n-l|>2L$.

We first estimate the terms~\eqref{e_second der c12} to~\eqref{e_second der c25} except for the term~\eqref{e_second der c22}. We postpone the estimation of the terms~\eqref{e_second der c11} and~\eqref{e_second der c22}.\\

From now on, we set~$\left(\xi_1, \xi_2 \right) = \left( \frac{\tilde{\xi}}{\sqrt{K}} , \frac{\xi}{\sqrt{K}} \right)$ such that~$\frac{\left|\tilde{\xi}\right|}{\sqrt{K}}\leq \frac{\left|\xi\right|}{\sqrt{K}}\leq \delta \leq 1$. Let us begin with the estimation of~\eqref{e_second der c12}. We have by using Lemma~\ref{c_amgm} and Lemma~\ref{l_exponential decay of correlations} that
\begin{align}
\left|T_{\eqref{e_second der c12}}\right|
&\lesssim \left|\cov \left( Y_n, Y_l \right)\right|\mathbb{E} \left|  \left(  \sum_{ i \in F_1 ^{n,l}} Y_i   \right)^2 \right| \\
&\lesssim L^2 \exp\left(-CL\right) .
\end{align}
Let us turn to the estimation of~\eqref{e_second der c13} and~\eqref{e_second der c23}. We combine these terms to make a covariance and do the Taylor expansion with respect to the first variable. Then we get \\
\begin{align}
T_{\eqref{e_second der c13}}+T_{\eqref{e_second der c23}}
&=2  \mathbb{E} \left[ Y_n  \sum_{i_1 \in F_1 ^{n,l} } Y_{i_1}  \right] \cov\left(Y_l \sum_{ i_2 \in F_1 ^{n,l}} Y_{i_2} ,  e\left(\frac{\tilde{\xi}}{\sqrt{K}}, \frac{\xi}{\sqrt{K}}\right) \right)\\
& = 2  \mathbb{E} \left[ Y_n \sum_{i_1 \in F_1 ^{n,l} } Y_{i_1} \right] \cov\left(Y_l \sum_{ i_2 \in F_1 ^{n,l}} Y_{i_2} ,  e\left(0, \frac{\xi}{\sqrt{K}}\right) \right) \label{e_second der c13+c23-1} \\
& \quad + 2  \mathbb{E} \left[ Y_n  \sum_{i_1 \in F_1 ^{n,l} } Y_{i_1}  \right] \cov\left(Y_l \sum_{ i_2 \in F_1 ^{n,l}} Y_{i_2} ,  \sum_{ i_3 \in F_1 ^{n,l}} Y_{i_3}e\left(\frac{\hat {\xi}}{\sqrt{K}}, \frac{\xi}{\sqrt{K}}\right) \right) i \frac{\tilde{\xi}}{\sqrt{K}},\label{e_second der c13+c23-2}
\end{align}
where~$\frac{\hat {\xi}}{\sqrt{K}} $ is a real number between~$0$ and~$\frac{\tilde{\xi}}{\sqrt{K}}$. In particular,~$\left|\frac{\hat {\xi}}{\sqrt{K}}\right| \leq \left|\frac{\tilde{\xi}}{\sqrt{K}}\right|\leq \frac{\left|\xi\right|}{\sqrt{K}}$. Let us consider~\eqref{e_second der c13+c23-1}: We have by using Lemma~\ref{c_amgm} and Lemma~\ref{l_exponential decay of correlations} that
\begin{align}
\left|T_{\eqref{e_second der c13+c23-1}} \right| &\leq \left|2  \mathbb{E} \left[ Y_n  \sum_{i_1 \in F_1 ^{n,l} } Y_{i_1}  \right] \right| \left| \cov\left(Y_l \sum_{ i_2 : \ |i_2 -l |\leq L/2 } Y_{i_2} ,  e\left(0, \frac{\xi}{\sqrt{K}}\right) \right) \right| \\
& \quad +  \left|2  \mathbb{E} \left[ Y_n  \sum_{i_1 \in F_1 ^{n,l} } Y_{i_1}  \right] \right| \left| \mathbb{E}\left[Y_l \sum_{ \substack{i_2 \in F_1 ^{n,l} \\ |i_2 -l | > L/2 }} Y_{i_2}   e\left(0, \frac{\xi}{\sqrt{K}}\right) \right] \right| \\
& \quad +  \left|2  \mathbb{E} \left[ Y_n  \sum_{i_1 \in F_1 ^{n,l} } Y_{i_1}  \right] \right|\left|\mathbb{E} \left[Y_l \sum_{ \substack{i_2 \in F_1 ^{n,l} \\ |i_2 -l | > L/2 }} Y_{i_2}  \right] \mathbb{E} \left[ e\left(0, \frac{\xi}{\sqrt{K}}\right) \right] \right| \\
&  = \left|2  \mathbb{E} \left[ Y_n  \sum_{i_1 \in F_1 ^{n,l} } Y_{i_1}  \right] \right| \left| \cov\left(Y_l \sum_{ i_2 : \ |i_2 -l |\leq L/2 } Y_{i_2} ,  e\left(0, \frac{\xi}{\sqrt{K}}\right) \right) \right| \\
& \quad +  \left|2  \mathbb{E} \left[ Y_n  \sum_{i_1 \in F_1 ^{n,l} } Y_{i_1}  \right] \right| \left| \cov\left(Y_l ,  \sum_{ \substack{i_2 \in F_1 ^{n,l} \\ |i_2 -l | > L/2 }} Y_{i_2}   e\left(0, \frac{\xi}{\sqrt{K}}\right) \right) \right| \\
& \quad +  \left|2  \mathbb{E} \left[ Y_n  \sum_{i_1 \in F_1 ^{n,l} } Y_{i_1}  \right] \right|\left|\cov \left(Y_l ,\sum_{ \substack{i_2 \in F_1 ^{n,l} \\ |i_2 -l | > L/2 }} Y_{i_2}  \right) \right| \left| \mathbb{E} \left[ e\left(0, \frac{\xi}{\sqrt{K}}\right) \right] \right| \\
& \lesssim L^2 \left|\xi\right|\exp\left(-CL\right) + L^2 \left(1+ \left|\xi\right| \right) \exp\left(-CL\right) + L^2 \exp\left(-CL\right) \\
& \lesssim L^2 \left(1+ \left|\xi\right| \right) \exp\left(-CL\right).
\end{align}
Let us consider~\eqref{e_second der c13+c23-2}: By definition of covariances, it follows that
\begin{align}
T_{\eqref{e_second der c13+c23-2}} & = 2  \mathbb{E} \left[ Y_n  \sum_{i_1 \in F_1 ^{n,l} } Y_{i_1}  \right] \mathbb{E}\left[Y_l \sum_{ i_2 \in F_1 ^{n,l}} Y_{i_2}   \sum_{ i_3 \in F_1 ^{n,l}} Y_{i_3}e\left(\frac{\hat {\xi}}{\sqrt{K}}, \frac{\xi}{\sqrt{K}}\right) \right]\frac{\tilde {\xi}}{\sqrt{K}}  \label{e_second der c13+c23-21}\\
& \quad -  2  \mathbb{E} \left[ Y_n  \sum_{i_1 \in F_1 ^{n,l} } Y_{i_1}  \right]\mathbb{E}\left[Y_l \sum_{ i_2 \in F_1 ^{n,l}} Y_{i_2} \right] \mathbb{E} \left[  \sum_{ i_3 \in F_1 ^{n,l}} Y_{i_3}e\left(\frac{\hat {\xi}}{\sqrt{K}}, \frac{\xi}{\sqrt{K}}\right) \right]\frac{\tilde {\xi}}{\sqrt{K}} \label{e_second der c13+c23-22}.
\end{align}
We have by taking the conditional expectation with respect to~$\mathscr{G }_{n,l} := \sigma \left( X_k , k \in F_1 ^{n,l} \right)$ and applying Lemma~\ref{c_amgm} and Lemma~\ref{l_exchanging exponential terms},
\begin{align}
&\left|\mathbb{E}\left[Y_l \sum_{ i_2 \in F_1 ^{n,l}} Y_{i_2}   \sum_{ i_3 \in F_1 ^{n,l}} Y_{i_3}e\left(\frac{\hat {\xi}}{\sqrt{K}}, \frac{\xi}{\sqrt{K}}\right) \right] \right|\\
& \qquad = \left|\mathbb{E}\left[Y_l \sum_{ i_2 \in F_1 ^{n,l}} Y_{i_2}   \sum_{ i_3 \in F_1 ^{n,l}} Y_{i_3}e\left(\frac{\hat {\xi}}{\sqrt{K}},0 \right)\mathbb{E}\left[ e\left(0,\frac{ \xi}{\sqrt{K}}\right) \mid \mathscr{G}_{n, l} \right] \right] \right| \\
& \qquad \overset{Lemma~\ref{l_exchanging exponential terms}}{\lesssim} \mathbb{E}\left[\left| Y_l \sum_{ i_2 \in F_1 ^{n,l}} Y_{i_2}   \sum_{ i_3 \in F_1 ^{n,l}} Y_{i_3}e\left(\frac{\hat {\xi}}{\sqrt{K}},0 \right)\right| \right] \left(1+ \xi^2 \right) \exp\left(-C\xi^2 \right)  \\
& \qquad \overset{Lemma~\ref{c_amgm}}{\lesssim} L^2 \left(1+\xi^2 \right)\exp\left(-C\xi^2 \right),  \label{e_second der c13+c23-21 aux}
\end{align}
and similarly one gets
\begin{align}
\left|\mathbb{E} \left[  \sum_{ i_3 \in F_1 ^{n,l}} Y_{i_3}e\left(\frac{\hat {\xi}}{\sqrt{K}}, \frac{\xi}{\sqrt{K}}\right) \right]\right| \lesssim L\left(1+\xi^2 \right)\exp\left(-C\xi^2 \right). \label{e_second der c13+c23-22 aux}
\end{align}
Plugging~\eqref{e_second der c13+c23-21 aux} and~\eqref{e_second der c13+c23-22 aux} into~\eqref{e_second der c13+c23-21} and~\eqref{e_second der c13+c23-22} and applying Lemma~\ref{c_amgm} yield
\begin{align}
\left| T_{\eqref{e_second der c13+c23-2}} \right| \lesssim L^3\frac{\left|\xi\right|}{\sqrt{K}}\left(1+\xi^2 \right)\exp\left(-C\xi^2 \right).
\end{align}
Therefore we have
\begin{align}
\left|T_{\eqref{e_second der c13}}+T_{\eqref{e_second der c23}} \right| &\leq \left| T_{\eqref{e_second der c13+c23-1}} \right|+\left| T_{\eqref{e_second der c13+c23-2}} \right|\\
& \lesssim L^2 \left(1+ \left|\xi\right| \right)\exp\left(-CL\right) + L^3\frac{\left|\xi\right|}{\sqrt{K}}\left(1+\xi^2 \right)\exp\left(-C\xi^2 \right).
\end{align}
Let us turn to the estimation of~\eqref{e_second der c14}. This can also be estimated by taking the conditional expectation with respect to~$\mathscr{G }_{n,l}$ and applying Lemma~\ref{c_amgm} and Lemma~\ref{l_exchanging exponential terms}. Indeed,
\begin{align}
\left|T_{\eqref{e_second der c14}}\right| &= \left|\frac{\tilde {\xi}}{\sqrt{K}}\right| \left| \mathbb{E}\left[ Y_n \sum_{i_2 \in F_1 ^{n,l}} Y_{i_2} \right]\right|\left| \mathbb{E} \left[ Y_l \left( \sum_{ i_1 \in F_1 ^{n,l}} Y_{i_1}  \right)^2 e\left(\frac{\tilde {\xi}}{\sqrt{K}}, \frac{\xi}{\sqrt{K}}\right) \right] \right| \\
 & \lesssim  \frac{ \left| \xi \right|}{\sqrt{K}}L \left| \mathbb{E} \left[ Y_l \left( \sum_{ i_1 \in F_1 ^{n,l}} Y_{i_1}  \right)^2   e\left(\frac{\tilde {\xi}}{\sqrt{K}}, 0\right) \mathbb{E}\left[ e\left(0,\frac{\xi}{\sqrt{K}}\right) \mid \mathscr{G}_{n,l} \right] \right] \right|   \\
& \lesssim L^3 \frac{ \left| \xi \right|}{\sqrt{K}} \left(1+ \xi^2 \right) \exp\left(-C \xi^2\right).
\end{align}
Let us turn to the estimation of~\eqref{e_second der c15}. We have by using Lemma~\ref{c_amgm} and Lemma~\ref{l_exponential decay of correlations} that
\begin{align}
\left|T_{\eqref{e_second der c15}}\right| 
&= \frac{ \left|\xi\right|}{\sqrt{K}} \left|\mathbb{E}\left[ Y_n \sum_{j \in F_2 ^{n,l}} Y_j \right]\right|  \left|\mathbb{E} \left[ Y_l \left( \sum_{ i \in F_1 ^{n,l}} Y_i  \right)^2 e\left(\frac{\tilde {\xi}}{\sqrt{K}}, \frac{\xi}{\sqrt{K}} \right) \right]\right|\\
& \lesssim \frac{ \left|\xi \right|}{\sqrt{K}} \left|\cov\left( Y_n ,  \sum_{j \in F_2 ^{n,l}} Y_j \right)\right|\mathbb{E} \left|Y_l \left( \sum_{ i \in F_1 ^{n,l}} Y_i  \right)^2 \right|\\
& \lesssim L^2 \left|\xi\right| \exp\left(-CL\right),
\end{align}
Let us turn to the estimation of~\eqref{e_second der c21}. This also follows from dividing into two parts and applying Lemma~\ref{c_amgm} and Lemma~\ref{l_exponential decay of correlations} together. More precisely, we have
\begin{align}
\left|T_{\eqref{e_second der c21}}\right| &\leq \left| 2\sum_{i_1 : \ |i_1 -n| \leq L }\mathbb{E}\left[ Y_n Y_l  Y_{i_1} \right] \mathbb{E} \left[ \left( \sum_{i_2 \in F_1 ^{n,l}} Y_{i_2}  \right)  e\left(\frac{\tilde {\xi}}{\sqrt{K}}, \frac{\xi}{\sqrt{K}}\right) \right] \right|\label{e_dsigman c211}\\
& \quad + \left|2\sum_{i_1 : \ |i_1 -l| \leq L }\mathbb{E}\left[ Y_n Y_l  Y_{i_1} \right]\mathbb{E} \left[ \left( \sum_{i_2 \in F_1 ^{n,l}} Y_{i_2}  \right)  e\left(\frac{\tilde {\xi}}{\sqrt{K}}, \frac{\xi}{\sqrt{K}}\right) \right] \right| \label{e_dsigman c212}\\
&\leq 2\sum_{i_1 : \ |i_1 -n| \leq L }\left| \cov\left( Y_l, Y_n Y_{i_1} \right) \right| \mathbb{E} \left|  \sum_{i_2 \in F_1 ^{n,l}} Y_{i_2}    \right|\\
&\quad + 2\sum_{i_1 : \ |i_1 -l| \leq L }\left| \cov\left( Y_n, Y_l Y_{i_1} \right) \right| \mathbb{E} \left|  \sum_{i_2 \in F_1 ^{n,l}} Y_{i_2}    \right|\\
& \lesssim L^2\exp\left(-CL\right) .
\end{align}
Let us turn to the estimation of~\eqref{e_second der c24}. A similar argument as for~\eqref{e_second der c14} implies
\begin{align}
\left|T_{\eqref{e_second der c24}}\right| \lesssim L^3 \frac{ \left| \xi \right|}{\sqrt{K}} \left(1+\xi^2\right) \exp\left(-C \xi^2\right).
\end{align}
Let us turn to the estimation of~\eqref{e_second der c25}. A similar argument as for~\eqref{e_second der c15} yields
\begin{align}
\left|T_{\eqref{e_second der c25}}\right| \lesssim L^2 \left|\xi\right| \exp\left(-CL\right).
\end{align}
It remains to estimate~\eqref{e_second der c11} and~\eqref{e_second der c22}. Recalling the definition~\eqref{e_F1nl} of~$F_1 ^{n,l}$, a direct computation yields
\begin{align}
T_{\eqref{e_second der c11}} & = - \mathbb{E} \left[ Y_n Y_l  \left(  \sum_{ i_1 : \ |i_1-l|\leq L }  Y_{i_1} + \sum_{i_2 : \ |i_2 -n| \leq L} Y_{i_2}  \right)^2 e\left(\frac{\tilde {\xi}}{\sqrt{K}}, \frac{\xi}{\sqrt{K}}\right) \right] \\
& = - \mathbb{E} \left[ Y_n Y_l \left(  \sum_{ i_1 : \ |i_1-l|\leq L } Y_{i_1} \right)^2 e\left(\frac{\tilde {\xi}}{\sqrt{K}}, \frac{\xi}{\sqrt{K}}\right) \right]  
\label{e_dsigman c111} \\
&  \quad  - \mathbb{E} \left[ Y_n Y_l \left(  \sum_{ i_2 : \ |i_2 -l|\leq L } Y_{i_2}  \right)^2 e\left(\frac{\tilde {\xi}}{\sqrt{K}}, \frac{\xi}{\sqrt{K}}\right) \right]
  \label{e_dsigman c112}\\
&  \quad - 2\mathbb{E} \left[ Y_n Y_l \sum_{i_1 : \ |i_1 -n| \leq L} Y_{i_1} \sum_{ i_2 : \ |i_2 -l|\leq L } Y_{i_2}   e\left(\frac{\tilde {\xi}}{\sqrt{K}}, \frac{\xi}{\sqrt{K}} \right) \right], \label{e_dsigman c113} 
\end{align}
and
\begin{align}
T_{\eqref{e_second der c22}} & = 2\mathbb{E}\left[ Y_l \sum_{ i_1 : \ |i_1 -n|\leq L} Y_{i_1} \right] \mathbb{E}\left[ Y_n \left( \sum_{ i_2 \in F_1 ^{n,l}} Y_{i_2}  \right)  e\left(\frac{\tilde {\xi}}{\sqrt{K}}, \frac{\xi}{\sqrt{K}}\right) \right]  \\
&  \quad + 2\mathbb{E}\left[ Y_l \sum_{ i_1 : \ |i_1 -l|\leq L} Y_{i_1} \right] \mathbb{E}\left[ Y_n \left( \sum_{ i_2 \in F_1 ^{n,l}} Y_{i_2}  \right)  e\left(\frac{\tilde {\xi}}{\sqrt{K}}, \frac{\xi}{\sqrt{K}}\right) \right] \\
& = 2\mathbb{E}\left[ Y_l \sum_{ i_1 : \ |i_1 -n|\leq L} Y_{i_1} \right] \mathbb{E}\left[ Y_n \left( \sum_{ i_2 \in F_1 ^{n,l}} Y_{i_2}  \right)  e\left(\frac{\tilde {\xi}}{\sqrt{K}}, \frac{\xi}{\sqrt{K}}\right) \right] \label{e_dsigman c221}\\
& \quad + 2\mathbb{E}\left[ Y_l \sum_{ i_1 : \ |i_1 -l|\leq L} Y_{i_1} \right] \mathbb{E}\left[ Y_n  \sum_{ i_2 : \ |i_2 - l | \leq L} Y_{i_2}    e\left(\frac{\tilde {\xi}}{\sqrt{K}}, \frac{\xi}{\sqrt{K}}\right) \right] \label{e_dsigman c222}\\
& \quad + 2\mathbb{E}\left[ Y_l \sum_{ i_1 : \ |i_1 -l|\leq L} Y_{i_1} \right] \mathbb{E}\left[ Y_n  \sum_{ i_2 : \ |i_2 - n | \leq L} Y_{i_2}    e\left(\frac{\tilde {\xi}}{\sqrt{K}}, \frac{\xi}{\sqrt{K}}\right) \right] . \label{e_dsigman c223}
\end{align}
The terms~\eqref{e_dsigman c111},~\eqref{e_dsigman c112},~\eqref{e_dsigman c221} and~\eqref{e_dsigman c222} can be estimated using the arguments from above. Let us turn to the estimation of~\eqref{e_dsigman c111}. Taylor expansion with respect to the first variable yields
\begin{align}
T_{\eqref{e_dsigman c111}}
& = - \mathbb{E} \left[ Y_n Y_l \left(  \sum_{ i_1 : \ |i_1 -l|\leq L } Y_{i_1} \right)^2 e\left(0, \frac{\xi}{\sqrt{K}}\right) \right]
 \label{c111 taylor 0}\\
& \quad -  \mathbb{E} \left[ Y_n Y_l \left(  \sum_{ i_1 : \ |i_1 -l|\leq L } Y_{i_1}  \right)^2 \left(  \sum_{ i_2 \in F_1 } Y_{i_2}  \right)e\left(\frac{\hat {\xi}}{\sqrt{K}}, \frac{\xi}{\sqrt{K}}\right) \right]i\frac{\tilde {\xi}}{\sqrt{K}}.
\label{c111 taylor 1}
\end{align}
Let us consider~\eqref{c111 taylor 0}: Note that for~$i_1$ with~$\left| i_1 -l \right| \leq L$, we have~$\left|i_1 -n \right| \geq \left|n-l\right| - \left| i_1 -l \right| > 2L - L = L$. Thus applying Lemma~\ref{c_amgm} and Lemma~\ref{l_exponential decay of correlations} yields
\begin{align}
\left|T_{\eqref{c111 taylor 0}} \right| &   = \left| \cov\left( Y_n, Y_l \left(  \sum_{ i_1 : \ |i_1 -l|\leq L } Y_{i_1}  \right)^2 e\left(0, \frac{\xi}{\sqrt{K}}\right)   \right)
\right| \lesssim L^2 \left(1+ \left|\xi\right|\right) \exp\left(-CL\right).
\end{align}
Let us consider~\eqref{c111 taylor 1}: A similar argument as for~\eqref{e_second der c14} yields
\begin{align}
\left|T_{\eqref{c111 taylor 1}}  \right| 
 \lesssim L^3  \frac{\left|\xi\right|}{\sqrt{K}}\left(1+ \xi^2 \right) \exp\left(-C\xi^2\right),
\end{align}
and thus
\begin{align}
\left|T_{\eqref{e_dsigman c111}} \right| & = \left|T_{\eqref{c111 taylor 0}}+T_{\eqref{c111 taylor 1}}\right| \\
& \leq \left|T_{\eqref{c111 taylor 0}}\right| +\left|T_{\eqref{c111 taylor 1}}\right|\\
&\lesssim L^2 \left(1+ \left|\xi\right|\right)\exp\left(-CL\right) + L^3  \frac{\left|\xi\right|}{\sqrt{K}} \left(1+ \xi^2 \right) \exp\left(-C\xi^2\right).
\end{align}
Let us turn to the estimation of~\eqref{e_dsigman c112}. A similar argument as for~\eqref{e_dsigman c111} yields
\begin{align}
\left|T_{\eqref{e_dsigman c112}} \right| \lesssim L^2 \left(1+ \left|\xi\right|\right)\exp\left(-CL\right) + L^3  \frac{\left|\xi\right|}{\sqrt{K}} \left(1+ \xi^2 \right) \exp\left(-C\xi^2\right).
\end{align}
Let us turn to the estimation of~\eqref{e_dsigman c221}. It follows from Lemma~\ref{c_amgm} and Lemma~\ref{l_exponential decay of correlations} that
\begin{align}
\left|T_{\eqref{e_dsigman c221}} \right| & = \left| 2\cov\left( Y_l ,  \sum_{ i_1 : \ |i_1 -n|\leq L} Y_{i_1} \right)\right| \left| \mathbb{E}\left[ Y_n \left( \sum_{ i_2 \in F_1 ^{n,l}} Y_{i_2}  \right)  e\left(\frac{\tilde {\xi}}{\sqrt{K}}, \frac{\xi}{\sqrt{K}}\right) \right] \right| \\
& \lesssim L^{3/2} \exp\left(-CL\right) \lesssim L^2  \exp\left(-CL\right).
\end{align}
Let us turn to the estimation of~\eqref{e_dsigman c222}. A similar argument as for~\eqref{e_dsigman c111} yields
\begin{align}
\left|T_{\eqref{e_dsigman c222}} \right| &\lesssim L^2 \left(1+ \left|\xi\right|\right) \exp\left(-CL\right)+ L^3 \frac{\left|\xi\right|}{\sqrt{K}}\left(1+ \xi^2 \right) \exp\left(-C \xi^2\right).
\end{align}
Let us turn to the estimation of~\eqref{e_dsigman c113} and~\eqref{e_dsigman c223}. Combining these terms and applying Taylor expansion with respect to the first variable yield
\begin{align}
&T_{\eqref{e_dsigman c113}}+ T_{\eqref{e_dsigman c223}}\\
&\qquad = - 2\cov\left( Y_l \sum_{ i_1 : \ |i_1 -l|\leq L} Y_{i_1} ,  Y_n  \sum_{ i_2 : \ |i_2 -n|\leq L }Y_{i_2}   e\left(\frac{\tilde {\xi}}{\sqrt{K}}, \frac{\xi}{\sqrt{K}}\right) \right) \\
& \qquad   = -2\cov\left( Y_l \sum_{ i_1 : \ |i_1 -l|\leq L}Y_{i_1} , Y_n  \sum_{ i_2 : \ |i_2 -n|\leq L } Y_{i_2}   e\left(0, \frac{\xi}{\sqrt{K}}\right) \right) \label{e_dsigman c113+c223-1} \\
& \qquad \quad - 2\cov\left( Y_l \sum_{ i_1 : \ |i_1 -l|\leq L} Y_{i_1}, Y_n \sum_{ i_2 : \ |i_2 -n|\leq L } Y_{i_2} \sum_{ i_3 \in F_1 ^{n,l} } Y_{i_3}   e\left(\frac{\hat {\xi}}{\sqrt{K}}, \frac{\xi}{\sqrt{K}}\right) \right) i\frac{\tilde {\xi}}{\sqrt{K}}.\label{e_dsigman c113+c223-2} 
\end{align}
Let us consider~\eqref{e_dsigman c113+c223-1}: Dividing into two parts, we get
\begin{align}
T_{\eqref{e_dsigman c113+c223-1}}&=  -2\cov\left( Y_l \sum_{ i_1 : \ |i_1 -l|\leq L/2} Y_{i_1} , Y_n \left( \sum_{ i_2 : \ |i_2 -n|\leq L } Y_{i_2}  \right)  e\left(0, \frac{\xi}{\sqrt{K}}\right) \right) \label{e_dsigman c113+c223-11}\\
&  \quad - 2\cov\left(Y_l \sum_{ i_1 : \ L/2 \leq |i_1 -l|\leq L} Y_{i_1},  Y_n \left( \sum_{ i_2 : \ |i_2 -n|\leq L } Y_{i_2}  \right)  e\left(0, \frac{\xi}{\sqrt{K}}\right) \right). \label{e_dsigman c113+c223-12}
\end{align}
Note that for~$\left(i_1, i_2\right)$ with~$|i_1 -l|\leq \frac{L}{2}$ and~$|i_2 - n | \leq L$, it holds that
\begin{align}
\left|i_1 -i_2 \right| \geq \left|l-n\right| - \left|i_1-l\right|-\left|i_2-n\right| \geq 2L-\frac{L}{2}-L \geq \frac{L}{2} .
\end{align}
Therefore Lemma~\ref{c_amgm} and Lemma~\ref{l_exponential decay of correlations} implies
\begin{align}
\left| T_{\eqref{e_dsigman c113+c223-11}} \right| &\lesssim  L^2 \left(1+ \left|\xi\right|\right) \exp\left(-CL\right) ,
\end{align}
and also by definition of covariances, it holds that
\newpage
\begin{align}
\left|T_{\eqref{e_dsigman c113+c223-12}}\right|
&\leq \left|2\mathbb{E}\left[ Y_l \left( \sum_{ i_1 : \ L/2 \leq |i_1 -l|\leq L} Y_{i_1}\right)Y_n \left( \sum_{ i_2 : \ |i_2 -n|\leq L } Y_{i_2}  \right)  e\left(0, \frac{\xi}{\sqrt{K}}\right) \right] \right|\\
& \quad + \left|2\mathbb{E}\left[ Y_l \sum_{ i_1 : \ L/2 \leq |i_1 -l|\leq L} Y_{i_1}\right] \right|\left| \mathbb{E}\left[Y_n \left( \sum_{ i_2 : \ |i_2 -n|\leq L } Y_{i_2} \right)  e\left(0, \frac{\xi}{\sqrt{K}}\right) \right]\right| \\
& = \left| 2\cov\left( Y_l ,  \left(\sum_{ i_1 : \ L/2 \leq |i_1 -l|\leq L} Y_{i_1}\right) Y_n  \left( \sum_{ i_2 : \ |i_2 -n|\leq L } Y_{i_2}\right)  e\left(0, \frac{\xi}{\sqrt{K}}\right) \right) \right| \\
& \quad + \left|2\cov\left( Y_l , \sum_{ i_1 : \ L/2 \leq |i_1 -l|\leq L} Y_{i_1}\right) \right|\mathbb{E}\left|Y_n \left( \sum_{ i_2 : \ |i_2 -n|\leq L } Y_{i_2}  \right)  \right| \\
& \lesssim L^2 \left(1+ \left|\xi \right| \right) \exp\left(-CL\right) + L^{3/2} \exp\left(-CL\right) \lesssim L^2 \left(1+ \left|\xi\right| \right) \exp\left(-CL\right).
\end{align}
Thus we have
\begin{align}
\left|T_{\eqref{e_dsigman c113+c223-1}}\right| \leq \left| T_{\eqref{e_dsigman c113+c223-11}}\right| + \left| T_{\eqref{e_dsigman c113+c223-12}}\right| \lesssim L^2 \left(1+ \left|\xi\right| \right) \exp\left(-CL\right).
\end{align}
Let us consider~\eqref{e_dsigman c113+c223-2}: A similar argument as for~\eqref{e_second der c13+c23-2} yields
\begin{align}
\left|T_{\eqref{e_dsigman c113+c223-2}}\right| \lesssim L^3 \frac{ \left|\xi\right|}{\sqrt{K}}\left(1+ \xi^2 \right) \exp\left(-C\xi^2\right).
\end{align}
Therefore we conclude that
\begin{align}
\left| T_{\eqref{e_dsigman c113}}+ T_{\eqref{e_dsigman c223}} \right| &\leq  \left|T_{\eqref{e_dsigman c113+c223-1}}\right|+\left|T_{\eqref{e_dsigman c113+c223-2}}\right| \\
&\lesssim L^2 \left(1+ \left|\xi\right| \right) \exp\left(-CL\right) + L^3 \frac{ \left|\xi\right|}{\sqrt{K}}\left(1+ \xi^2 \right) \exp\left(-C\xi^2\right).
\end{align}
To sum up the estimation of~\eqref{e_second der c11} and~\eqref{e_second der c22}, we have
\begin{align}
\left|T_{\eqref{e_second der c11}} + T_{\eqref{e_second der c22}}\right| &=\left|\left(T_{\eqref{e_dsigman c111}}+T_{\eqref{e_dsigman c112}}+T_{\eqref{e_dsigman c113}}\right) + \left(T_{\eqref{e_dsigman c221}}+T_{\eqref{e_dsigman c222}}+ T_{\eqref{e_dsigman c223}}\right) \right| 
\\ & \leq \left|T_{\eqref{e_dsigman c111}}\right|+\left|T_{\eqref{e_dsigman c112}}\right|+ \left|T_{\eqref{e_dsigman c221}}\right|+ \left|T_{\eqref{e_dsigman c222}}\right| + \left| T_{\eqref{e_dsigman c113}} +  T_{\eqref{e_dsigman c223}} \right|  \\
&\lesssim  L^2 \left(1+ \left|\xi\right| \right) \exp\left(-CL\right)  + L^3 \frac{ \left|\xi\right|}{\sqrt{K}} \left( 1+ \xi^2 \right) \exp\left(-C \xi^2 \right) \\
&\lesssim L^2 \left(1+ \xi^2 \right) \exp\left(-CL\right)  + L^3 \frac{ \left|\xi\right|}{\sqrt{K}} \left( 1+ \xi^2 \right) \exp\left(-C \xi^2 \right) .
\end{align}
Overall we have proven the desired estimates that for~$|n-l|>2L$,
\begin{align}
&\left|T_{\eqref{e_dsigman c1}} + T_{\eqref{e_dsigman c2}}\right| \\
&  \leq \left|T_{\eqref{e_second der c11}}+T_{\eqref{e_second der c22}} \right|+\left|T_{\eqref{e_second der c12}} \right|+ \left|T_{\eqref{e_second der c13}} + T_{\eqref{e_second der c23}}\right| + \left| T_{\eqref{e_second der c14}}\right| + \left| T_{\eqref{e_second der c15}}\right|+ \left| T_{\eqref{e_second der c21}}\right|+ \left| T_{\eqref{e_second der c24}}\right| +\left| T_{\eqref{e_second der c25}}\right|   
\\
& \lesssim L^2\left(1+ \xi^2 \right) \exp\left(-CL\right) + L^3 \frac{\left|\xi\right|}{\sqrt{K}}\left(1+ \xi^2\right) \exp\left(-C \xi^2\right).
\end{align}
\medskip
\newpage
\noindent \textbf{Case 2}. \ $|n-l|\leq 2L$.

Let us set~$\left(\xi_1, \xi_2 \right) = \left( \frac{\tilde {\xi}}{\sqrt{K}}, \frac{\xi}{\sqrt{K}}\right)$ with~$\left|\frac{\tilde {\xi}}{\sqrt{K}}\right|\leq \left|\frac{\xi}{\sqrt{K}}\right| \leq \delta \leq 1$. Except for~\eqref{e_second der c15} and~\eqref{e_second der c25}, we may apply~\eqref{e_lemma exchanging 2nd der} in Lemma~\ref{l_exchanging exponential terms} followed by~\eqref{e_amgm}. For example, we have
\begin{align}
\left|T_{\eqref{e_second der c11}}\right| &  = \left| \mathbb{E} \left[ Y_n Y_l \left(  \sum_{ i \in F_1 ^{n,l}} Y_i  \right)^2 e\left(\frac{\tilde {\xi}}{\sqrt{K}}, 0\right) \mathbb{E} \left[ e\left(0, \frac{\xi}{\sqrt{K}} \right) \mid \mathscr{G}_{n,l} \right] \right]\right|\\
&  \lesssim \mathbb{E} \left|Y_n Y_l \left(  \sum_{ i \in F_1 ^{n,l}} Y_i  \right)^2 \right| \left(1+ \frac{\left|\xi\right|^3}{\sqrt{K}}\right) \exp\left(-C \xi^2\right)  \\
& \lesssim  L^2\left(1+ \xi^2 \right) \exp\left(-C \xi^2\right).
\end{align}
Similar computations yield
\begin{align}
\left|T_{\eqref{e_second der c12}}\right|, \left|T_{\eqref{e_second der c13}}\right|, \left|T_{\eqref{e_second der c21}} \right| , \left|T_{\eqref{e_second der c22}} \right|, \left|T_{\eqref{e_second der c23}}\right| \lesssim L^2 \left(1 + \xi^2 \right) \exp\left(-C\xi^2 \right),
\end{align}
\begin{align}
\left|T_{\eqref{e_second der c14}}\right|, \left|T_{\eqref{e_second der c24}}\right|\lesssim  L^3 \frac{\left|\xi\right|}{\sqrt{K}}\left(1 + \xi^2 \right) \exp\left(-C\xi^2 \right).
\end{align}
Let us turn to the estimation of~\eqref{e_second der c15} and~\eqref{e_second der c25}: Because those terms involve a sum over~$F_2 ^{n,l}$, one needs to apply the decay of correlations (cf.~Lemma~\ref{l_exponential decay of correlations}). It follows that (see also the computation for the case~$|n-l|>2L$)
\begin{align}
\left|T_{\eqref{e_second der c15}}\right|, \left|T_{\eqref{e_second der c25}}\right| \lesssim L^2 \left|\xi\right| \exp\left(-CL\right).
\end{align}
Therefore we have for~$|n-l|\leq 2L$,
\begin{align}
&\left|T_{\eqref{e_dsigman c1}} + T_{\eqref{e_dsigman c2}}\right| \\
& \leq \left|T_{\eqref{e_second der c11}}\right| +\left|T_{\eqref{e_second der c12}}\right|+\left|T_{\eqref{e_second der c13}}\right| +\left|T_{\eqref{e_second der c14}}\right| +\left|T_{\eqref{e_second der c15}} \right|+\left|T_{\eqref{e_second der c21}}\right| +\left|T_{\eqref{e_second der c22}}\right|+\left|T_{\eqref{e_second der c23}}\right|+\left|T_{\eqref{e_second der c24}}\right|+\left|T_{\eqref{e_second der c25}} \right|  \\
& \lesssim L^2 \left|\xi\right| \exp\left(-CL\right) + L^2\left(1 + \xi^2 \right) \exp\left(-C\xi^2 \right) +L^3 \frac{\left|\xi\right|}{\sqrt{K}}\left(1 + \xi^2\right) \exp\left(-C\xi^2 \right) \\
& \lesssim L^2 \left(1+ \xi^2 \right) \exp\left(-CL\right) + L^3 \left(1 + \xi^2 \right) \exp\left(-C\xi^2 \right).
\end{align}
\medskip

To conclude, we sum up all the bounds we have proven so far. That is, for~$\left|n-l\right|> 2L$, we have
\begin{align}
&\left| \frac{d^2}{d \xi_1^2}  \frac{d}{d \sigma_n }\frac{d}{d \sigma_l } G_{n,l} \left(\tilde{ \frac{\xi}{\sqrt{K}}}, \frac{\xi}{\sqrt{K}}\right) \right|\\
&\qquad \leq \left|T_{\eqref{e_dsigman c1}} +T_{\eqref{e_dsigman c2}} \right|+\left|T_{\eqref{e_dsigman c3}} \right|+\left|T_{\eqref{e_dsigman c4}} \right|\\
&\qquad  \lesssim L^3 \left(1+ \xi^2 \right) \exp\left(-CL\right)  + L^4 \frac{ \left|\xi\right|}{\sqrt{K}} \left(1+ \xi^2 \right) \exp\left(-C \xi^2 \right)
\end{align}
and for~$\left|n-l\right|\leq 2L$,
\begin{align}
&\left| \frac{d^2}{d \xi_1^2}  \frac{d}{d \sigma_n }\frac{d}{d \sigma_l } G_{n,l} \left(\tilde{ \frac{\xi}{\sqrt{K}}}, \frac{\xi}{\sqrt{K}}\right) \right|\\
& \qquad \leq \left|T_{\eqref{e_dsigman c1}} \right|+\left|T_{\eqref{e_dsigman c2}} \right|+\left|T_{\eqref{e_dsigman c3}} \right|+\left|T_{\eqref{e_dsigman c4}} \right|\\
& \qquad  \lesssim L^3 \left(1+\xi^2 \right)\exp\left(-CL\right) + L^3 \left(1+ \xi^2 \right) \exp\left(-C\xi^2 \right)\\
& \qquad \quad + L^4 \frac{\left|\xi\right|}{\sqrt{K}} \left(1+ \xi^2 \right)\exp\left(-C \xi^2 \right) \\
& \qquad \lesssim L^3 \left(1+\xi^2 \right)\exp\left(-CL\right) + L^4 \left(1+ \xi^2 \right) \exp\left(-C\xi^2 \right). 
\end{align}
This finishes the proof of Lemma~\ref{l_proof_of_2nd_der_(c)}. \qed
\medskip

\begin{remark} The proof of Lemma~\ref{l_proof_of_2nd_der_(b)} is almost the same as the proof of Lemma~\ref{l_proof_of_2nd_der_(c)}. Carrying out a Taylor expansion with respect to the second variable in~$G_{n,l}$ yields
\begin{align}
&\frac{d}{d \xi_1} \frac{d}{d \sigma_n}\frac{d}{d \sigma_l} G_{n,l} (0, \xi_2) \xi_1  \\
&\qquad  = \underbrace{\frac{d}{d \xi_1}\frac{d}{d \sigma_n} \frac{d}{d \sigma_l} G_{n,l} (0, 0)}_{=0} \xi_1  + \frac{d}{d \xi_2}\frac{d}{d \xi_1} \frac{d}{d \sigma_n}\frac{d}{d \sigma_l} G_{n,l} (0, \tilde \xi_2) \xi_1 \xi_2. \label{second derivative 2nd term taylor again}
\end{align}
The first term in~\eqref{second derivative 2nd term taylor again} vanishes because
\begin{align}
\frac{d}{d \xi_1}\frac{d}{d \sigma_n} \frac{d}{d \sigma_l} G_{n,l} (0, 0) = \frac{d}{d \sigma_n} \frac{d}{d \sigma_l} \frac{d}{d \xi_1}G_{n,l} (0, 0) = \frac{d}{d \sigma_n} \frac{d}{d \sigma_l} \mathbb{E} \left[ i \sum_{i \in F_1 ^{n,l} }Y_i \right] = 0.
\end{align}
Then similar arguments applied to the second term in~\eqref{second derivative 2nd term taylor again} gives the desired bound.
\end{remark}

\section*{Acknowledgment}

The authors want to thank Felix Otto and H.T.~Yau for bringing this problem to their attention. The authors are also thankful to many people discussing the problem and helping to improve the preprint. Among them are Tim Austin, Frank Barthe, Marek Biskup, Pietro Caputo, Jean-Dominique Deuschel, Max Fathi, Andrew Krieger, Michel Ledoux, Sangchul Lee, Thomas Liggett, Felix Otto, Daniel Ueltschi, and Tianqi Wu. The authors want to thank Marek Biskup, UCLA and KFAS for financial support. Last but not least, the authors thank the anonymous referee for the detailed report that helped to further improve the article.

\bibliographystyle{spmpsci}      
\bibliography{bib}


%



\end{document}